\documentclass[pdflatex,sn-mathphys-num]{sn-jnl}

\usepackage{graphicx}%
\usepackage{multirow}%
\usepackage{amsmath,amssymb,amsfonts}%
\usepackage{amsthm}%
\usepackage{mathrsfs}%
\usepackage[title]{appendix}%
\usepackage{xcolor}%
\usepackage{textcomp}%
\usepackage{manyfoot}%
\usepackage{booktabs}%
\usepackage{algorithm}%
\usepackage{algorithmicx}%
\usepackage{algpseudocode}%
\usepackage{listings}%

\usepackage{bbm}
\usepackage{amsfonts}
\usepackage{multirow}
\usepackage{pgfplots}
\usepackage{algorithm}
\usepackage{algpseudocode}
\usepackage{enumitem}
\usepackage[capitalize]{cleveref}
\usepackage{aliascnt}
 \usepackage{tikz}
\usepackage{pgfplots}
\usepackage{pgfplotstable}
\usepackage{subcaption}
\usepackage{longtable}
\usepackage{colortbl}
\usepackage{pdflscape}
\usepackage{adjustbox}
\usepackage{mathrsfs}
\usepackage{breakurl}
\usetikzlibrary{patterns}
\usepgfplotslibrary{groupplots}
\usepgfplotslibrary{colormaps, colorbrewer}

\pgfplotsset{
    colormap={magma}{
        rgb255=(15,15,40)    
        rgb255=(130,70,140)  
        rgb255=(190,80,100)  
        rgb255=(230,160,90)  
        rgb255=(240,200,140) 
    }
}

\pgfplotsset{compat=1.18}

\newtheorem{example}{Example}
\newtheorem{assumption}{Assumption}
\DeclareMathAlphabet\mathbfcal{OMS}{cmsy}{b}{n}
\newcommand{\vecTT}{\mathop{\mathrm{vec}}\nolimits_{\mathrm{TT}}}
\newcommand{\diagTT}{\mathop{\mathrm{diag}}\nolimits_{\mathrm{TT}}}
\newcommand{\rankTT}{\mathop{\mathrm{rank}}\nolimits_{\mathrm{TT}}}

\theoremstyle{thmstyleone}%

\theoremstyle{thmstyletwo}%

\theoremstyle{thmstylethree}%
\newtheorem{theorem}{Theorem}
\crefname{theorem}{theorem}{theorems}
\Crefname{theorem}{Theorem}{Theorems}
\newaliascnt{proposition}{theorem}
\newtheorem{proposition}[proposition]{Proposition}%
\aliascntresetthe{proposition}
\crefname{proposition}{proposition}{propositions}
\Crefname{proposition}{Proposition}{Propositions}
\newtheorem{definition}{Definition}%
\newtheorem{remark}{Remark}%
\newtheorem{lemma}{Lemma}%
\crefname{lemma}{lemma}{lemmas}
\Crefname{lemma}{Lemma}{Lemmas}
\newtheorem{corollary}{Corollary}%

\begin{document}

\title{An Inexact Tensor-Train Primal-Dual Interior-Point Method for Semidefinite Programs}


\author*[1]{\fnm{Frederik} \sur{Kelbel}}\email{f.kelbel@ic.ac.uk}

\author[3]{\fnm{Sergey} \sur{Dolgov}}\email{s.dolgov@bath.ac.uk}
\equalcont{These authors contributed equally to this work.}

\author[2]{\fnm{Dante} \sur{Kalise}}\email{d.kalise-balza@imperial.ac.uk}
\equalcont{These authors contributed equally to this work.}

\author[1]{\fnm{Alessandra} \sur{Russo}}\email{a.russo@imperial.ac.uk}
\equalcont{These authors contributed equally to this work.}

\affil*[1]{\orgdiv{Department of Computing}, \orgname{Imperial College London}, \country{UK}}

\affil[2]{\orgdiv{Department of Mathematics}, \orgname{Imperial College London}, \country{UK}}

\affil[3]{\orgdiv{Department of Mathematics}, \orgname{University of Bath}, \country{UK}}


\abstract{In this work, we introduce an interior-point method that employs tensor decompositions to efficiently represent and manipulate the variables and constraints of semidefinite programs, targeting problems where the solutions may not be low-rank but admit low-tensor-train rank approximations. Our method leverages a primal-dual infeasible interior-point framework { with global and local convergence estimates}. In experiments on Maximum Cut, Maximum Stable Set, and Correlation Clustering, the tensor-train interior-point method handles problems up to size $2^{12}$ {with certified primal-dual accuracy, reported gaps of approximate order $10^{-4}$, and remains the only primal-dual SDP solver reported at the largest dimensions across all three benchmarks}. Moreover, numerical evidence indicates that tensor-train ranks of the iterates {typically remain structured and moderate along the interior-point trajectory, with instance-dependent growth for harder high-rank cases, supporting} the scalability of the approach. Tensor-train interior point methods offer a promising avenue for problems that lack traditional sparsity or low-rank structure, exploiting tensor-train structures instead.
}




\maketitle


\section{Introduction}

Semidefinite programs (SDPs) constitute a class of convex optimisation problems that form the foundation of many of the most powerful relaxations in combinatorial optimisation. Notably, SDPs can form the tightest convex relaxations of linear optimisation problems of any algebraic varieties \cite{Bleckherman2012, Laserre2001}. Under standard assumptions, these programs satisfy strong duality and admit polynomial-time solvability, rendering them a highly effective framework for addressing a broad class of NP-hard problems. Despite their theoretical appeal, the practical applicability of SDP relaxations is often hindered by the rapid growth in problem dimensionality associated with increasing problem size or tightening relaxations. This manifests in substantial computational and memory demands. Interior-Point Methods (IPMs), which are the predominant algorithms employed for solving SDPs, require manipulation of large, dense matrices. This limits their scalability in high-dimensional settings. To mitigate these challenges, considerable research has focused on first-order methods and the exploitation of low-rank structures inherent to the solutions of certain SDPs \cite{Arora2005, Renegar2019, Henrion2009, Yurtsever2017, Yurtsever2019, Yurtsever2021}. Maintaining a low-rank representation throughout the optimisation process has been shown to significantly reduce both memory requirements and computational overhead. This approach is particularly effective for problems characterised by relatively simple polytope structures, although it is well-documented that the solution rank tends to increase with the number of constraints imposed \cite{Pataki1998}. Moreover, first-order methods suffer from sublinear convergence rates. The Alternating Direction Method of Multipliers (ADMM) \cite{Wen2010}, in particular, is a scalable first-order method that tackles large-scale SDPs by decomposing them into simpler subproblems. While this approach substantially reduces per-iteration cost and memory usage, it often suffers from slow convergence, particularly near the solution, and generally yields lower-accuracy solutions compared to Interior-Point Methods. Alternative strategies, such as low-rank Burer–Monteiro factorisations \cite{Burer2003} combined with Riemannian Trust Region algorithms \cite{Boumal2016}, offer promising avenues but are applicable only to an even more restricted subset of SDPs.

The present work proposes a different model for approximating SDP solutions based on the Tensor-Train (TT) decomposition \cite{Oseledets2011}. This decomposition provides a compact representation of dense matrices by first reshaping them into higher-order tensors and then approximating them through a sequence of low-rank tensor contractions. It replaces the conventional notion of low matrix rank with the concept of low tensor-train rank. IPMs are of particular interest due to their superlinear convergence rates and high accuracy. This motivates our focus on tensor-train-based IPMs (TT IPM). The applicability of TT decompositions within IPMs has been previously established in the context of large-scale quadratic programming \cite{Garreis2017, Miculescu2018}.

This paper develops a tensor-train-based interior-point method grounded in three central hypotheses. First, we posit that there exist semidefinite programs whose \emph{solutions} are inherently of \emph{low tensor-train rank}, even when they do not admit low-rank solutions. Second, we show that tensor-train interior-point methods produce \emph{interior paths} of \emph{low tensor-train rank}, leading to substantially reduced memory requirements. Finally, we {prove a global convergence estimate, with local convergence under
a bounded tensor-train residual similarly to standard inexact Newton's forcing conditions, 
and observe that the tensor-train format does not preclude fast interior-point convergence on the reported instances}.

We investigate these hypotheses through three canonical problems. The Maximum Cut (MaxCut) problem is selected due to its known low-rank solutions, enabling assessment of the extension to low TT rank. The Maximum Stable Set problem \cite{Lovasz1979} is considered owing to its typically high-rank solutions, which pose challenges to traditional low-rank approaches. Finally, we examine the Correlation Clustering Problem's \cite{Bansal2004} relaxation introduced by \cite{Chaitanya2004}, which combines structures from MaxCut and Maximum Stable Set.

In our experiments, the TT based interior-point method solved large-scale problems {with good primal-dual accuracy, while runtime and memory remained favourable for many structured low- and moderate-rank instances}. Classical solvers such as SDPA \cite{Makoto2003} and SCS \cite{SCSSolver2016} perform well on small instances but struggle as problem size or dimensionality increases. The TT IPM, on the other hand, {reaches the largest reported dimensions in all three benchmark classes, including instances where SDPA and SCS have no completed run, and remains competitive on several higher-rank instances}, showing strong potential for semidefinite programs that are challenging for traditional methods.

The rest of the paper is structured as follows. Section~\ref{sec:preliminaries} reviews the foundations of semidefinite programming together with the tensor-train decomposition framework. Building on this background, Section~\ref{sec:ttipm} develops the proposed TT based interior-point method in detail, including the representation of decision variables, the formulation and solution of the KKT system within TT arithmetic, {the inexact-Newton and convergence analysis,} and a range of practical considerations essential for numerical stability and efficiency. Section~\ref{sec:conceptualisation} formalises the three benchmark problem classes, namely, Maximum Cut, Maximum Stable Set, and Correlation Clustering, which are used to assess the method. The computational study in Section~\ref{sec:comp_study} outlines the experimental setup and problem instance generation, followed by an extensive performance evaluation, with particular attention to the interplay between problem rank, TT rank growth, and solver efficiency. Section~\ref{sec:discussion} offers an interpretation of these findings, summarises the main contributions, and identifies promising avenues for further research.

\section{Preliminaries}\label{sec:preliminaries}

In this section, we provide an overview of the fundamental concepts underlying semidefinite programming and tensor-train representations.

\subsection{Notation and Conventions}

For any matrices $A, B \in \mathbb{R}^{m \times n}$, the Frobenius inner product is defined as \[ \langle A, B \rangle := A \bullet B = \text{tr}(A^\top B). \] The symbol $\preceq$ denotes the L\"owner partial order defined by the cone of positive semidefinite matrices; that is, $A \preceq B$ (or $A \prec B$ for strict inequality). To every linear operator $A : \mathbb{R}^{n \times n} \to \mathbb{R}^m$, there corresponds an associated linear operator, denoted by $A^\top$, which is defined as the unique operator satisfying the adjoint relation
\[
\langle A(X), y \rangle = \langle X, A^\top(y) \rangle \quad \text{for all } X \in \mathbb{R}^{n \times n} \text{ and } y \in \mathbb{R}^m.
\]
Throughout this paper, $A \odot B$ denotes the Hadamard product, meaning the elementwise product of matrices or vectors of the same dimensions. Let $v \in \mathbb{R}^n$ be a vector with non-zero components ($v_i \neq 0$ for all $i$). We define the Hadamard inverse of $v$, denoted by $v^{\odot -1}$, as the vector whose $i$-th component is given by $(v^{\odot -1})_i = v_i^{-1}$. Likewise, this notation extends to any matrix $A$ with non-zero entries, by defining $(A^{\odot -1})_{ij} = A_{ij}^{-1}$.

We use the convention of multi-indices, where an index into the full tensor is represented as a collection of indices for each dimension. For a $d$-th order tensor $\mathcal{T} \in \mathbb{R}^{n_1 \times n_2 \times \dots \times n_d}$, a specific entry is denoted by $\mathcal{T}(\alpha) = \mathcal{T}(\alpha_1, \alpha_2, \dots, \alpha_d)$. Here, we use a multi-index, i.e., $\alpha = (\alpha_1, \dots, \alpha_d)$, where each entry has the index range $\alpha_k = 1, \ldots, n_k$.

\subsection{Semidefinite Programming and Interior-Point Methods}

Semidefinite programming constitutes a class of convex optimisation problems, distinguished by a matrix variable constrained to the cone of positive semidefinite matrices. The canonical form of such a problem is presented below.

\begin{definition}[SDP] \label{def:sdpprimal}
    Let $\mathbb{S}_n(\mathbb{F})$ be the set of symmetric $n \times n$ matrices over the field $\mathbb{F}$. An SDP is of the form
    \begin{align*}
        \begin{array}{lll}
            {\max_{X}} &  C \bullet X \\
             \textnormal{s.t. } & A_i \bullet X = b_i, &\textnormal{for } i = 1, ..., m \\
                                & A_i \bullet X \leq b_i, &\textnormal{for } i = m+1, ..., m+q \\
             & X \succeq 0,
        \end{array}
    \end{align*}
    where $X, C, A_i \in \mathbb{S}_n(\mathbb{F})$ and $b_i \in \mathbb{F}$.
\end{definition}

Associated with the primal problem in \Cref{def:sdpprimal} is a dual problem. It is given as follows.

\begin{definition}[The Dual SDP]\label{def:sdpdual}
    For the primal SDP with equality and inequality constraints, define the concatenated dual vector $ \lambda = \begin{bmatrix} y \\ t \end{bmatrix}$ {, where $y \in \mathbb{R}^{m}$ is the dual variable corresponding to the equality constraints and $t \in \mathbb{R}_+^{q}$ is the dual variable corresponding to the inequality constraints.} The corresponding Lagrangian dual problem is
    \begin{align*}
        \begin{array}{ll}
            {\min_{\lambda, Z}} & \lambda^\top b \\
            \textnormal{s.t. } & {\sum_{i=1}^{m+q} \lambda_i A_i - Z = C} \\
                               & Z \succeq 0.
        \end{array}
    \end{align*}
    Here, $ Z \in \mathbb{S}_n(\mathbb{F}) $ is the slack variable and  again $C, A_i \in \mathbb{S}_n(\mathbb{F})$ and $b_i \in \mathbb{F}$.
\end{definition}

The dual problem serves as a certificate of optimality: if a feasible primal-dual pair $(X, \lambda, Z)$ satisfies the complementary slackness conditions $\langle X, Z \rangle = 0$ and {$t^\top(b_\textnormal{ineq}-A_\textnormal{ineq}(X))=0$}, then $X$ and $(\lambda, Z)$ are optimal for the primal and dual problems, respectively.
The convergence and efficiency of many SDP algorithms, particularly interior-point methods, rely on certain regularity conditions. We make the two standard assumptions.

\begin{assumption}
    There exists a primal feasible point $X \succ 0$ for \cref{def:sdpprimal}, and a dual feasible point {$(\lambda, Z)$ with $Z \succ 0$ and $t > 0$} for \cref{def:sdpdual}.
\end{assumption}

\begin{assumption}
    The matrices $A_i$ of \cref{def:sdpprimal} are linearly independent.
\end{assumption}

The expressive power of the SDP framework, as defined in \cref{def:sdpprimal} and \ref{def:sdpdual}, allows for the modelling and relaxation of numerous computationally hard problems. We evaluate on the following problem relaxations.

\begin{example}[Maximum Cut]\label{ex:orig_maxcut}
    The MaxCut problem involves an undirected graph $G=(V, E)$ with a weighted adjacency matrix $C$ with $C_{ij} \in \mathbb{R}_+$. The task is to find a partitioning of the graph such that the sum of the weights of the edges going from one partition to the other is maximal. Let $n = |V|$ and assign a variable $x_i$ to every node. We say $x_i$ is in one part if $x_i=1$. Otherwise, if $x_i=-1$, it is in the other part, i.e., an edge $(i, j)$ is in the cut if $x_i= \pm 1$ and $x_j = \mp 1$. The optimisation problem is
    \begin{align*}
        \begin{array}{ll}
            \max_x & \frac{1}{2} \sum_{1 \leq i < j \leq n} C_{ij} (1-x_i x_j) \\
            \text{s.t.} & x_i \in \{-1, 1\}, \textnormal{ for } i = 1, ..., n.
        \end{array}
    \end{align*}
    The equivalent SDP relaxation is given as
    \begin{align*}
    \begin{array}{ll}
        \max_{X \in \mathbb{S}_n(\mathbb{R})} & \frac{1}{4} L \bullet X \\
        \text{s.t.} & \textnormal{diag}(X) = \mathbf{1}, \\
                    & \textnormal{rank}(X) = 1, \\
                    & X \succeq 0,
    \end{array}
    \end{align*}
    where $L \in \mathbb{S}_n(\mathbb{R})$ is the weighted Laplacian of $G$. The rank constraint is often relaxed or even removed.
\end{example}

\begin{example}[Maximum Stable Set]\label{ex:orig_stable_set}
A stable set of a graph $G = (V, E)$, $n=|V|$ is identified with a vector $x \in \textnormal{STAB}(G) = \{x \in \{0, 1\}^n : x_i + x_j \leq 1, \forall (i, j) \in E\}$. The stability number of a graph $G$ is the size of the largest of these vectors and is denoted by $\alpha(G)$, i.e.,
\begin{align*}
    \begin{array}{ll}
        \max_x & \sum_{i=1}^n x_i \\
        \textnormal{s.t. } & x \in \textnormal{STAB}(G).
    \end{array}
\end{align*}
The computation of  $\alpha(G)$ is NP-hard. However, one can compute an upper bound of $\alpha(G)$ by calculating the Lovász theta $\vartheta(G)$ as the optimal objective value to
    \begin{align*}
        \begin{array}{ll}
            \max_{X \in \mathbb{S}_n(\mathbb{R})} & J \bullet X \\
            \text{s.t. } & \textnormal{tr}(X) = 1 \\
             & X_{ij} = 0, \quad (i, j) \in E \\
             & X \succeq 0.
        \end{array}
    \end{align*}
    Here, $J$ is a matrix of all-one entries of dimension $n \times n$.
\end{example}

\begin{example}[Correlation Clustering]\label{ex:corr_clust}
    Correlation Clustering is a variation of the Max $k$-Cut problem. It involves an undirected graph $G=(V, E)$ with a weighted adjacency matrix $C$ with $C_{ij} \in \mathbb{R}$. The task is to find a multi-partitioning of the graph such that the maximal amount of similar nodes is put in the same cluster and dissimilar nodes are put into different clusters. Let $n=|V|$ and {$e_k \in \mathbb{R}^n$ be the vector with $1$ in the $k$-th coordinate and $0$s everywhere else. Assign a variable $x_v \in \{e_1, ..., e_n\}$ to every node $v \in V$. We say $x_v$ is in cluster $k$ if $x_v=e_k$.} The optimisation problem is
    \begin{align*}
        \begin{array}{lll}
            \max_{x} & {\sum_{1 \leq u < v \leq n} C^+_{uv} (x_u \cdot x_v) + C^-_{uv} (1-x_u \cdot x_v)} \\
            \text{s.t.} & x_v \cdot x_v = 1, \quad \textnormal{ for all }  v = 1, ..., n, \\
            & {x_u \cdot x_v \geq 0, \quad (u, v) \in E.}
        \end{array}
    \end{align*}
    The respective SDP relaxation is given as
    \begin{align*}
    \begin{array}{ll}
        \max_{X \in \mathbb{S}_n(\mathbb{R})} & (C^+ + L^-) \bullet X \\
        \text{s.t.} & \textnormal{diag}(X) = \mathbf{1}, \\
                    & X_{uv} \geq 0, \quad (u, v) \in E, \\
                    & X \succeq 0,
    \end{array}
    \end{align*}
    where $L^- \in \mathbb{S}_n(\mathbb{R})$ is the Laplacian of $C^-$.
\end{example}

Interior-point methods are the most prevalent and robust algorithms for solving SDPs. These methods operate by iteratively solving the Karush-Kuhn-Tucker (KKT) optimality conditions.
For each $ \mu > 0 $, we define the corresponding \emph{log-barrier Lagrangian}:
\begin{align}
\mathcal{L}_\mu(X, y, t, Z) =\
    &\ b_\textnormal{eq}^\top y + b_\textnormal{ineq}^\top t
    - \mu \left( \ln \det Z + \mathbf{1}^\top \ln t \right) \nonumber \\
    &+ \left\langle Z + C - A_\textnormal{eq}^\top(y) - A_\textnormal{ineq}^\top(t),\ X \right\rangle,
\label{eq:logbarrier-lagrangian}
\end{align}
where $\ln t$ is the element-wise natural logarithm and $\mathbf{1} \in \mathbb{R}^q$ denotes the vector of all ones. We aggregate the equality constraint matrices indexed from $1$ to $m$ into a linear operator $A_{\textnormal{eq}}$ acting on the appropriate vector space. Similarly, the inequality constraint matrices indexed from $m+1$ to $m+q$ are collected into the operator $A_{\textnormal{ineq}}$, along with their respective right-hand side vectors.
The first-order optimality (KKT) conditions of this Lagrangian are:
\begin{align}
\nabla_X \mathcal{L}_\mu &= Z + C - A_\textnormal{eq}^\top(y) - A_\textnormal{ineq}^\top(t) = 0 \label{eq:kkt-grad-x} \\
\nabla_y \mathcal{L}_\mu &= b_\textnormal{eq} - A_\textnormal{eq}(X) = 0 \label{eq:kkt-grad-y} \\
\nabla_t \mathcal{L}_\mu &= b_\textnormal{ineq} - A_\textnormal{ineq}(X) - \mu t^{\odot -1} = 0 \label{eq:kkt-grad-t} \\
\nabla_Z \mathcal{L}_\mu &= X - \mu Z^{-1} = 0 \label{eq:kkt-grad-z}
\end{align}
Using these KKT conditions, an interior-point method proceeds as shown in \cref{alg:pdipm_sdp}.

\begin{algorithm}[t]
\caption{Primal-Dual Interior-Point Method for SDPs}
\label{alg:pdipm_sdp}
\begin{algorithmic}[1]
\Require Symmetric matrix $C$, Constraint operators $A_\mathrm{eq}, A_\mathrm{ineq}$, right-hand side $b_\mathrm{eq}, b_\mathrm{ineq}$, initial $(X_0, Z_0) \succ 0, t_0 > 0$, tolerance $\epsilon_c$
\Ensure Approximate solution $(X,Z)$
\While{$ \mu := \frac{\langle X, Z \rangle + t^\top (b_\mathrm{ineq} - A_\mathrm{ineq}(X))}{n + q} > \epsilon_c$, residuals of \cref{eq:kkt-grad-x,eq:kkt-grad-y} large}
    \State Compute residuals of perturbed KKT conditions \cref{eq:kkt-grad-x,eq:kkt-grad-y,eq:kkt-grad-t,eq:kkt-grad-z} \label{algline:residuals}
    \State Compute Newton's directions $(\Delta X, \Delta y, \Delta t, \Delta Z).$ \label{algline:kkt_solve}
    \State Find maximum step size $\alpha_p, \alpha_d \in (0,1]$ s.t. \label{algline:step_size_1}
    \State $\quad X+\alpha_p \Delta X \succ 0, \quad t + \alpha_d \Delta t > 0, \quad Z+\alpha_d \Delta Z \succ 0$ \label{algline:step_size_2}
    \State Update variables: \label{algline:update_vars_1}
    \State $ \quad
        X \gets X + \alpha_p \Delta X, \quad
        (y, t, Z) \gets (y, t, Z) + \alpha_d (\Delta y, \Delta t, \Delta Z)$ \label{algline:update_vars_2}
\EndWhile
\State \Return $(X, Z)$
\end{algorithmic}
\end{algorithm}

Primal-dual interior-point methods follow a central path defined by the perturbed KKT conditions \cref{eq:kkt-grad-t} and \cref{eq:kkt-grad-z}
\[
XZ = \mu I, \quad t \odot (b_\textnormal{ineq} - A_\textnormal{ineq}(X)) = \mu \mathbf{1},
\]
where $\mu > 0$ is the centring parameter. The central path converges smoothly to an optimal solution as $\mu \to 0$. At each iteration, the algorithm computes a Newton search direction $(\Delta X, \Delta y, \Delta t, \Delta Z)$ by solving the KKT system and updates the iterates along this direction while ensuring that $X$ and $Z$ remain positive semidefinite and $t > 0$. The determination of the maximal permissible step length is therefore a critical subproblem, which is solved as the following eigenvalue problem.

\begin{theorem}[Optimal Step Size, \cite{Alizadeth1998}] \label{th:step_size}
    If the problem of selecting the optimal step size is posed as
    \begin{align*}
        \Hat{\alpha} = \sup_{0 \leq \alpha \leq 1} \{ \alpha : M + \alpha \Delta M \succeq 0 \},
    \end{align*}
    where $M \succ 0$, $\Delta M$ is symmetric, and $L = \textnormal{cholesky}(M)$, then the optimal step size is given by
    \begin{align*}
        \Hat{\alpha} = \begin{cases}
            1 & \textnormal{if } \Delta M \succeq 0 \\
            \lambda^{-1}_\textnormal{max}(-L^{-1} \Delta M L^{-T}) & \textnormal{ otherwise.}
        \end{cases}
    \end{align*}
\end{theorem}

\subsection{Tensor-Train Decomposition}\label{sec:ttd}
Let us first recall a general definition of the tensor-train decomposition \cite{Oseledets2011}. Its application to the SDP will be detailed in the next section.
A $d$-th order tensor $\mathcal{T} \in \mathbb{R}^{n_1 \times n_2 \times \dots \times n_d}$ requires storing $\prod_{k=1}^d n_k$ elements, a number that grows exponentially with the order $d$. This phenomenon is known as the \emph{curse of dimensionality}. Tensor decompositions offer a powerful solution by approximating a high-order tensor with a structured, low-parametric representation. In particular, the TT decomposition expresses a tensor as a product of smaller third-order cores.

\begin{definition}[Tensor-Train]\label{def:tt}
    A tensor $\mathcal{T} \in \mathbb{R}^{\times^d_{k=1} n_k}$ in the tensor-train format has entries
    \begin{align*}
        \mathcal{T}(\alpha_1, \ldots, \alpha_d) = T_1(\alpha_1) \cdot T_2(\alpha_2) \cdots T_d(\alpha_d),
    \end{align*}
    where $T_k(\alpha_k) \in \mathbb{F}^{r_{T_k} \times r_{T_{k+1}}}$ with some $r_{T_k} \in \mathbb{N}^+$, $r_{T_1} = r_{T_{d+1}} = 1,$ such that $T_k \in \mathbb{R}^{r_{T_k} \times n_k \times r_{T_{k+1}}},$ $k=1,\ldots,d.$ For brevity, this is often also denoted as $\mathcal{T} = [\![ T_1, T_2, \ldots, T_d ]\!]$.
    The factors $T_1, T_2, \ldots, T_d$ are referred to as TT cores, and $\rankTT(\mathcal{T}) := [r_{T_2}, \ldots, r_{T_{d}}]^\top$ are called TT ranks.
\end{definition}

\begin{definition}[TT Vector]\label{def:tt_vector}
    Given a tensor $\mathcal{T} \in \mathbb{R}^{\times^d_{k=1} n_k}$ in the TT format, its vectorisation $\mathrm{vec}(\mathcal{T}) \in \mathbb{R}^{\prod_{k=1}^d n_k}$ is called a TT vector.
\end{definition}
Note that the TT tensor is linear in the elements of one TT core at a time. To see this, we can introduce the partial TT formats.
\begin{definition}[Interface Matrix]\label{def:interface}
Let $\mathcal{T} = [\![ T_1, T_2, \ldots, T_d ]\!] \in \mathbb{R}^{\times^d_{j=1} n_j}$ be a tensor in the TT format.
The matrices $\mathcal{T}^{<k} \in \mathbb{R}^{\prod_{j=1}^{k-1}n_j \times r_{T_k}}$ and $\mathcal{T}^{>k} \in \mathbb{R}^{\prod_{j=k+1}^{d}n_j \times r_{T_{k+1}}}$ with rows
\begin{align*}
    \mathcal{T}^{<k}(\alpha_1,\ldots,\alpha_{k-1}) & = T_1(\alpha_1) \cdots T_{k-1}(\alpha_{k-1}), & k&=2,\ldots,d, \\
    \mathcal{T}^{>k}(\alpha_{k+1},\ldots,\alpha_d) & = T_{k+1}(\alpha_{k+1}) \cdots T_d(\alpha_d), & k&=1,\ldots,d-1,
\end{align*}
are called left resp. right interface matrices.
Without loss of generality, we extend the notation to include $\mathcal{T}^{<1} = \mathcal{T}^{>d} = 1.$
\end{definition}
\begin{definition}[Frame Matrix]\label{def:frame}
Let $\mathcal{T} = [\![ T_1, T_2, \ldots, T_d ]\!] \in \mathbb{R}^{\times^d_{j=1} n_j}$ be a tensor in the TT format. The frame matrix reads
\begin{align*}
\mathcal{T}_{\neq k} &= \mathcal{T}^{<k} \otimes I_{n_k} \otimes \mathcal{T}^{>k} \in \mathbb{R}^{\prod_{j=1}^dn_j \times r_{T_k}n_kr_{T_{k+1}}}, & k&=1,\ldots,d.
\end{align*}
\end{definition}
{Throughout this section, $\mathrm{vec}(\mathcal{T})$ denotes the big-endian lexicographic vectorisation with multi-index order $(\alpha_1, \ldots, \alpha_d)$. All Kronecker products in the interface and frame matrices are written with respect to this ordering. Under a different vectorisation convention, the same identities hold up to the corresponding permutation.}
\begin{remark}
One can verify that
\begin{align*}
\mathrm{vec}(\mathcal{T}) & = \mathcal{T}_{\neq k} \mathrm{vec}(T_k), & k&=1,\ldots,d.
\end{align*}
{If a little-endian vectorisation convention is used, the ordering of the Kronecker factors in $\mathcal{T}_{\neq k}$ must be permuted accordingly.}
This linear embedding, with the frame matrix playing the role of the basis in the tensor space, will be crucial for alternating linear solvers for the KKT equations.
\end{remark}

To apply TT methodology to SDPs, we also require linear operators that efficiently act on tensors. Since a tensor $\mathcal{T}$ is bijectively mapped to $\mathrm{vec}(\mathcal{T})$, such operators can be written as matrices acting on $\mathrm{vec}(\mathcal{T})$.
They are represented in a slightly different TT format to keep the TT ranks low (in particular, to represent the identity matrix with TT rank 1).
\begin{definition}[TT Matrix]
    A matrix $\mathcal{A}= [\![ A_1, A_2, \ldots, A_d ]\!]  \in \mathbb{R}^{\mathrm{M} \times \mathrm{N}}$ where $\mathrm{M} = \prod^d_{k=1} m_{k}$ and $\mathrm{N}=\prod^d_{k=1} n_k$ is said to be in the tensor-train format if its entries are expressed as follows,
    \begin{align*}
        \mathcal{A}(\alpha_1, \ldots, \alpha_d;~\beta_1,\ldots,\beta_d) = A_1(\alpha_1,\beta_1) \cdot A_2(\alpha_2,\beta_2) \cdots A_d(\alpha_d,\beta_d),
    \end{align*}
    where $A_k(\alpha_k,\beta_k) \in \mathbb{R}^{r_{A_k} \times r_{A_{k+1}}}$ with TT ranks $r_{A_k} \in \mathbb{N}^+$, $r_{A_1} = r_{A_{d+1}} = 1,$ $\alpha_k=1,\ldots,m_k,$ $\beta_k=1,\ldots,n_k,$ such that $A_k \in \mathbb{R}^{r_{A_k} \times (m_{k} \times n_k) \times r_{A_{k+1}}},$ $k=1,\ldots,d.$
\end{definition}

Many fundamental operations on tensors in the tensor-train format, such as addition, matrix multiplication or Hadamard products, naturally yield a resultant tensor with larger TT ranks. This rank growth inflates storage and computational costs, undermining the core efficiency benefits of the decomposition. TT rounding \cite{Oseledets2011} constitutes a fundamental compression step applied subsequent to such operations. It may be used either to obtain a quasi-optimal tensor of prescribed TT ranks or, alternatively, to construct an approximation subject to a specified error tolerance $\delta > 0$. In the latter case, we denote the tensor-train rounding of a tensor $\mathcal{T}$ by $\widetilde{\mathcal{T}} = \mathrm{round}_\mathrm{TT}(\mathcal{T}, \delta)$.

One of the first iterative schemes for solving equations in the TT format was
the Density Matrix Renormalisation Group (DMRG) \cite{White1992} for eigenvalue problems. Similar algorithms for linear systems were developed under the name of Alternating Linear Schemes (ALS) \cite{Holtz2012}.

Consider a linear system $\mathcal{A} \, \mathrm{vec}(\mathcal{X}) = \mathrm{vec}(\mathcal{B})$ in TT format with a symmetric $\mathcal{A} \succ 0$. Its solution can be found as the minimiser of the energy function
\[
    \mathscr{J}(\mathcal{X}) := \langle \mathrm{vec}(\mathcal{X}), \mathcal{A} \, \mathrm{vec}(\mathcal{X}) \rangle - 2 \langle \mathcal{X}, \mathcal{B} \rangle.
\]
Substituting $\mathcal{X} = [\![X_1,\ldots, X_d]\!]$ in the TT format and using the associated frame matrices, we obtain the restricted energy function on the $k$th TT core, $k=1,\ldots,d$
\[
    \mathscr{J}(\mathcal{X}) = \mathscr{\hat J}(X_k) := \langle \mathcal{X}^\top_{\neq k} \mathcal{A} \mathcal{X}_{\neq k} \, \mathrm{vec}(X_k), \mathrm{vec}(X_k)\rangle - 2 \langle \mathcal{X}^\top_{\neq k} \mathrm{vec}(\mathcal{B}), \mathrm{vec}(X_k) \rangle.
\]
This energy is minimised by solving the projected linear system
$(\mathcal{X}^\top_{\neq k} \mathcal{A} \mathcal{X}_{\neq k}) \mathrm{vec}(X_k) = \mathcal{X}^\top_{\neq k} \mathrm{vec}(\mathcal{B})$ iteratively for each $k$. However, when employing an ALS, it is important to ensure that the projectors are not degenerate. { A projector $\mathcal{X}_{\neq k}$ is called degenerate if $\mathcal{X}^{\top}_{\neq k} \mathcal{A} \mathcal{X}_{\neq k}$ is singular, which causes the projected linear system to be unsolvable or ill-conditioned. This can occur even when $\mathcal{A}$ itself is nonsingular.
\begin{example}\label{ex:singular_projection}
Let
\begin{align*}
    \mathcal{A} = \begin{bmatrix}0 & 1 \\ 1 & 0\end{bmatrix}, \qquad \mathcal{X}_{\neq k} = \begin{bmatrix}0 \\ 1\end{bmatrix}.
\end{align*}
One can readily verify that $\mathcal{X}_{\neq k}^\top \mathcal{A} \mathcal{X}_{\neq k} = 0$, so the projected system $\mathcal{X}_{\neq k}^\top \mathcal{A} \mathcal{X}_{\neq k}\, \mathrm{vec}(X_k) = \mathcal{X}_{\neq k}^\top \mathrm{vec}(\mathcal{B})$ is singular despite $\mathcal{A}$ being nonsingular.
\end{example}
}

This issue arises in particular when dealing with indefinite block matrices, like those in the KKT system. Consequently, standard ALS solvers are inadequate for solving the KKT system arising in an interior-point method. A way to overcome this limitation is by employing a Block-AMEn scheme, akin to the approach in~\cite{Benner2016, Benner2020}. This is accomplished by representing the solution in a block tensor-train format
and by jointly optimising over multiple solution components simultaneously, which prevents the projected matrix from becoming singular.

\begin{definition}[Block Tensor-Train]
    A set of tensors $\mathcal{W}_{\ell} \in \mathbb{R}^{\times^d_{j=1} n_j}, \ell \in \mathbb{N}^+$ in block tensor-train format with some $k=1,\ldots,d$ has entries
    \begin{align*}
        \mathcal{W}_{\ell}(\alpha_1, \ldots, \alpha_d) = W_1(\alpha_1) \cdot W_2(\alpha_2) \cdots \hat W_k(\alpha_k, \ell) \cdots W_d(\alpha_d),
    \end{align*}
    where $\hat W_k \in \mathbb{R}^{r_{W_k} \times n_k \times \ell \times r_{W_{k+1}}}$ with some $r_{W_k} \in \mathbb{N}^+$, $r_{W_1} = r_{W_{d+1}} = 1,$ and $W_s \in \mathbb{R}^{r_{W_s} \times n_s \times r_{W_{s+1}}},$ $s \neq k$.
\end{definition}

This representation allows us to introduce the interface matrices $\mathcal{W}^{<k}, \mathcal{W}^{>k}$ and the frame matrix $\mathcal{W}_{\neq k}$ using Definitions~\ref{def:interface}--\ref{def:frame} verbatim.
However, for a different $k=1,\ldots,d$, we need to move the block index $\ell$ to the corresponding core. This can be accomplished using SVD, by realising one step of the TT rounding procedure to move $\ell$ to $W_{k+1}$ or $W_{k-1}$. We can repeat this until $\ell$ appears in the desired core. This also facilitates a DMRG-like rank adaptation. However, in our experiments, we found that additional enrichment of the TT cores by TT cores of the approximate residual, as in the AMEn algorithm \cite{Dolgov2014AMeN}, can improve the convergence.

\section{A Tensor-Train Interior-Point Method for SDPs}\label{sec:ttipm}
In this section, we present the methodology behind the TT IPM. We first describe how to represent the problems in the tensor-train format, then explain how to solve the resulting KKT system, and finally discuss the computation of the optimal step size and practical considerations for implementing the method.

\subsection{Representation of Matrices and Operators in the TT Format}
Throughout the TT IPM, we are concerned with square matrices with $n := n_k = n_{k+d}$ for all $k$.

The solution of the KKT conditions in the TT setting predominantly involves working with \emph{suitably} vectorised TT matrices.

\begin{definition}[TT Vectorisation]\label{def:tt_vectorisation}
    For given TT matrix $\mathcal{M} = [\![ M_1, ..., M_d ]\!]$ with $M_k \in \mathbb{R}^{r_{M_k} \times (n \times n) \times r_{M_{k+1}}}$, its vectorisation is defined as a TT vector $\vecTT(\mathcal{M}) := \mathrm{vec}([\![ V_1, ..., V_d ]\!]),$ where
    \begin{align*}
    V_k(s_k,:,s_{k+1}) & := \mathrm{vec}(M_k(s_k, :, :, s_{k+1})), && s_k=1,\ldots,r_{M_k}, \quad k=1,\ldots,d,
    \end{align*}
    and $\mathrm{vec}(\cdot)$ on the right-hand side is the standard vectorisation of a matrix.
\end{definition}

\begin{remark}
Both $\mathrm{vec}(\mathcal{M})$ and $\vecTT(\mathcal{M})$ yield vectors of dimension $n^{2d}$, but they differ in their ordering of entries. Specifically, $\mathrm{vec}(\mathcal{M})$ is indexed as
\[
\mathrm{vec}(\mathcal{M})(\alpha_1,\ldots,\alpha_d,\beta_1,\ldots,\beta_d),
\]
whereas $\vecTT(\mathcal{M})$ is indexed in an interleaved manner as
\[
\vecTT(\mathcal{M})(\alpha_1,\beta_1,\ldots,\alpha_d,\beta_d).
\]
\begin{example}
Consider $d=n=2$ and
$$
\mathcal{M} = \begin{bmatrix}
1 & 1 & 2 & 2 \\
1 & 1 & 2 & 2 \\
3 & 3 & 4 & 4 \\
3 & 3 & 4 & 4
\end{bmatrix}.
$$
Assuming the big-endian convention as above, the two vectorisations read
\begin{align*}
\mathrm{vec}(\mathcal{M}) & = \begin{bmatrix}1 \ 1 \ 2 \ 2 \ 1 \ 1 \ 2 \ 2 \ 3 \ 3 \ 4 \ 4 \ 3 \ 3 \ 4 \ 4\end{bmatrix} \\
\vecTT(\mathcal{M}) & = \begin{bmatrix}1 \ 1 \ 1 \ 1 \ 2 \ 2 \ 2 \ 2 \ 3 \ 3 \ 3 \ 3 \ 4 \ 4 \ 4 \ 4\end{bmatrix}.
\end{align*}
\end{example}

Empirically, for matrices $\mathcal{M}$ arising in semidefinite programming, the tensor-train ranks associated with $\vecTT(\mathcal{M})$ are typically much lower than those corresponding to $\mathrm{vec}(\mathcal{M})$.
\end{remark}

This results in the tensorised SDP problem
\[
{\max_{\mathcal{X}}} \langle \mathcal{C}, \mathcal{X} \rangle \quad \text{s.t.} \quad \mathbfcal{A}_{\text{eq}}\vecTT(\mathcal{X}) = \mathcal{B}_{\text{eq}}, \quad \mathbfcal{A}_{\text{ineq}}\vecTT(\mathcal{X}) \leq \mathcal{B}_{\text{ineq}}, \quad \mathcal{X} \succeq 0,
\]
where $\mathcal{X},\mathcal{C} \in \mathbb{R}^{n^d \times n^d}$, $\mathcal{B}_{\text{eq}}, \mathcal{B}_{\text{ineq}}\in \mathbb{R}^{n^{2d}}$, and $\mathbfcal{A}_\textnormal{eq}, \mathbfcal{A}_\textnormal{ineq} \in \mathbb{R}^{n^{2d} \times n^{2d}}$. The SDP equality constraints $1, \ldots, m$ and the inequality constraints $1, \ldots, q$ are embedded into TT matrices and TT vectors of the right-hand sides of the form
\begin{align*}
    \mathbfcal{A}_{\nu} &= [\![ \mathbf{A}_{\nu,1}, \ldots, \mathbf{A}_{\nu,d} ]\!], & \mathbf{A}_{\nu,k} &\in \mathbb{R}^{r_{\mathbf{A}_{\nu,k}} \times (n^2 \times n^2) \times r_{\mathbf{A}_{\nu,k+1}}}, \\ \mathcal{B}_{\nu} &= [\![ B_{\nu,1}, \ldots, B_{\nu,d} ]\!], & B_{\nu,k} &\in \mathbb{R}^{r_{B_{\nu,k}} \times n^2 \times r_{B_{\nu,k+1}}},
\end{align*}
where $\nu \in \{\text{eq}, \text{ineq}\}$. We highlight the TT operators in bold to emphasise the dimensionality differences. Likewise, the objective matrix is
\begin{align*}
    \mathcal{C} &= [\![C_1, C_2, \ldots C_d ]\!], & C_k \in \mathbb{R}^{r_{C_k} \times (n \times n) \times r_{C_{k+1}}}.
\end{align*}

\subsection{Representation of SDP Variables in the TT format}
The dual slack variable is of the same size as the primal variable, and can thus be stored as a TT matrix $\mathcal{Z} \in \mathbb{R}^{n^d \times n^d}$.
In contrast, the Lagrange multipliers have different sizes ($m$ and $q$).
The Block-AMEn scheme requires all solution components to be of the same size.
Thus, we embed the Lagrange multipliers into larger TT matrices $\mathcal{Y}, \mathcal{T} \in \mathbb{R}^{n^d \times n^d}$.

\begin{definition}[Idle Lagrange Multiplier]
    {An entry of the Lagrange multiplier TT matrices $\mathcal{Y}$ or $\mathcal{T}$ is referred to as \emph{idle} if it is not associated with any equality or inequality constraint, respectively.}
\end{definition}

{
A heterogeneous formulation could avoid this embedding by keeping the equality and inequality multipliers as TT vectors on their own tensorised constraint index sets, of dimensions $m$ and $q$, rather than embedding them into the common $n^d \times n^d$ TT matrix index space. However, the KKT unknowns would then belong to different TT spaces, so the standard Block-AMEn formulation with one shared frame matrix $\mathcal{W}_{\neq k}$ would no longer apply directly. The projected KKT blocks would involve mixed-frame contractions, such as $\mathcal{W}_{Y,\neq k}^{\top}\mathbfcal{A}_{\mathrm{eq}}\mathcal{W}_{X,\neq k}$, and may invalidate both the theoretical invertibility of the projected matrix and the convergence of the overall scheme.
We therefore use the padded common-index representation for reasons of computational efficiency and to retain the shared-frame projected KKT structure analysed below.

The placement of the idle multipliers is not unique. Since $\mathcal{X}$ is symmetric, at most $n^d(n^d+1)/2$ independent constraint multipliers are required, so they could be stored in one triangular part of $\mathcal{Y}$ and $\mathcal{T}$. In practice, this triangular placement is unfavourable in TT format: the resulting masks are less aligned with the graph adjacency operators used in Section~\ref{sec:conceptualisation} and typically lead to larger TT ranks. We therefore store mirrored copies of the multipliers on both sides of the diagonal. This keeps the associated constraint operators self-adjoint and preserves TT rank bounds close to those of the underlying adjacency matrices. In summary, we represent the SDP variables respectively as
}
\begin{align*}
    \mathcal{X} &= [\![X_1, X_2, \ldots X_d ]\!], & X_k \in \mathbb{R}^{r_{X_k} \times (n \times n) \times r_{X_{k+1}}}, \\
    \mathcal{Y} &= [\![Y_1, Y_2, \ldots Y_d ]\!], & Y_k \in \mathbb{R}^{r_{Y_k} \times (n \times n) \times r_{Y_{k+1}}}, \\
    \mathcal{T} &= [\![T_1, T_2, \ldots T_d ]\!], & T_k \in \mathbb{R}^{r_{T_k} \times (n \times n) \times r_{T_{k+1}}}, \\
    \mathcal{Z} &= [\![Z_1, Z_2, \ldots Z_d ]\!], & Z_k \in \mathbb{R}^{r_{Z_k} \times (n \times n) \times r_{Z_{k+1}}}.
\end{align*}

\subsection{Solution of the KKT system in the TT format}\label{sec:sol_kkt_tt_format}
We now present additional tensor-train constructions that are pertinent to the methodology developed in the following sections.

\begin{definition}[TT Diagonalisation]
    For given TT vector $\mathcal{V} = [\![ V_1, ..., V_d ]\!]$, its diagonalisation is
    defined as a TT matrix $\diagTT(\mathcal{V}) := [\![ M_1, ..., M_d ]\!],$ where
    \begin{align*}
    M_k(s_k,:,:,s_{k+1}) & := \mathrm{diag}(V_k(s_k, :, s_{k+1})), && s_k=1,\ldots,r_{V_k}, \quad k=1,\ldots,d,
    \end{align*}
    and $\mathrm{diag}(\cdot)$ on the right-hand side is the standard diagonal matrix with the elements of the given vector on the diagonal.
\end{definition}

\begin{remark}
    The element-wise product of given TT matrices $\mathcal{M} = [\![ M_1, ..., M_d ]\!]$ and $\mathcal{N} = [\![ N_1, ..., N_d ]\!]$ can be represented in the TT format as $\vecTT(\mathcal{M} \odot \mathcal{N}) = \diagTT(\vecTT(\mathcal{M})) \vecTT(\mathcal{N})$.
\end{remark}

\begin{definition}[Kronecker product of TT matrices]
    For given TT matrices $\mathcal{A} = [\![ A_1, ..., A_d ]\!]$ with $A_k \in \mathbb{R}^{r_{A_k} \times (n_{k+d} \times n_k) \times r_{A_{k+1}}}$ and $\mathcal{B} = [\![ B_1, ..., B_d ]\!]$ with $B_k \in \mathbb{R}^{r_{B_k} \times (n_{k+d} \times n_k) \times r_{B_{k+1}}}$, their Kronecker product is a TT matrix $\mathcal{A} \otimes \mathcal{B} := [\![ C_1, ..., C_d ]\!]$, where
    \begin{align*}
    C_k &= A_k \otimes B_k \in \mathbb{R}^{r_{A_k}r_{B_k} \times (n_{k+d}^2 \times n_k^2) \times r_{A_{k+1}}r_{B_{k+1}}}, && k=1,\ldots,d,
    \end{align*}
    using the standard Kronecker product of 4-dimensional tensors $A_k$ and $B_k$.
\end{definition}

{
For the inequality constraints, define the slack vector and its associated diagonal TT operator by
\begin{align*}
    s_{\mathrm{ineq}}(\mathcal X)
    &:=
    \mathcal{B}_{\mathrm{ineq}}
    -
    \mathbfcal{A}_{\mathrm{ineq}}\vecTT(\mathcal X),\\
    \mathcal D_{\mathrm{ineq}}(\mathcal X)
    &:=
    \diagTT\!\left(s_{\mathrm{ineq}}(\mathcal X)\right).
\end{align*}
The operator $\mathcal D_{\mathrm{ineq}}(\mathcal X)$ is the diagonal operator formed from the inequality slack vector and is distinct from the SDP dual slack matrix $\mathcal Z$. 
In the particular case of entry-wise masked inequalities appearing in our numerical examples,
the above simplify to
\[
    s_{\mathrm{ineq}}(\mathcal X)
    =
    \vecTT(M_T\odot\mathcal X),
    \qquad
    \mathcal D_{\mathrm{ineq}}(\mathcal X)
    =
    \diagTT\!\left(\vecTT(M_T\odot\mathcal X)\right),
\]
where $M_T$ is the inequality mask.
}

Using TT matrices, the tensorised KKT conditions \eqref{eq:kkt-grad-x}--\eqref{eq:kkt-grad-z}
can be written as follows:
\begin{align*}
    \nabla_\mathcal{X} \mathcal{L}_\mu(\mathcal{X}, \mathcal{Y}, \mathcal{T}, \mathcal{Z}) &= \vecTT(\mathcal{Z}) + \vecTT(\mathcal{C}) - \mathbfcal{A}_\textnormal{eq}^\top \vecTT(\mathcal{Y}) - \mathbfcal{A}_\textnormal{ineq}^\top \vecTT(\mathcal{T}) = 0 \\
    \nabla_\mathcal{Y} \mathcal{L}_\mu(\mathcal{X}, \mathcal{Y}, \mathcal{T}, \mathcal{Z}) &= \mathcal{B}_\textnormal{eq} - \mathbfcal{A}_\textnormal{eq}\vecTT(\mathcal{X}) = 0 \\
    \nabla_\mathcal{T} \mathcal{L}_\mu(\mathcal{X}, \mathcal{Y}, \mathcal{T}, \mathcal{Z}) &= {s_{\mathrm{ineq}}(\mathcal X)} - \mu\, \vecTT(\mathcal{T})^{\odot -1} = 0 \\
    \nabla_\mathcal{Z} \mathcal{L}_\mu(\mathcal{X}, \mathcal{Y}, \mathcal{T}, \mathcal{Z}) &= \mathcal{X} - \mu\, \mathcal{Z}^{-1} = 0
\end{align*}
The barrier parameter $\mu$ is chosen as
\[
    \mu = \frac{1}{n^d + q} \left( \langle \mathcal{Z}, \mathcal{X} \rangle + \langle \vecTT(\mathcal{T}), {s_{\mathrm{ineq}}(\mathcal X)} \rangle \right).
\]

However, the idle Lagrange multipliers cannot be resolved from the above equations.
To circumvent this, we introduce
\begin{align*}
    \mathbfcal{K}_\mathcal{Y} &= [\![ \mathbf{K}^{(\mathcal{Y})}_1, \ldots, \mathbf{K}^{(\mathcal{Y})}_{d} ]\!], & \mathbf{K}^{(\mathcal{Y})}_k \in \mathbb{R}^{r_{\mathbf{K}^{(\mathcal{Y})}_k} \times (n^2 \times n^2) \times r_{\mathbf{K}^{(\mathcal{Y})}_{k+1}}}, \\
    \mathbfcal{K}_\mathcal{T} &= [\![ \mathbf{K}^{(\mathcal{T})}_1, \ldots, \mathbf{K}^{(\mathcal{T})}_{d} ]\!], & \mathbf{K}^{(\mathcal{T})}_k \in \mathbb{R}^{r_{\mathbf{K}^{(\mathcal{T})}_k} \times (n^2 \times n^2) \times r_{\mathbf{K}^{(\mathcal{T})}_{k+1}}},
\end{align*}
with identity submatrices at the positions of zero padding rows in $\mathbfcal{A}_\textnormal{eq}, \mathbfcal{A}_\textnormal{ineq}$.
This allows us to introduce the additional constraints
$\mathbfcal{K}_\mathcal{Y}\vecTT( \mathcal{Y}) = \mathbf{0}$ and $\mathbfcal{K}_\mathcal{T} \vecTT(\mathcal{T}) = \mathbf{0}$.
The zeros in the right-hand sides are embedded into the corresponding elements of $\mathcal{B}_{\text{eq}},\mathcal{B}_{\text{ineq}}$.

We employ the AHO (Alizadeh–Haeberly–Overton) search direction~\cite{Alizadeth1998}, as alternative search directions generally require at least approximate eigendecompositions, which are not only computationally prohibitive in the TT format but also prone to significant rank inflation. Based on the above KKT conditions, the resulting Newton system to be solved is
\begin{align}
    \left[ \begin{array}{cccc}
        0 & - \mathbfcal{A}_\textnormal{eq}^\top & - \mathbfcal{A}_\textnormal{ineq}^\top & \mathcal{I} \otimes \mathcal{I} \\
        - \mathbfcal{A}_\textnormal{eq} & \mathbfcal{K}_\mathcal{Y} & 0 & 0 \\
         - \mathbfcal{T} \mathbfcal{A}_\textnormal{ineq} & 0 & {\mathcal D_{\mathrm{ineq}}(\mathcal X)} + \mathbfcal{K}_\mathcal{T}  & 0 \\
        \mathbfcal{L}_\mathcal{Z} & 0 & 0 & \mathbfcal{L}_\mathcal{X}
    \end{array} \right] \mathcal{W} = -\begin{bmatrix}
       \mathcal{R}_d \\
        \mathcal{R}_p \\
        \mathcal{R}_t \\
       \mathcal{R}_c
    \end{bmatrix},
    \label{eq:tt-kkt-newton-system}
\end{align}
Here, $\mathbfcal{T} := \diagTT(\vecTT(\mathcal{T}))$ with the residuals given as
\begin{align*}
    \begin{bmatrix}
       \mathcal{R}_d \\
        \mathcal{R}_p \\
        \mathcal{R}_t \\
       \mathcal{R}_c
    \end{bmatrix} := \begin{bmatrix}
       \vecTT(\mathcal{Z}) + \vecTT(\mathcal{C}) - \mathbfcal{A}_\textnormal{eq}^\top\vecTT(\mathcal{Y}) - \mathbfcal{A}_\textnormal{ineq}^\top\vecTT(\mathcal{T})\\
        \mathcal{B}_{\text{eq}} - \mathbfcal{A}_\textnormal{eq} \vecTT(\mathcal{X})  \\
        \vecTT(\mathcal{T}) \odot {s_{\mathrm{ineq}}(\mathcal X)} - \mu \, \vecTT(\mathbf{1}_\mathcal{T}) \\
       \mathbfcal{L}_\mathcal{Z}\vecTT(\mathcal{X}) - 2 \mu \, \vecTT(\mathcal{I})
    \end{bmatrix},
\end{align*}
where $\mathbf{1}_\mathcal{T}$ is the matrix of ones except for zeros at the positions corresponding to idle multipliers in $\mathcal{T}$,
and $\mathcal{W} := \begin{bmatrix} \vecTT(\Delta \mathcal{X}) & \vecTT(\Delta \mathcal{Y}) & \vecTT(\Delta \mathcal{T}) & \vecTT(\Delta \mathcal{Z}) \end{bmatrix}^\top$.

{
The displayed Newton system is written for the AHO direction. In this case $\mathbfcal{L}_\mathcal{Z} := \mathcal{I}\otimes\mathcal{Z} + \mathcal{Z}\otimes\mathcal{I}$, $\mathbfcal{L}_\mathcal{X} := \mathcal{I}\otimes\mathcal{X} + \mathcal{X}\otimes\mathcal{I}$, and the complementarity residual is $\mathcal{R}_c = \mathbfcal{L}_\mathcal{Z}\vecTT(\mathcal{X}) - 2\mu\vecTT(\mathcal{I})$. During the initial $XZ$-direction warm-up, the implementation instead uses the unsymmetrised operators and the matching residual $\mathcal{R}_c = (\mathcal{Z}\otimes\mathcal{I})\vecTT(\mathcal{X}) - \mu\vecTT(\mathcal{I})$. This direction is cheaper in TT arithmetic, while the AHO form is the one used for the nonsingularity argument below under the standard regularity assumptions \cite{Alizadeth1998}.
}

Let $\mathcal{W}_{\neq k}\in\mathbb{R}^{n^{2d}\times p_k}$, $p_k:=r_{W_k}n^2r_{W_{k+1}}$, be the frame matrix of the vectorised solution matrix in the block TT format $\mathcal{W}:= [\![ W_1, \ldots, W_d ]\!]$. Then, we begin by eliminating the block associated with $\Delta \mathcal{Z}$ via
\begin{align*}
    \textnormal{vec}(\Delta Z_k) = \mathcal{W}_{\neq k}^\top \mathbfcal{A}_\textnormal{eq}^\top \mathcal{W}_{\neq k} \textnormal{vec}(\Delta Y_k) +  \mathcal{W}_{\neq k}^\top \mathbfcal{A}_\textnormal{ineq}^\top \mathcal{W}_{\neq k} \textnormal{vec}(\Delta T_k) - \mathcal{W}_{\neq k}^\top \mathcal{R}_d.
\end{align*}

Note that $\Delta Z_k$ can therefore always be computed matrix-free. Substituting this expression into the remaining block equations eliminates $\Delta Z_k$, and the resulting projected system for the remaining updates becomes
\begin{align}
\label{eq:reduced_kkt_system}
\resizebox{\linewidth}{!}{$
    \begin{aligned}
        &\left[
        \begin{array}{ccc}
            - \mathcal{W}_{\neq k}^\top \mathbfcal{A}_\textnormal{eq} \mathcal{W}_{\neq k} &
            \mathcal{W}_{\neq k}^\top \mathbfcal{K}_\mathcal{Y} \mathcal{W}_{\neq k} & 0 \\
            - \mathcal{W}_{\neq k}^\top \mathbfcal{T}\mathbfcal{A}_\textnormal{ineq} \mathcal{W}_{\neq k} & 0 &
            \mathcal{W}_{\neq k}^\top({\mathcal D_{\mathrm{ineq}}(\mathcal X)} + \mathbfcal{K}_\mathcal{T}) \mathcal{W}_{\neq k} \\
            \mathcal{W}_{\neq k}^\top \mathbfcal{L}_\mathcal{Z} \mathcal{W}_{\neq k} &
            \mathcal{W}_{\neq k}^\top \mathbfcal{L}_\mathcal{X} \mathcal{W}_{\neq k}
            \mathcal{W}_{\neq k}^\top \mathbfcal{A}_\textnormal{eq}^\top \mathcal{W}_{\neq k} &
            \mathcal{W}_{\neq k}^\top \mathbfcal{L}_\mathcal{X} \mathcal{W}_{\neq k}
            \mathcal{W}_{\neq k}^\top \mathbfcal{A}_\textnormal{ineq}^\top \mathcal{W}_{\neq k}
        \end{array}
        \right]
        \begin{bmatrix}
            \textnormal{vec}(\Delta X_k) \\
            \textnormal{vec}(\Delta Y_k) \\
            \textnormal{vec}(\Delta T_k)
        \end{bmatrix} \\[1ex]
        &\qquad = -\begin{bmatrix}
            \mathcal{W}_{\neq k}^\top \mathcal{R}_p \\
            \mathcal{W}_{\neq k}^\top \mathcal{R}_t \\
            \mathcal{W}_{\neq k}^\top \mathcal{R}_c
            - \mathcal{W}_{\neq k}^\top \mathbfcal{L}_\mathcal{X} \mathcal{W}_{\neq k}
              \mathcal{W}_{\neq k}^\top \mathcal{R}_d
        \end{bmatrix}.
    \end{aligned}
$}
\end{align}
{
If the projected blocks are too large, we solve the local system by GMRES in a matrix-free form. During an AMEn sweep, the left and right TT environments of each operator block are accumulated and reused. A product with a formally projected block, such as $\mathcal{W}_{\neq k}^\top \mathbfcal{L}_\mathcal{Z}\mathcal{W}_{\neq k}$, is then evaluated by inserting the local GMRES vector into the active core and contracting only this core with the cached environments and the corresponding operator core. Thus each Krylov product uses the local projected action without forming the full frame $\mathcal{W}_{\neq k}$, the full operator $\mathbfcal{L}_\mathcal{Z}$, or the dense projected matrix.} In the absence of a suitable preconditioner for the matrix-free solution of the local linear systems, the tensor-train solution may become unstable and exhibit excessive rank growth. To mitigate this, we employ the accelerated GMRES method~\cite{Baker2005}. When the dimensions remain manageable, however, we can proceed with the block elimination strategy by computing a Cholesky factorisation of $\mathcal{W}_{\neq k}^\top \mathbfcal{L}_\mathcal{Z} \mathcal{W}_{\neq k}$ together with two LU decompositions of
\begin{align*}
     \Phi = \; &\mathcal{W}_{\neq k}^\top \mathbfcal{T} \mathbfcal{A}_\mathrm{ineq} \mathcal{W}_{\neq k} (\mathcal{W}_{\neq k}^\top \mathbfcal{L}_\mathcal{Z} \mathcal{W}_{\neq k})^{-1}
         \mathcal{W}_{\neq k}^\top \mathbfcal{L}_\mathcal{X} \mathcal{W}_{\neq k} \mathcal{W}_{\neq k}^\top \mathbfcal{A}^\top_\mathrm{ineq} \mathcal{W}_{\neq k} \\
         & +  \mathcal{W}_{\neq k}^\top({\mathcal D_{\mathrm{ineq}}(\mathcal X)} + \mathbfcal{K}_\mathcal{T})\mathcal{W}_{\neq k},
\end{align*}
and the Schur complement
\begin{align*}
        S = \; &\mathcal{W}_{\neq k}^\top \mathbfcal{A}_\mathrm{eq} \mathcal{W}_{\neq k} (\mathcal{W}_{\neq k}^\top \mathbfcal{L}_\mathcal{Z} \mathcal{W}_{\neq k})^{-1}
         \mathcal{W}_{\neq k}^\top \mathbfcal{L}_\mathcal{X} \mathcal{W}_{\neq k} \mathcal{W}_{\neq k}^\top \mathbfcal{A}_\mathrm{eq}^\top \mathcal{W}_{\neq k} \\
         & + \mathcal{W}_{\neq k}^\top \mathbfcal{K}_\mathcal{Y} \mathcal{W}_{\neq k} -  \mathcal{W}_{\neq k}^\top \mathbfcal{A}_\mathrm{eq} \mathcal{W}_{\neq k} \Phi^{-1} \mathcal{W}_{\neq k}^\top \mathbfcal{T} \mathbfcal{A}_\mathrm{ineq} \mathcal{W}_{\neq k}.
\end{align*}
Using these factorisations, the updates to the system variables are obtained by a sequence of triangular solves, with the first being
\begin{align*}
    S \mathrm{vec}(\Delta Y_k) = \; &
    \mathcal{W}_{\neq k}^\top \mathbfcal{A}_\mathrm{eq} \mathcal{W}_{\neq k}
    (\mathcal{W}_{\neq k}^\top \mathbfcal{L}_\mathcal{Z} \mathcal{W}_{\neq k})^{-1} (\mathcal{W}_{\neq k}^\top \mathbfcal{L}_\mathcal{X} \mathcal{W}_{\neq k}
        \mathcal{W}_{\neq k}^\top \mathcal{R}_d
        - \mathcal{W}_{\neq k}^\top \mathcal{R}_c) \\
    & - \mathcal{W}_{\neq k}^\top \mathcal{R}_p.
\end{align*}
With $\Delta Y_k$ known, we proceed to find the update for the inequality multipliers
\begin{align*}
    \Phi \mathrm{vec}(\Delta T_k) = \; &
    \mathcal{W}_{\neq k}^\top \mathbfcal{T} \mathbfcal{A}_\mathrm{ineq} \mathcal{W}_{\neq k}
    \Big( (\mathcal{W}_{\neq k}^\top \mathbfcal{L}_\mathcal{Z} \mathcal{W}_{\neq k})^{-1} \\
    & \quad \big(\mathcal{W}_{\neq k}^\top \mathcal{R}_c
      - \mathcal{W}_{\neq k}^\top \mathbfcal{L}_\mathcal{X} \mathcal{W}_{\neq k}
        \mathcal{W}_{\neq k}^\top \mathcal{R}_d \big)
      - \mathrm{vec}(\Delta Y_k) \Big) \\
    & - \mathcal{W}_{\neq k}^\top \mathcal{R}_t.
\end{align*}
Finally, once all multiplier updates are computed, the update for the primal variable $\Delta X_k$ is retrieved as
\begin{align*}
    \mathrm{vec}(\Delta X_k) = -(\mathcal{W}_{\neq k}^\top \mathbfcal{L}_\mathcal{Z} \mathcal{W}_{\neq k})^{-1} (\mathcal{W}_{\neq k}^\top \mathcal{R}_c + \mathcal{W}_{\neq k}^\top \mathbfcal{L}_\mathcal{X} \mathcal{W}_{\neq k} \mathrm{vec}(\Delta {Z}_k)).
\end{align*}
Notably, we avoid eliminating $\vecTT(\Delta \mathcal{Z})$ prior to projection, as doing so would introduce matrices $\mathbfcal{L}_\mathcal{X} \mathbfcal{A}_\textnormal{eq}^\top$ and $\mathbfcal{L}_\mathcal{X} \mathbfcal{A}_\textnormal{ineq}^\top$, which lead to prohibitively high TT ranks and risk compromising the low-rank structure of the solution. In the absence of inequalities, the system matrix simplifies accordingly, only requiring one Cholesky and LU decomposition.

Note that this results in a cost per AMEn-sweep of $\mathcal{O}(d r_\mathrm{max}^6 n^3 )$ with $r_\mathrm{max} = \max_k \{ r_{\mathcal{W}_k} \}_{k=2}^d$ compared to the cost of a regular IPM's KKT system solve of $\mathcal{O}(n^{3d})$. The computational cost can be further reduced by applying the matrix-free approach to solve the projected system. Each matrix-block multiplication with a vector reduces to three matrix products: two with cost $\mathcal{O}(r_\mathrm{max}^3 R_\mathrm{max} n)$ and one with cost $\mathcal{O}(r_\mathrm{max}^2 R_\mathrm{max}^2 n^2)$, where $R_\mathrm{max}$ is the maximum TT rank of the TT matrices in the KKT system. This is significantly cheaper than the original block-multiplication cost of $\mathcal{O}(n^{2d})$.

{

\begin{proposition}[Nonsingularity of the projected KKT system]
\label{prop:projected-kkt-nonsingular}
Let $\mathcal{W}_{\neq k}\in\mathbb{R}^{N\times p}$ be the enriched AMEn frame at
core $k$, with $p$ denoting the local frame dimension, after orthogonalisation, so that $\mathcal{W}_{\neq k}^{\top}\mathcal{W}_{\neq k}=\mathcal{I}$. Assume $\mathcal{X}\succ0$, $\mathcal{Z}\succ 0$, and strict positivity of the active inequality slacks and multipliers. Throughout this proposition, hats denote the corresponding projections with $\mathcal{W}_{\neq k}$. Set
\begin{align*}
    \widehat{\mathbfcal L}_{\mathcal Z}
    &:=
    \mathcal{W}_{\neq k}^{\top}\mathbfcal L_{\mathcal Z}\mathcal{W}_{\neq k},
    &
    \widehat{\mathbfcal L}_{\mathcal X}
    &:=
    \mathcal{W}_{\neq k}^{\top}\mathbfcal L_{\mathcal X}\mathcal{W}_{\neq k},
    &
    \widehat{\mathbfcal M}
    &:=
    \widehat{\mathbfcal L}_{\mathcal X}^{-1/2}
    \widehat{\mathbfcal L}_{\mathcal Z}
    \widehat{\mathbfcal L}_{\mathcal X}^{-1/2}.
\end{align*}
Define
\begin{align*}
    \widehat{\mathbfcal A}
    &:=
    \operatorname{diag}(\mathcal{W}_{\neq k},\mathcal{W}_{\neq k})^\top
    \begin{bmatrix}\mathbfcal A_{\mathrm{eq}}\\ \mathbfcal A_{\mathrm{ineq}}\end{bmatrix}
    \mathcal{W}_{\neq k}, \quad
    \widehat{\mathbfcal B}
    :=
    \operatorname{diag}(\mathcal{W}_{\neq k},\mathcal{W}_{\neq k})^\top
    \begin{bmatrix}\mathbfcal A_{\mathrm{eq}}\\ \mathbfcal T\mathbfcal A_{\mathrm{ineq}}\end{bmatrix}
    \mathcal{W}_{\neq k},\\
    \widehat{\mathbfcal D}
    &:=
    \operatorname{diag}(\mathcal{W}_{\neq k},\mathcal{W}_{\neq k})^\top
    \operatorname{diag}\bigl(\mathbfcal K_{\mathcal Y},
    \mathcal D_{\mathrm{ineq}}(\mathcal X)+\mathbfcal K_{\mathcal T}\bigr) \operatorname{diag}(\mathcal{W}_{\neq k},\mathcal{W}_{\neq k}).
\end{align*}
Suppose that
\begin{itemize}
    \item[(A1)] The projection retains the equality-multiplier kernel separation on the local frame:
    \begin{align*}
        \ker\!\left(
            \left(
                \mathcal{W}_{\neq k}^{\top}
                \mathbfcal A_{\mathrm{eq}}
                \mathcal{W}_{\neq k}
            \right)^{\top}
        \right)
        \cap
        \ker\!\left(
            \mathcal{W}_{\neq k}^{\top}
            \mathbfcal K_{\mathcal Y}
            \mathcal{W}_{\neq k}
        \right)
        =
        \{0\}.
    \end{align*}

    \item[(A2)] The transformed off-diagonal coupling is relatively bounded:
    \begin{align*}
        \left|
        u^\top
        \left(
            \widehat{\mathbfcal A}\widehat{\mathbfcal L}_{\mathcal X}^{1/2}
            -
            \widehat{\mathbfcal B}\widehat{\mathbfcal L}_{\mathcal X}^{-1/2}
        \right)x
        \right|
        \le
        \gamma
        \bigl(u^\top\widehat{\mathbfcal D}u\bigr)^{1/2}
        \bigl(x^\top\widehat{\mathbfcal M}x\bigr)^{1/2}
    \end{align*}
    for some $0\le\gamma<2$ and all $u=[y,t]^\top$ with
    $x,y,t\in\mathbb R^p$.
\end{itemize}
Then the projected AMEn KKT system is nonsingular. \\
\textit{Proof.} The proof is given in Appendix~\ref{app:proof_projected_kkt}.
\end{proposition}
For the benchmark formulations in Section~\ref{sec:conceptualisation}, (A1)
follows from the projected mask structure. Assumption (A2) is a local compatibility condition for the block scalings. It is favoured when the KKT scaling keeps the projected AHO blocks generated by
$\mathcal X$ and $\mathcal Z$ spectrally comparable, so that
$\widehat{\mathbfcal M}$ is well conditioned on the local frame, and when the
active multiplier scaling $\mathbfcal T$ is moderate on the projected inequality
rows. In concrete terms, the active entries of $\mathcal T$ should remain bounded
and should not vary sharply relative to the slack diagonal
$\mathcal D_{\mathrm{ineq}}(\mathcal X)+\mathbfcal K_{\mathcal T}$; otherwise
$\widehat{\mathbfcal B}$ can become a strongly rescaled version of
$\widehat{\mathbfcal A}$. Under such scaling, a bound with some $\gamma<2$ is
natural. In the absence of active inequalities, $\widehat{\mathbfcal B}
=\widehat{\mathbfcal A}$, so only the $\mathcal X$/$\mathcal Z$ scaling remains.
}

{

Let $\mathcal K$ denote the current block KKT operator in \eqref{eq:tt-kkt-newton-system}, and let $\mathcal R :=\begin{bmatrix} \mathcal R_d & \mathcal R_p & \mathcal R_t & \mathcal R_c \end{bmatrix}^\top$.
Let $\widetilde{\mathcal K}$ and $\widetilde{\mathcal R}$ denote the TT-rounded operator and residual used in the AMEn solve, and let $\Delta \mathcal U^{\mathrm{tt}}$ be the computed TT direction. We measure its exact-equation residual by $\mathcal E := \mathcal K\Delta \mathcal U^{\mathrm{tt}}+\mathcal R$. Thus an exact Newton direction has $\mathcal E=0$. The standard inexact-Newton forcing condition requires this residual to be small relative to the current KKT residual:
\begin{align}
    \|\mathcal E\| \leq \eta \|\mathcal R\|,
    \label{eq:forcing-condition}
\end{align}
for some forcing term $\eta \in [0,1)$ \cite{Dembo1982}. The next lemma gives a sufficient rounded-system condition under which the computed TT direction satisfies \eqref{eq:forcing-condition}.

\begin{lemma}[AMEn Directions as Inexact Newton Directions]
\label{lem:tt-inexact-newton-direction}
Let $\mathcal R\neq 0$. Assume that $\mathcal K$ is
nonsingular and $\widetilde{\mathcal R} = \mathrm{round}_\mathrm{TT}(\mathcal R, \delta_R)$. Suppose that the AMEn direction $\Delta \mathcal U^{\mathrm{tt}}$ satisfies the rounded-system
forcing condition
\begin{align*}
    \|\widetilde{\mathcal K} \Delta \mathcal U^{\mathrm{tt}}+\widetilde{\mathcal R}\|
    \leq
    \widetilde{\eta}\|\widetilde{\mathcal R}\|,
\end{align*}
then there exists a rounding tolerance $0\leq \delta_{\mathcal K} \leq \frac{\eta-\delta_R-(1+\delta_R)\widetilde{\eta}}{(1+\delta_R)(1+\widetilde{\eta})}$ with forcing parameter $\delta_R + (1+\delta_R)\widetilde{\eta} \leq \eta < 1$ such that $\widetilde{\mathcal K} = \mathrm{round}_\mathrm{TT}(\mathcal K, \frac{\delta_{\mathcal K}}{(1+\delta_{\mathcal K}) \kappa_F(\mathcal K)})$ and the inexact Newton forcing condition holds:
\begin{align*}
    \frac{\|\mathcal K \Delta \mathcal U^{\mathrm{tt}}+\mathcal R\|}{\|\mathcal R\|} \leq \eta .
\end{align*}
\textit{Proof.} The proof is given in Appendix~\ref{app:proof_tt_inexact_newton}.
\end{lemma}
For the AHO Newton system, the KKT operator is locally uniformly conditioned under the standard regularity assumptions. With strict complementarity and primal-dual nondegeneracy, the limiting AHO Jacobian is nonsingular \cite[Theorem~3.1]{Alizadeth1998}. This means $\kappa_F(\mathcal K) \leq C_\mathrm{AHO}$ is bounded, and a small enough rounding tolerance is sufficient for the forcing condition to hold. In practice, the blocks are TT rounded independently. For a block matrix with four blocks, choosing a rounding tolerance of $\delta/2$ for each block guarantees that the Frobenius-norm rounding error of the full block matrix is at most $\delta$. The arguments below are in the spirit of the analysis for primal-dual interior-point methods \cite{Bellavia1998}.

\begin{theorem}[Local convergence under TT forcing condition]
\label{th:tt-local-convergence}
Let $\mathcal U_k=(\mathcal X_k,\mathcal Y_k,\mathcal T_k,\mathcal Z_k)$,
$\mathcal R_k:=\mathcal R_{\mu_k}(\mathcal U_k)$, and let $\mathcal K_k$ be the
corresponding block KKT operator. Let $\mathcal U^\star$ be a locally regular
KKT point for the AHO direction. Suppose that there is a sequence
$\{\rho_k\}\subset[0,\infty)$ such that:
\begin{itemize}
    \item[(L1)] The iterates remain in the local AHO neighbourhood of
    $\mathcal U^\star$.
    \item[(L2)] The accepted step lengths satisfy
    $\alpha_k\in(0,1]$ and $\alpha_k\geq\underline\alpha$ for some
    $\underline\alpha>0$.
    \item[(L3)] With
    $\widehat{\mathcal U}_{k+1}:=
    \mathcal U_k+\alpha_k\Delta\mathcal U_k^{\mathrm{tt}}$, the rounded and
    spectrally corrected iterate satisfies
    \begin{align*}
        \left\|
        \mathcal R_{\mu_k}(\mathcal U_{k+1})
        -
        \mathcal R_{\mu_k}(\widehat{\mathcal U}_{k+1})
        \right\|
        \leq \rho_k .
    \end{align*}
\end{itemize}
Assume further that, for all sufficiently large $k$, the rounded AMEn direction
satisfies the conditions of \Cref{lem:tt-inexact-newton-direction} with resulting
forcing term $\eta_k\leq\bar\eta<1$.
Then there exist constants $q_0\in(0,1)$ and $C>0$ such that, for all sufficiently large
$k$,
\begin{align*}
    \|\mathcal R_0(\mathcal U_{k+1})\|
    \leq
    q_0\|\mathcal R_0(\mathcal U_k)\|
    +C(\mu_k+\rho_k).
\end{align*}
Consequently,
\begin{align*}
    \limsup_{k\to\infty}\|\mathcal R_0(\mathcal U_k)\|
    \leq
    \frac{C}{1-q_0}
    \limsup_{k\to\infty}(\mu_k+\rho_k).
\end{align*}
In particular, if $\mu_k=O(\rho_{\min})$, $\rho_k=O(\rho_{\min})$, and the
rounded-solve tolerances keep $\eta_k\leq\bar\eta<1$, then the TT IPM converges
locally to an $O(\rho_{\min})$ neighbourhood of the KKT conditions. \\
\textit{Proof.} The proof is given in Appendix~\ref{app:proof_tt_local_convergence}.
\end{theorem}

\begin{remark}
\Cref{th:tt-local-convergence} is a local convergence result for the implemented TT IPM with finite tensor-train tolerances. If the tolerance floor is active, the perturbation term need not vanish and convergence is only to a tolerance-dependent neighbourhood. If instead $\mu_k+\rho_k\to0$, the same estimate gives $\|\mathcal R_0(\mathcal U_k)\|\to0$. Moreover, under the standard additional inexact-Newton assumptions that the forcing terms vanish and full Newton steps are eventually accepted, the usual local superlinear convergence conclusion is recovered.
\end{remark}

The next result gives a global residual bound for a safeguarded accepted-step variant under a finite-window residual control condition.

\begin{theorem}[Global convergence of a safeguarded TT IPM]
\label{th:tt-global-convergence}
Let $\{\mathcal U_k\}$ be a bounded sequence generated by a safeguarded TT IPM, and define $\Psi_k:=\|\mathcal R_{\mu_k}(\mathcal U_k)\|+\mu_k$. Suppose that:
\begin{itemize}
    \item[(G1)] Writing $\mathcal U_k=(\mathcal X_k,\mathcal Y_k,\mathcal T_k,\mathcal Z_k)$, the iterates remain strictly interior: $\mathcal X_k\succ0$, $\mathcal Z_k\succ0$, $\mathcal T_k>0$, and {$s_{\mathrm{ineq},k}>0$}.
    \item[(G2)] Whenever $\mathcal R_{\mu_k}(\mathcal U_k)\neq0$, the rounded AMEn direction satisfies the rounded-system conditions of \Cref{lem:tt-inexact-newton-direction}, with resulting forcing term $\eta_k\leq\bar\eta<1$.
    \item[(G3)] With $\widehat\Psi_{k+1}:=\|\mathcal R_{\mu_{k+1}}(\mathcal U_k+\alpha_k\Delta\mathcal U_k^{\mathrm{tt}})\|+\mu_{k+1}$, the accepted iterate satisfies $|\Psi_{k+1}-\widehat\Psi_{k+1}|\leq\rho_k$.
    \item[(G4)] There exist an integer $M\ge1$ and constants $\tau\in(0,1)$ and $c_\Omega>0$ such that, with
    \begin{align*}
        \Theta_k := \max_{i=k-M,\ldots,k}\Psi_i,
        \qquad
        \Omega_k := \max_{i=k-M,\ldots,k}\rho_i,
    \end{align*}
    the accepted iterates satisfy, for all sufficiently large $k$,
    \begin{align*}
        \Theta_{k+M+1}
        \leq
        (1-\tau)\Theta_k+c_\Omega\Omega_k .
    \end{align*}
\end{itemize}
Then there is a constant $C>0$ such that
\begin{align*}
    \limsup_{k\to\infty}\bigl(\|\mathcal R_0(\mathcal U_k)\|+\mu_k\bigr)
    \leq
    C\limsup_{k\to\infty}\rho_k .
\end{align*}
Thus a fixed perturbation floor $\rho_{\min}$ gives convergence to an $O(\rho_{\min})$ KKT neighbourhood. If $\rho_k\to0$, then $\|\mathcal R_0(\mathcal U_k)\|\to0$ and $\mu_k\to0$, and every accumulation point of $\{\mathcal U_k\}$ is a KKT point. If a tail of the sequence additionally satisfies the assumptions of \Cref{th:tt-local-convergence}, then the local estimate in \Cref{th:tt-local-convergence} applies to that tail. \\
\textit{Proof.} The proof is given in Appendix~\ref{app:proof_tt_global_convergence}.
\end{theorem}
The practical implementation recomputes the KKT residual after each rounded and spectrally corrected update and monitors the finite-window residual measure in (G4). Steps that leave this windowed control regime are damped and retried. Thus, TT rounding, approximate step-size computation, and spectral correction enter the accepted iteration only through the perturbation bound $\rho_k$ used in \Cref{th:tt-global-convergence}.

}

\subsection{Maintaining Positive Definiteness and Optimal Step Length Calculation}

 The iterative updates (see \cref{sec:practical_considerations} for details) are paired with a subsequent tensor-train rounding step. This potentially leads to a violation of the positive definiteness of $\mathcal{X}$, $\mathcal{Z}$, and the strict positivity of $\mathcal{T}$. To prevent this violation from occurring, we need to estimate the perturbation to the (smallest) eigenvalue of those matrices, introduced by the TT rounding.

\begin{corollary}\label{cor:psdness_rounding}
    Suppose we have a symmetric TT matrix $\mathcal{M} \in \mathbb{R}^{n^d \times n^d}$ with a symmetric TT rounded approximation $\widetilde{\mathcal{M}} = \mathrm{round}_\mathrm{TT}(\mathcal{M}, \delta)$
    with the error $\|\mathcal{M} - \widetilde{\mathcal{M}}\|_F = \zeta \le \delta.$
    Then their minimum eigenvalues differ by at most $\zeta$, i.e.,
    \begin{align*}
        |\lambda_{\min}(\mathcal{M}) - \lambda_{\min}(\widetilde{\mathcal{M}})| \leq \zeta.
    \end{align*}
\textit{Proof.} {The proof is given in Appendix~\ref{app:psdness_rounding}.}
\end{corollary}

In principle, one may symmetrise each TT core of the resulting TT matrix with respect to its $i_k,j_k$ indices to enforce exact symmetry. In practice, this is rarely done since it doubles the rank, and the asymmetry error is negligible compared to the truncation error.
Using the result of this corollary, we can alleviate the issue of maintaining positive definiteness. {Let $\delta$ be the desired perturbation budget and let $\mathcal{S}_\mathcal{M}$ denote the stabilising shift direction, with $\mathcal{S}_\mathcal{M}=\mathcal{I}$ for $\mathcal{M}\in\{\mathcal{X},\mathcal{Z}\}$ and $\mathcal{S}_\mathcal{M}=\mathbf{1}_\mathcal{T}$ for $\mathcal{M}=\mathcal{T}$. We first round with tolerance $\bar\delta=\delta/(1+\|\mathcal{S}_\mathcal{M}\|_F)$, obtain the actual truncation error $\zeta=\|\mathcal{M}-\widetilde{\mathcal{M}}\|_F\leq\bar\delta$, and then use the corrected tensor $\widetilde{\mathcal{M}}+\zeta\mathcal{S}_\mathcal{M}$. Hence
\[
    \|\mathcal{M}-(\widetilde{\mathcal{M}}+\zeta\mathcal{S}_\mathcal{M})\|_F
    \leq \zeta(1+\|\mathcal{S}_\mathcal{M}\|_F)\leq \delta,
\]
while the shift restores positive definiteness for $\mathcal{X},\mathcal{Z}$ and positivity on the active inequality mask for $\mathcal{T}$. In the unlikely case that this does not result in a lower TT rank, we abstain from TT rounding.}
{ This spectral correction is not a projection onto the affine constraint set. For instance, in MaxCut it can perturb $\operatorname{diag}(\mathcal{X})=\mathbf{1}$. Since the perturbation is bounded by $\delta$, the induced change in any affine residual is also bounded, e.g., by $\|\mathbfcal{A}_{\mathrm{eq}}\|\,\delta$ for the equality residual. Our method is an infeasible primal-dual IPM: after every rounded or shifted update, the KKT residuals are recomputed, and the next Newton system corrects the resulting primal and dual residuals. Thus the shift is accounted for as a controlled inexactness in the sense of \Cref{lem:tt-inexact-newton-direction}. Exact affine feasibility is not required at intermediate iterates, while final feasibility is enforced by the stopping criteria.}

Tensor-train eigenvalue problems \cite{Dolgov2014} are sparse, extreme eigenvalue problems. {The step length is determined by the requirement $\mathcal{M}+\alpha\Delta\mathcal{M}\succeq0$. For every vector $\mathcal{E}$ with $\langle \mathcal{E},\Delta\mathcal{M}\mathcal{E}\rangle<0$, admissibility requires $\alpha\leq-\langle \mathcal{E},\mathcal{M}\mathcal{E}\rangle/\langle \mathcal{E},\Delta\mathcal{M}\mathcal{E}\rangle$. Thus the limiting step is an extreme generalised Rayleigh quotient.} We can pose the optimal step length calculation as a generalised eigenvalue problem adapted to the TT format by solving $d$ generalised eigenvalue problems with eigenvector $\mathcal{E} = \mathrm{vec}([\![E_1, ..., E_d ]\!])$ for given symmetric matrices $\mathcal{M} \succeq 0, \Delta \mathcal{M} \in \mathbb{R}^{n^d \times n^d}$. In the terminology of ALS, this means we minimise $\mathscr{J}(\mathcal{E}) := -\frac{\langle \mathcal{E},  \mathcal{M} \mathcal{E}\rangle}{\langle \mathcal{E}, \Delta \mathcal{M} \mathcal{E}\rangle}$ via the local eigenvalue problem
\[
\min_{E_k} \mathrm{vec}(E_k)^\top \mathcal{E}_{\neq k}^\top \Delta \mathcal{M} \mathcal{E}_{\neq k} \mathrm{vec}(E_k) \quad \text{s.t.} \quad \mathrm{vec}(E_k)^\top \mathcal{E}_{\neq k}^\top \mathcal{M} \mathcal{E}_{\neq k} \mathrm{vec}(E_k) = 1.
\]
We solve the analogous problems for $\mathcal{M} = \mathcal{X}, \mathcal{Z}$.
If $\mathcal{M}$ is not positive semidefinite, we can usually detect this even without accurate frame matrices and raise the appropriate errors.
A sparse eigenvalue solver with a Cholesky-based reduction finds the largest eigenvalue { of the projected system}, and the primal and dual step sizes are computed as
\begin{align*}
\alpha_p &= \left\{\begin{array}{ll}\min\left(1, -\frac{\mathcal{E}_{\mathcal{X} + \Delta \mathcal{X}}^\top \mathcal{X} \mathcal{E}_{\mathcal{X} + \Delta \mathcal{X}}}{\mathcal{E}_{\mathcal{X} + \Delta \mathcal{X}}^\top \Delta \mathcal{X} \mathcal{E}_{\mathcal{X} + \Delta \mathcal{X}}}\right) & \Delta \mathcal{X} \prec 0, \\ 1 & \Delta \mathcal{X} \succeq 0, \end{array}\right. \\ \alpha_d &= \left\{\begin{array}{ll}\min\left(1, -\frac{\mathcal{E}_{\mathcal{Z} + \Delta \mathcal{Z}}^\top \mathcal{Z} \mathcal{E}_{\mathcal{Z} + \Delta \mathcal{Z}}}{\mathcal{E}_{\mathcal{Z} + \Delta \mathcal{Z}}^\top \Delta \mathcal{Z} \mathcal{E}_{\mathcal{Z} + \Delta \mathcal{Z}}}\right) & \Delta \mathcal{Z} \prec 0, \\ 1 & \Delta \mathcal{Z} \succeq 0.\end{array}\right.
\end{align*}

{ For this eigenvalue subproblem, we use a custom two-site DMRG-type TT eigensolver.} This yields a per-DMRG sweep complexity of $\mathcal{O}\!\left(d \, r_\mathrm{max}^6 \, n^6 \right), \quad r_\mathrm{max} = \max_{k=2,\dots,d} r_{\mathcal{E}_k}$, as compared to the $\mathcal{O}(n^{3d})$ cost of the Cholesky factorisation in a standard IPM eigenproblem solve. The larger exponent of $n$ in the DMRG case reflects the fact that the algorithm operates on pairs of neighbouring cores simultaneously to achieve rank-adaptation.
However, it can also be implemented in a matrix-free fashion by using iterative eigenvalue solvers (such as LOBPCG \cite{Knyazev2001}) with $\mathcal{O}(d R_{\max}^4 n^2 + d R_{\max}^3 n^3)$ complexity, where $R_{\max}$ is the maximum TT rank over $\mathcal{M}$ and $\Delta \mathcal{M}$.

{The IPM also ensures that $s_{\mathrm{ineq}}(\mathcal X+\alpha_p\Delta\mathcal X)\ge0$ and $\mathcal{T} + \alpha_d \Delta \mathcal{T} \geq 0$. This is verified by checking the minimum eigenvalues of the diagonal matrices $\diagTT(s_{\mathrm{ineq}}(\mathcal X+\alpha_p\Delta\mathcal X))$ and $\diagTT(\vecTT(\mathcal{T} + \alpha_d \Delta \mathcal{T}))$. Since these operators are diagonal, their smallest eigenvalues are exactly the minimum entries of the slack and multiplier vectors. The diagonal system stays block diagonal when projected, so the projected eigenproblem is substantially less computationally demanding.} The corresponding eigenvectors $\mathcal{E}_{\mathcal{X} + \alpha_p \Delta \mathcal{X}}$ and $\mathcal{E}_{\mathcal{T} + \alpha_d \Delta \mathcal{T}}$ are usually TT rank one. {If the denominator in either ratio is nonnegative, the corresponding correction factor is set to one.} They give the step lengths
\begin{align*}
\alpha'_p &= \alpha_p \cdot \min \Big(1, - \frac{\langle \mathcal{E}_{\mathcal{X} + \alpha_p\Delta \mathcal{X}}, \mathcal{X} \mathcal{E}_{\mathcal{X} + \alpha_p\Delta \mathcal{X}} \rangle}{\langle \mathcal{E}_{\mathcal{X} + \alpha_p\Delta \mathcal{X}}, \Delta \mathcal{X} \mathcal{E}_{\mathcal{X} + \alpha_p\Delta \mathcal{X}} \rangle} \Big), \\
\alpha'_d &= \alpha_d \cdot \min \Big(1, -\frac{\langle \mathcal{E}_{\mathcal{T} + \alpha_d \Delta \mathcal{T}}, \mathcal{T} \mathcal{E}_{\mathcal{T} + \alpha_d  \Delta \mathcal{T}} \rangle}{\langle \mathcal{E}_{\mathcal{T} + \alpha_d  \Delta \mathcal{T}}, \Delta \mathcal{T} \mathcal{E}_{\mathcal{T} + \alpha_d \Delta \mathcal{T}} \rangle} \Big).
\end{align*}

We enhance the DMRG sweeps with random kick ranks. When inequalities are involved in the problem formulation, there are several variables in $\mathbfcal{A}_{\text{ineq}}\vecTT(\mathcal{X} + \alpha_p \Delta \mathcal{X})$ that are mapped to zero as they are not involved in any inequalities. Similarly, $\mathcal{T} + \alpha_d \Delta \mathcal{T}$ will have only nonzero entries on $\mathbf{1}_\mathcal{T}$. This leads to numerical issues as it may prune $\alpha_d$ based on extremely small value entries beyond the precision of tensor-trains. Should the respective entry fall outside the precision of the tensor-trains, i.e., $\langle \mathcal{E}_{\mathcal{X} + \alpha_p\Delta \mathcal{X}}, (\mathcal{X} + \alpha_p \Delta \mathcal{X}) \mathcal{E}_{\mathcal{X} + \alpha_p\Delta \mathcal{X}} \rangle$ and $\langle \mathcal{E}_{\mathcal{T} + \alpha_d \Delta \mathcal{T}}, (\mathcal{T} + \alpha_d  \Delta \mathcal{T}) \mathcal{E}_{\mathcal{T} + \alpha_d  \Delta \mathcal{T}} \rangle$ fall below the numerical precision of the tensor-train, we know that it is safe to take the whole step $\alpha_p, \alpha_d$ without a reduction, as the minimum entry will be one of the uninvolved zero entries. We warm-start all eigen-problems with the computed eigenvector from the previous IPM iteration.

\subsection{Practical Considerations}\label{sec:practical_considerations}

{
\paragraph{Combined Multiplier Storage}
The algorithm generally uses separate multiplier blocks $\mathcal{Y}$ and $\mathcal{T}$.
Let $\mathbfcal P_{\mathcal Y}:=\mathcal I-\mathbfcal K_{\mathcal Y}$ and
$\mathbfcal P_{\mathcal T}:=\mathcal I-\mathbfcal K_{\mathcal T}$ denote the
corresponding non-idle multiplier masks. In the implementation, if these masks are
disjoint, i.e., if $\mathbfcal P_{\mathcal Y}\mathbfcal P_{\mathcal T}=0$, the two updates can
be stored in one padded TT block $\Delta\Lambda$ and recovered by
\begin{align*}
    \Delta\mathcal{Y}=\mathbfcal P_{\mathcal Y}\Delta\Lambda,\qquad
    \Delta\mathcal{T}=\mathbfcal P_{\mathcal T}\Delta\Lambda .
\end{align*}
Conversely, on the non-idle supports, $\Delta\Lambda=\Delta\mathcal{Y}+\Delta\mathcal{T}$.
Thus the combined multiplier block is only a masked coordinate compression of the
separate $(\Delta\mathcal{Y},\Delta\mathcal{T})$ variables.

Since $\mathbfcal P_{\mathcal Y}\mathbfcal P_{\mathcal T}=0$, the map
$\Delta\Lambda\mapsto(\mathbfcal P_{\mathcal Y}\Delta\Lambda,
\mathbfcal P_{\mathcal T}\Delta\Lambda)$ is a bijection between the combined
non-idle support and the direct sum of the separate non-idle supports, with inverse
$(\Delta\mathcal{Y},\Delta\mathcal{T})\mapsto\Delta\mathcal{Y}+\Delta\mathcal{T}$.
Substituting this parametrisation in the separate KKT system changes only the
multiplier blocks. Suppressing vectorisation, the separate multiplier column
\begin{align*}
    -\mathbfcal A_{\mathrm{eq}}^\top\Delta\mathcal{Y}
    -
    \mathbfcal A_{\mathrm{ineq}}^\top\Delta\mathcal{T}
\end{align*}
becomes
\begin{align*}
    -\bigl(
        \mathbfcal A_{\mathrm{eq}}^\top\mathbfcal P_{\mathcal Y}
        +
        \mathbfcal A_{\mathrm{ineq}}^\top\mathbfcal P_{\mathcal T}
    \bigr)\Delta\Lambda .
\end{align*}
The two separate multiplier rows are replaced by the masked sum
\begin{align*}
    -\bigl(
        \mathbfcal P_{\mathcal Y}\mathbfcal A_{\mathrm{eq}}
        +
        \mathbfcal P_{\mathcal T}\mathbfcal T\mathbfcal A_{\mathrm{ineq}}
    \bigr)\Delta\mathcal X
    +
    \bigl[
        \mathbfcal P_{\mathcal Y}\mathbfcal K_{\mathcal Y}\mathbfcal P_{\mathcal Y}
        +
        \mathbfcal P_{\mathcal T}
        \bigl(\mathcal D_{\mathrm{ineq}}(\mathcal X)+\mathbfcal K_{\mathcal T}\bigr)
        \mathbfcal P_{\mathcal T}
    \bigr]\Delta\Lambda ,
\end{align*}
with right-hand side
$\mathbfcal P_{\mathcal Y}\mathcal R_p+\mathbfcal P_{\mathcal T}\mathcal R_t$.
The $\Delta\mathcal X$- and $\Delta\mathcal Z$-equations are otherwise unchanged.

After AMEn projection, a combined multiplier frame $\mathcal{W}_\Lambda$ induces
the separate masked frames $\mathbfcal P_{\mathcal Y}\mathcal{W}_\Lambda$ and
$\mathbfcal P_{\mathcal T}\mathcal{W}_\Lambda$; these are not independently
chosen $\mathcal{Y}$- and $\mathcal{T}$-frames. Hence the projected combined
matrix is exactly the projected separate matrix assembled on this induced masked
subspace, written in the $\Delta\Lambda$ coordinates. The nonsingularity
assumptions are therefore checked on the same finite-dimensional projected
operator in either representation. Accordingly, assumptions (A1)--(A2) in
\Cref{prop:projected-kkt-nonsingular} are imposed on the projected blocks of the
operator actually solved. In the combined branch this means the compressed
$\Delta\Lambda$ formulation, or equivalently the separate formulation restricted
to the induced masked frames $\mathbfcal P_{\mathcal Y}\mathcal{W}_\Lambda$ and
$\mathbfcal P_{\mathcal T}\mathcal{W}_\Lambda$. The proof therefore applies
without modification. The combined storage is only a small performance
optimisation that reduces repeated multiplier contractions.
}

\paragraph{Error Tracking and TT Rounding Tolerance}
Both SDPA and SCS utilise the infinity norm of the primal and dual residuals as part of their stopping criteria. However, in the context of tensor-train representations, the evaluation of the infinity norm is computationally infeasible. Consequently, employing the infinity norm for error monitoring is impractical. Thus, we track the relative errors and employ the criteria
\begin{align*}
    \frac{||\mathcal{R}_p||_F}{1 + ||\mathcal{B}_\text{eq}||_F} \leq \epsilon_p, \quad \frac{||\mathcal{R}_d||_F}{1 + ||\mathcal{C}||_F} \leq \epsilon_d,
\end{align*}
\begin{align*}
    \frac{\langle \mathcal{X}, \mathcal{Z} \rangle + \langle \vecTT(\mathcal{T}), {s_{\mathrm{ineq}}(\mathcal X)} \rangle}{(2^d + q)(1+\langle \mathcal{C}, \mathcal{X} \rangle)} \leq \epsilon_c.
\end{align*}
The thresholds $\epsilon_p$, $\epsilon_d$, and $\epsilon_c$ are typically chosen to be in the interval $[10^{-4}, 10^{-3}]$.

When a desired relative error threshold is achieved, we set the respective residual to zero in the KKT system. As we do not contract over zero-blocks, this reduces the computational load, also because the associated residual TT rank of a zero-block is one and is more likely to yield lower TT rank solutions.
{
Before each AMEn solve, the KKT block system is row-scaled using the current feasibility and complementarity residual magnitudes. These scalings are internal to the linear solve. All residuals and stopping criteria are recomputed in the original variables.
}

The error tolerance for the AMEn-solver is reduced in accordance with inexact IPMs via $\eta \leftarrow \max \{\min\{\eta, C \mu \}, \delta \}$, where we set $C=2, \delta = 10^{-4}$. TT rounding is carried out at three stages of the IPM. The primal variables and dual variables of \cref{algline:update_vars_1,algline:update_vars_2} are rounded at $\eta \cdot (1 + ||\mathcal{B}_\text{eq}||_F)/10$ and $\eta \cdot (1 + ||\mathcal{C}||_F)/10$, respectively. The operators $\mathbfcal{L}_\mathcal{X}, \mathbfcal{L}_\mathcal{Z}$ at \cref{algline:kkt_solve} are TT rounded at the same respective errors. The primal, dual, and centrality residuals of \cref{algline:residuals} are accordingly rounded at $\eta \cdot (1 + ||\mathcal{B}_\text{eq}||_F)/100$, $\eta \cdot (1 + ||\mathcal{C}||_F)/100$, and $\eta \cdot (1 + |\langle \mathcal{C}, \mathcal{X} \rangle|)/100$. { The theory of the AHO direction does not prescribe a separate $\mu$-dependent tightening rate for the TT rounding of the KKT operator. A fixed sufficiently small relative rounding tolerance would be sufficient. In the implementation, however, we tie the KKT operator rounding tolerance to the current forcing term $\eta$. This keeps operator-compression errors subordinate to the current inexact Newton tolerance while $\eta$ is decreased, and it becomes a fixed tolerance once the prescribed floor $\rho_{\min}$ is reached.}

\paragraph{Relaxing the Positive Definiteness Barrier}
After the relative stopping criteria are satisfied, we relax the positive definiteness constraint to $\mathcal{M} + \nu \mathcal{I} \succeq 0$ for some small $\nu > 0$, allowing slight negative eigenvalues in $\mathcal{M}$. This relaxation is performed gradually up to $\nu < 10^{-3}$ if necessary until the absolute tolerance criteria on the residual and duality gap are met, enabling convergence towards a potentially rank-deficient solution that better approximates the true optimum. At the very least, we aim to achieve absolute errors of $||\mathcal{R}_p||_F, ||\mathcal{R}_d||_F, \langle \mathcal{X}, \mathcal{Z} \rangle < 10^{-3}$, at which point the IPM terminates.
{
During this finishing phase, damped updates and postprocessing steps are accepted only if they reduce a merit quantity involving the final gap and primal/dual residuals without increasing the residuals beyond their current admissible limits. Otherwise the previous iterate is retained or the direction is recomputed.
}

\paragraph{Mehrotra's Predictor–Corrector Method}
We employ Mehrotra's predictor–corrector method \cite{Mehrotra1992}, a widely used enhancement of primal-dual interior-point methods. The method consists of two phases per iteration: a predictor step, which estimates the affine-scaling direction by ignoring the centrality conditions, and a corrector step, which refines this direction by reintroducing the complementarity conditions to improve convergence towards the central path. This approach typically leads to faster convergence and greater robustness, especially in the presence of ill-conditioning or poorly scaled problems. We warm-start all AMEn-solves with the computed solution from the previous IPM iteration.
{
The first iterations use the cheaper $XZ$-direction, which often admits lower TT-rank search directions. After this warm-up phase, the method switches to the AHO direction, but if the AHO AMEn solve cannot reach the required residual target within the prescribed rank and sweep budget, or if the resulting admissible step is nearly zero, the implementation temporarily retries the Newton solve with the $XZ$-direction. This fallback addresses the practical case where the main obstruction is not the IPM model but the TT rank required to represent the AHO direction at acceptable cost.
}

\section{Problem Formulations in TT format} \label{sec:conceptualisation}
Although one could construct the constraint operators in their standard matrix form and subsequently apply a TT SVD \cite{Oseledets2011} to obtain their tensor-train representation, this approach is often computationally inefficient. Constructing the full matrix explicitly requires extensive nested loops and results in significant memory overhead, particularly for large-scale problems. In contrast, designing the operators directly in TT format circumvents the need for materialising dense matrices. This direct formulation also ensures that the resulting representation is rank-aware by design.

To assemble the problems under consideration, we make use of the following matrices:
\begin{itemize}
    \item $\mathcal{I}$: the identity matrix. This matrix admits a tensor-train representation with TT ranks equal to 1. When $n$ is a power of 2, it can be expressed as $\mathcal{I} = I_2 \otimes \cdots \otimes I_2$, where $I_2$ denotes the $2 \times 2$ identity matrix.
    \item $\mathcal{J}$: the all-ones matrix. This matrix also admits a rank-1 TT representation, i.e.,
    $\mathcal{J} = J_2 \otimes \cdots \otimes J_2$,
    where $J_2$ is the $2 \times 2$ matrix with all entries equal to 1.
    \item $\mathcal{E}_{(\alpha; \beta)}$: the diagonal matrix with a single nonzero entry equal to 1 at the $(\alpha; \beta)$-th position. This matrix is likewise representable as a rank-1 TT matrix $\mathcal{E}_{(\alpha; \beta)} = E_{\alpha_1 \beta_1} \otimes \cdots \otimes E_{\alpha_d \beta_d}$.
\end{itemize}

\subsection{Maximum Cut}
Let $A_G$ denote the adjacency matrix of the graph $G$ under consideration, excluding self-connections, and let $\mathcal{A}_G \in \mathbb{R}^{n^d \times n^d}$ be its representation in TT format. The corresponding objective is then computed via the graph Laplacian as
\begin{align*}
    {\mathcal{C} = \frac{1}{4}\left(\diagTT(\mathcal{A}_G \cdot \mathbf{1}) - \mathcal{A}_G\right).}
\end{align*}
This also means that we have $\rankTT(\mathcal{C}) \leq 2 \cdot \rankTT(\mathcal{A}_G)$. {As in the standard MaxCut SDP relaxation, the nonconvex constraint $\textnormal{rank}(X)=1$ is dropped, while the semidefinite constraint $X\succeq0$ is retained.} The condition $\textnormal{diag}(X) = \mathbf{1}$, i.e., $X_{ii} = 1$, is enforced using the operator and right-hand side
\begin{align*}
    \mathbfcal{A}_\textnormal{eq} := \diagTT(\vecTT(\mathcal{I})), \quad \mathcal{B}_\textnormal{eq} := \vecTT(\mathcal{I}).
\end{align*}
The placement of non-idle Lagrangian multipliers in $\mathcal{Y}$ can be chosen freely. The most straightforward choice is to place them along the diagonal, realised as
\begin{align*}
    \mathbfcal{K}_\mathcal{Y} := \diagTT(\vecTT(\mathcal{J} - \mathcal{I})).
\end{align*}
This operator $\mathbfcal{K}_\mathcal{Y}$ acts by zeroing out the diagonal of a matrix while preserving all off-diagonal entries. By selecting the diagonal of $\mathcal{Y}$, the operator $\mathbfcal{A}_\textnormal{eq}$ becomes self-adjoint.

The TT ranks of the associated operators are given by $\rankTT(\mathbfcal{A}_\textnormal{eq}) = \rankTT(\mathbfcal{A}_\textnormal{eq}^\top) = 1$. For the right-hand side $\mathcal{B}$, we similarly have $\rankTT(\mathcal{B}_\textnormal{eq}) = \mathbf{1}$, while the rank of the masking operator is $\rankTT(\mathbfcal{K}_\mathcal{Y}) = \mathbf{2}$.

{
Let $y$ be a local coefficient vector satisfying the two projected kernel
equations in (A1). Testing these equations with $y$ gives
\begin{align*}
    y^\top\mathcal{W}_{\neq k}^{\top}\mathbfcal K_{\mathcal Y}
    \mathcal{W}_{\neq k}y
    =
    y^\top\mathcal{W}_{\neq k}^{\top}\mathbfcal A_{\mathrm{eq}}^\top
    \mathcal{W}_{\neq k}y
    =
    0.
\end{align*}
Since $\mathbfcal A_{\mathrm{eq}}$ and $\mathbfcal K_{\mathcal Y}$ are
complementary positive semidefinite diagonal projectors, this implies
$\mathcal{W}_{\neq k}y=0$ and hence $y=0$. Thus (A1) holds for the projected
MaxCut equality blocks.
}

\subsection{Maximum Stable Set}
Define $\mathcal{A}_G$ as before. The Maximum Stable Set problem imposes two types of constraints. The first constraint enforces a zero structure on the entries of $\mathcal{X}$ corresponding to the adjacency pattern of the graph. It is realised through
\begin{align*}
    \mathbfcal{A}_\text{part 1} := \diagTT(\vecTT(\mathcal{A}_G)), \quad \mathcal{B}_\text{part 1} := \mathbf{0},
\end{align*}
with $\rankTT(\mathbfcal{A}_\text{part 1}) = \rankTT(\mathcal{A}_G)$. We place the idle Lagrangian multipliers in $\mathcal{Y}$ at the locations of the zero entries in $\mathcal{A}_G$, which ensures that the operator is self-adjoint.

The second constraint enforces a fixed trace on $\mathcal{X}$. Storing the value of this trace in the first tensor entry, we can realise this constraint via
\begin{align*}
    \mathbfcal{A}_\text{part 2} := \vecTT(\mathcal{E}_{(1, \ldots, 1;1, \ldots, 1)})\vecTT(\mathcal{I})^\top, \quad \mathcal{B}_\text{part 2} := \vecTT(\mathcal{E}_{(1, \ldots, 1;1, \ldots, 1)})
\end{align*}
which have $\rankTT(\mathbfcal{A}_\text{part 2}) = \mathbf{1}$ and $\rankTT(\mathcal{B}) = \mathbf{1}$.
Combining both constraints, we obtain the overall operator $\mathbfcal{A}_\textnormal{eq} := \mathbfcal{A}_\text{part 1} + \mathbfcal{A}_\text{part 2}$ and $\mathcal{B}_\textnormal{eq} := \mathcal{B}_\text{part 1} + \mathcal{B}_\text{part 2}$, respectively, with rank bounded by
\begin{align*}
    \rankTT(\mathbfcal{A}_\textnormal{eq}) \leq \rankTT(\mathcal{A}_G) + \mathbf{1}, \quad \rankTT(\mathcal{B}_\textnormal{eq}) = \mathbf{1}.
\end{align*}
Finally, we define the masking operator that filters out idle multipliers in the dual space as
\begin{align*}
    \mathbfcal{K}_\mathcal{Y} := \diagTT(\vecTT(\mathcal{J} - \mathcal{A}_G - \mathcal{B}_\textnormal{eq})),
\end{align*}
which results in $\rankTT(\mathbfcal{K}_\mathcal{Y}) \leq \rankTT(\mathcal{A}_G) + \mathbf{2}$.

{
To verify (A1), use
$\mathbfcal A_{\mathrm{eq}}=\mathbfcal A_\text{part 1}+\mathbfcal A_\text{part 2}$
with
$\mathbfcal A_\text{part 2}=\mathcal B_\text{part 2}\vecTT(\mathcal I)^\top$.
Thus
$\mathbfcal A_{\mathrm{eq}}^\top=\mathbfcal A_\text{part 1}
+\vecTT(\mathcal I)\mathcal B_\text{part 2}^\top$. Let $y$ satisfy the two
projected kernel equations in (A1). Testing the projected
$\mathbfcal K_{\mathcal Y}$-equation with $y$ gives
$y^\top\mathcal{W}_{\neq k}^{\top}\mathbfcal K_{\mathcal Y}\mathcal{W}_{\neq k}y=0$.
Since $\mathbfcal K_{\mathcal Y}\succeq0$ is diagonal, this implies
$\mathbfcal K_{\mathcal Y}\mathcal{W}_{\neq k}y=0$, so
$\mathcal{W}_{\neq k}y$ is supported only on the edge mask and the stored trace
coordinate. These two remaining supports are disjoint:
$\mathbfcal A_\text{part 1}$ selects edge coordinates, while
$\mathcal B_\text{part 2}$ selects the stored trace coordinate. Testing the
projected $\mathbfcal A_{\mathrm{eq}}^\top$-equation with $y$ gives
\begin{align*}
    0
    =
    (\mathcal{W}_{\neq k}y)^\top
    \mathbfcal A_{\mathrm{eq}}^\top
    (\mathcal{W}_{\neq k}y)
    =
    \|(\mathbfcal A_\text{part 1})^{1/2}\mathcal{W}_{\neq k}y\|^2
    +
    |\mathcal B_\text{part 2}^\top\mathcal{W}_{\neq k}y|^2,
\end{align*}
because $\mathcal A_G$ has nonnegative elements, in particular zero diagonal, whilst the remaining trace support is
selected by $\mathcal B_\text{part 2}$. The two terms therefore vanish only if
the edge support and the stored trace coordinate of $\mathcal{W}_{\neq k}y$
vanish. Hence $\mathcal{W}_{\neq k}y=0$ and
then $y=0$. Thus (A1) holds for the projected Maximum Stable Set equality
blocks.
}

\subsection{Correlation Clustering}
As introduced in the previous sections, we denote by $\mathcal{A}_G$ the TT adjacency matrix, together with the corresponding weights $\mathcal{C}^+$ and $\mathcal{C}^-$ expressed in TT format. We create $\mathcal{C}$ as
\begin{align*}
    \mathcal{C} = \mathcal{C}^+ + \diagTT(\mathcal{C}^- \cdot \mathbf{1}) - \mathcal{C}^-.
\end{align*}
This means that we have $\rankTT(\mathcal{C}) \leq \rankTT(\mathcal{C}^+) + 2 \cdot \rankTT(\mathcal{C}^-)$. The condition $\textnormal{diag}(X) = \mathbf{1}$, i.e., $X_{ii} = 1$, is again enforced using the self-adjoint operator and right-hand side:
\begin{align*}
    \mathbfcal{A}_\textnormal{eq} := \diagTT(\vecTT(\mathcal{I})), \quad \mathcal{B}_\textnormal{eq} := \vecTT(\mathcal{I}).
\end{align*}
We place the idle multipliers as we did in MaxCut with:
\begin{align*}
    \mathbfcal{K}_\mathcal{Y} := \diagTT(\vecTT(\mathcal{J} - \mathcal{I})).
\end{align*}
The inequality constraint operator is the same as $\mathbfcal{A}_\text{part 1}$ from the Maximum Stable Set problem, i.e.,
\begin{align*}
    \mathbfcal{A}_\textnormal{ineq} := \diagTT(\vecTT(\mathcal{A}_G)), \quad \mathcal{B}_\textnormal{ineq} := \mathbf{0}.
\end{align*}
{
For Correlation Clustering, this positive masked operator is used with the lower-bound convention
\[
    \mathbfcal{A}_\textnormal{ineq}\vecTT(\mathcal{X}) \geq \mathcal{B}_\textnormal{ineq},
    \qquad
    s_\mathrm{ineq}(\mathcal{X})
    = \mathbfcal{A}_\textnormal{ineq}\vecTT(\mathcal{X})-\mathcal{B}_\textnormal{ineq}
    = \vecTT(\mathcal{A}_G\odot\mathcal{X}).
\]
Accordingly, the inequality terms use $+\mathbfcal{A}_\textnormal{ineq}^{\top}\vecTT(\mathcal{T})$ in stationarity and $+\mathbfcal{T}\mathbfcal{A}_\textnormal{ineq}$ in the complementarity linearisation, with the same signs propagated through the projected systems.
}
This ultimately means that we simply place the associated multipliers with
\begin{align*}
    \mathbfcal{K}_\mathcal{T} := \diagTT(\vecTT(\mathcal{J} - \mathcal{A}_G)).
\end{align*}
The TT ranks of the operators are given by $\rankTT(\mathbfcal{A}_\textnormal{eq}) = \mathbf{1}$. For the right-hand side $\mathcal{B}_\textnormal{eq}$, we similarly have $\rankTT(\mathcal{B}_\textnormal{eq}) = \mathbf{1}$, while the rank of the masking operator is $\rankTT(\mathbfcal{K}_\mathcal{Y}) = \mathbf{2}$. It is easy to see then that $\rankTT(\mathbfcal{A}_\textnormal{ineq}) \leq \rankTT(\mathcal{A}_G)$ and $\rankTT(\mathbfcal{K}_\mathcal{T}) \leq 1 + \rankTT(\mathcal{A}_G)$.

{
The equality mask is the same complementary diagonal projector pair as in
MaxCut. Therefore, the same projected argument proves (A1) for the Correlation
Clustering equality blocks. The additional inequality blocks enter (A2), not the
equality-multiplier kernel condition (A1).
}

\section{Numerical Experiments}\label{sec:comp_study}
We compare our TT IPM against SDPA, a high-precision interior-point method, and SCS, a first-order solver based on the alternating direction method of multipliers. SDPA is widely regarded as a state-of-the-art solver for small- to medium-scale SDPs, but it suffers from scalability limitations stemming from the cubic per-iteration cost characteristic of IPMs. SCS excels at handling large-scale problems due to its low per-iteration complexity and modest memory requirements, although it generally provides lower solution accuracy and slower convergence near optimality. We also tested SketchyCGAL, with sketch size set to $\lceil \sqrt{2(m+1)} \rceil$, in accordance with the Barvinok-Pataki rank bound \cite{Pataki1998}. However, on our problem instances beyond size $2^5$, both CGAL and SketchyCGAL repeatedly failed to close the surrogate optimality gap, which is also not directly comparable to the true duality gap. Consequently, we excluded them from the baseline comparisons. { For completeness, we also evaluated a more traditional low-rank baseline based on the Burer--Monteiro factorisation \cite{Boumal2016,pymanopt2016}, with the factorisation rank set according to the Barvinok-Pataki bound. It should be noted that this approach does not provide dual optimality certificates and only applies to MaxCut.}

\subsection{Experimental Setup}
All reported runtimes are wall-clock times measured on the same hardware: an Intel Core i7-1260P CPU with 32\,GB of RAM. Our method is implemented in Python, while both baseline solvers, SCS and SDPA, are compiled, highly parallelised C/C++ implementations linked against optimised numerical libraries. The baselines are accessed via { their respective official Python interfaces \cite{scs_python, sdpa_python}}. As such, raw timing comparisons do not reflect algorithmic efficiency alone and should be interpreted with care. We do not report CPU time, as the Python interpreter overhead and garbage collection distort CPU-based measurements.

All comparisons use SDPA configured to exploit sparsity in the constraint matrices, which improves performance on large-scale problems. {We run SCS with its matrix-free linear-system solver to prevent it from running out of memory early on.}

A solver is marked as out-of-time (OOT) if a single run exceeds 5 hours, and as out-of-memory (OOM) if memory usage exceeds the system's capacity. To ensure consistent performance and prevent thread oversubscription, we fixed the number of threads to 16. {Since runtimes and memory usage are most sensitive to seed-dependent tail behaviour, we report their median, interquartile range, and maximum. Iteration counts and final accuracy measures are reported as arithmetic means over completed seed runs.}

All errors are measured using the squared Frobenius norm to maintain consistency with the duality gap, which is expressed as an inner product. This choice aligns the scale of feasibility errors with that of the duality gap, making comparability easier.

In all experiments, the TT IPM is initialised with three warm-up iterations employing the $XZ$-direction, after which we switch to the AHO-direction. {The TT IPM relative stopping tolerances are derived from $\varepsilon_{\mathrm{gap}}$: $\epsilon_p=\epsilon_d=2\varepsilon_{\mathrm{gap}}$ and $\epsilon_c=\varepsilon_{\mathrm{gap}}/\sqrt{d}$, where $d$ is the TT order.} The AMEn solver tolerance is initially set to $\eta = 10^{-3}$ with $\delta = 10^{-4}$. {For the TT IPM,} the parameter {$\lambda^*$} scales the initial point. The variables are initialised as
\[
\mathcal{X} = \mathcal{Z} = {\lambda^*} \mathcal{I},
\quad
\mathcal{Y} = \mathbf{0}.
\]
The parameter {$\lambda^*$} is fixed at one for MaxCut but is set to two for Maximum Stable Set and Correlation Clustering. When inequality constraints are present, we instead initialise
\[
\mathcal{X} = {\lambda^*} \mathcal{I} + \frac{\alpha_p}{10} \,\mathbfcal{A}_{\mathrm{ineq}}(\mathcal{J}),
\quad
\mathcal{T} = {\lambda^*_{\mathrm{ineq}}} \mathbfcal{A}_{\mathrm{ineq}}(\mathcal{J}),
\]
where $\alpha_p$ is the maximum possible scalar while maintaining positive definiteness. {The variables $\mathcal{Z}$ and $\mathcal{Y}$ remain initialised as above.} The parameter {$\lambda^*_{\mathrm{ineq}}$} is also an initial scaling factor and is set to $10^{-3}$ for the Correlation Clustering problem.

{For SDPA, we set $\epsilon'=\epsilon^*=10^{-5}$. We set the SDPA step length discount factor to $\gamma^*=0.9$ for MaxCut; for Maximum Stable Set and Correlation Clustering, we use the more conservative values $\gamma^*=0.8$ and $0.75$, respectively. For SCS, both the absolute and relative solver tolerances are set to $\varepsilon=10^{-5}/d$, where $d$ is the TT order. These solver tolerances are not directly equivalent: TT IPM is tuned to mimic SDPA's separate feasibility and relative-gap controls, while SCS uses absolute and relative infinity-norm tests for its conic embedding; final accuracy is therefore compared through the reported Frobenius feasibility errors and duality gaps.}

{
For each benchmark, we generate synthetic graph instances whose adjacency tensors have prescribed input TT rank $r_{\mathrm{in}}$, are symmetric, and have zero diagonal. The construction forms binary TT matrices from random $r_{\mathrm{in}}$-dimensional orthonormal bases, enforces the adjacency structure at the level of the TT cores, and rejects samples whose realised TT rank differs from the target. Full details of the sampling procedure are given in Appendix~\ref{app:test_instance_generation}.
}

\subsection{Comparison to Baselines}

{
Across the three benchmark semidefinite relaxations, namely Maximum Cut, Maximum Stable Set, and Correlation Clustering, TT IPM is the only primal-dual SDP solver reported at the largest dimensions, while remaining competitive across several higher-rank instances. Figure~\ref{fig:rank-stratified-certificate-ratios} gives the main compact comparison against the primal-dual baselines SCS and SDPA. It combines certificate quality with runtime or memory usage; values below zero therefore indicate that TT IPM has the smaller certificate-weighted cost. The complete rank-stratified runtime, memory, iteration, feasibility, and duality-gap values, including the Burer--Monteiro reference for MaxCut, are reported in Tables~\ref{tab:comparison_maxcut_ranks}--\ref{tab:comparison_corr_clust_ranks}.
}

\input{figures/rank_stratified_certificate_ratios}

{
The solid SCS curves in Figure~\ref{fig:rank-stratified-certificate-ratios} show that first-order conic solves can be cheaper on some small instances, but this advantage is not uniform once certificate accuracy and memory are included. The dashed SDPA curves fall rapidly because the direct solver is accurate on small instances but becomes unavailable or memory-limited as the dimension grows. TT IPM is therefore most favourable in the regime targeted by the method: larger structured SDPs where a certified primal-dual method should remain available with controlled memory. Burer--Monteiro provides a useful low-rank reference for MaxCut and is very fast, but it solves a different primal-factor problem and does not provide dual feasibility or duality-gap certificates.
}

{
Overall, these results highlight TT IPM as a scalable and memory-efficient solver that preserves appropriate primal-dual accuracy across a wide range of problem sizes. Crucially, among the primal-dual solvers, it remains the most reliable source of SDP certificates at the largest reported dimensions.
}

\newpage
\subsection{Evolution of the TT Rank along the Interior Path}

{
In this section, we investigate the TT ranks of the variable tensors throughout the interior-point iterations. Figures~\ref{fig:rank_evolution_maxcut}--\ref{fig:rank_evolution_corr_clust} show compact rank-evolution strips for each problem. Each strip reports, over normalised optimisation progress, the median over seed runs of the maximum TT rank among the available variable tensors. The log-colour scale therefore encodes higher intermediate rank complexity.
}

{
Across the three problem classes, the TT ranks typically grow during the early or middle iterations and then stabilise or contract near convergence. Increasing the input TT rank generally shifts the strips toward larger scale values, reflecting higher but still structured intermediate TT complexity. The dependence on problem size is less monotone. MaxCut exhibits the most regular progression, whereas Maximum Stable Set and Correlation Clustering show stronger but still localised mid-path rank peaks for the higher-rank instances.
}

\input{figures/maxcut_rank_evolution}
\input{figures/maxstableset_rank_evolution}
\input{figures/corr_clust_rank_evolution}

\subsection{The Influence of Problem Rank on Solver Efficiency}

{
In this section, we examine how the TT rank of the input graph impacts the efficiency and behaviour of our solver. While the size of the problem controls the dimensionality of the variable tensors, the TT rank governs their structural complexity.
}

\begin{figure}[!htbp]
\centering
\begin{tikzpicture}
\begin{axis}[
    name=timeplot,
    ymode=log,
    width=0.47\linewidth,
    height=7cm,
    title={Runtime vs TT Rank},
    xlabel={TT Rank},
    ylabel={Time (s)},
    xtick={1, 2, 3, 4},
    enlarge y limits={value=0.1, upper},
    enlarge x limits=0.2,
    legend style={
        at={(0.02,0.98)},
        anchor=north west,
        font=\small,
        draw=none,
        fill=none,
    },
    label style={font=\small},
    tick label style={font=\small}
]

\addplot[only marks, color=blue, mark=*, mark options={fill=blue, draw=none}, mark size=2pt, opacity=0.7, xshift=-9.0pt, forget plot] coordinates {(1, 0.81) (1, 0.84) (1, 1.35) (1, 1.77) (1, 2.41) (2, 1.70) (2, 1.74) (2, 1.92) (2, 2.14) (2, 3.50) (3, 3.61) (3, 3.91) (3, 4.23) (3, 4.88) (3, 5.61) (4, 7.80) (4, 8.39) (4, 9.19) (4, 11.98) (4, 21.98)};
\addlegendimage{only marks, mark=*, mark options={fill=blue, draw=none}, color=blue}
\addplot[only marks, color=red, mark=square*, mark options={fill=red, draw=none}, mark size=2pt, opacity=0.7, xshift=-3.0pt, forget plot] coordinates {(1, 1.14) (1, 1.23) (1, 1.36) (1, 1.68) (1, 3.97) (2, 2.20) (2, 2.36) (2, 2.56) (2, 3.59) (2, 4.78) (3, 8.10) (3, 8.61) (3, 13.82) (3, 37.34) (3, 37.85) (4, 0.82) (4, 0.86) (4, 7.75) (4, 10.94) (4, 37.93)};
\addlegendimage{only marks, mark=square*, mark options={fill=red, draw=none}, color=red}
\addplot[only marks, color=green!60!black, mark=triangle*, mark options={fill=green!60!black, draw=none}, mark size=3pt, opacity=0.7, xshift=3.0pt, forget plot] coordinates {(1, 1.57) (1, 1.64) (1, 1.86) (1, 1.96) (1, 7.53) (2, 1.67) (2, 1.73) (2, 2.01) (2, 2.23) (2, 10.37) (3, 21.89) (3, 22.49) (3, 36.71) (3, 54.48) (3, 208.10) (4, 10.67) (4, 11.56) (4, 17.34) (4, 30.92) (4, 662.74)};
\addlegendimage{only marks, mark=triangle*, mark size=3pt, mark options={fill=green!60!black, draw=none}, color=green!60!black}
\addplot[only marks, color=orange, mark=diamond*, mark options={fill=orange, draw=none}, mark size=3pt, opacity=0.7, xshift=9.0pt, forget plot] coordinates {(1, 3.34) (1, 3.44) (1, 3.67) (1, 3.90) (1, 4.51) (2, 6.34) (2, 6.50) (2, 10.63) (2, 40.70) (2, 51.53) (3, 8.68) (3, 9.37) (3, 12.97) (3, 47.52) (3, 88.83) (4, 37.26) (4, 39.93) (4, 47.91) (4, 105.80) (4, 3750.51)};
\addlegendimage{only marks, mark=diamond*, mark size=3pt, mark options={fill=orange, draw=none}, color=orange}

\legend{$2^6$, $2^7$, $2^8$, $2^9$}

\end{axis}
\begin{axis}[
    at={(timeplot.south east)},
    anchor=south west,
    xshift=1.55cm,
    ymode=log,
    width=0.47\linewidth,
    height=7cm,
    title={Memory vs TT Rank},
    xlabel={TT Rank},
    ylabel={Memory (MB)},
    xtick={1, 2, 3, 4},
    ytick={1e1, 1e2, 1e3, 1e4},
    ymin=1e1,
    enlarge y limits={value=0.1, upper},
    enlarge x limits=0.2,
    label style={font=\small},
    tick label style={font=\small},
]

\addplot[only marks, color=blue, mark=*, mark options={fill=blue, draw=none}, mark size=2pt, opacity=0.7, xshift=-9.0pt, forget plot] coordinates {(1, 39.72) (1, 41.96) (1, 55.60) (1, 57.31) (1, 122.34) (2, 38.19) (2, 40.06) (2, 48.82) (2, 51.32) (2, 168.99) (3, 33.11) (3, 34.49) (3, 61.52) (3, 75.03) (3, 166.12) (4, 57.53) (4, 59.52) (4, 68.16) (4, 75.08) (4, 214.40)};
\addplot[only marks, color=red, mark=square*, mark options={fill=red, draw=none}, mark size=2pt, opacity=0.7, xshift=-3.0pt, forget plot] coordinates {(1, 43.88) (1, 45.45) (1, 46.13) (1, 46.56) (1, 150.17) (2, 28.28) (2, 29.16) (2, 58.31) (2, 116.79) (2, 168.41) (3, 65.65) (3, 67.22) (3, 156.94) (3, 219.03) (3, 258.54) (4, 98.84) (4, 106.67) (4, 120.36) (4, 155.97) (4, 258.72)};
\addplot[only marks, color=green!60!black, mark=triangle*, mark options={fill=green!60!black, draw=none}, mark size=3pt, opacity=0.7, xshift=3.0pt, forget plot] coordinates {(1, 45.12) (1, 48.00) (1, 50.82) (1, 126.79) (1, 131.96) (2, 4.80) (2, 5.07) (2, 45.60) (2, 91.70) (2, 230.24) (3, 79.79) (3, 83.70) (3, 196.73) (3, 229.03) (3, 446.27) (4, 25.47) (4, 26.53) (4, 79.59) (4, 371.91) (4, 703.06)};
\addplot[only marks, color=orange, mark=diamond*, mark options={fill=orange, draw=none}, mark size=3pt, opacity=0.7, xshift=9.0pt, forget plot] coordinates {(1, 43.19) (1, 45.12) (1, 46.25) (1, 46.48) (1, 158.13) (2, 20.08) (2, 20.99) (2, 155.11) (2, 353.33) (2, 407.91) (3, 38.86) (3, 40.34) (3, 261.30) (3, 266.84) (3, 275.55) (4, 214.22) (4, 220.85) (4, 310.04) (4, 399.23) (4, 4609.53)};

\end{axis}
\end{tikzpicture}
\caption{{
Runtime (left) and memory usage (right) for the MaxCut problem as a function of the TT rank, evaluated over multiple problem sizes (\(2^6\)–\(2^9\) vertices) using five random seeds; colours denote the problem size
}}
\label{fig:tt-rank-runtime-memory-maxcut}
\end{figure}

{
Figure~\ref{fig:tt-rank-runtime-memory-maxcut} presents the runtime (left) and memory usage (right) of our solver on MaxCut instances as a function of the TT rank of the input graph. Each point corresponds to a single random instance, with colours indicating problem sizes from $2^6$ to $2^9$. As the TT rank increases, both runtime and memory usage generally increase, with some seed-dependent outliers. This reflects the added computational and memory cost of handling more complex objective tensors in the MaxCut problem, induced by the higher TT rank of the input graph. The spread of runtimes and memory usage within each TT rank further reveals the impact of problem size, with larger instances incurring notably higher resource demands.
}

\begin{figure}[!htbp]
\centering
\begin{tikzpicture}
\begin{axis}[
    name=timeplot,
    ymode=log,
    width=0.47\linewidth,
    height=7cm,
    title={Runtime vs TT Rank},
    xlabel={TT Rank},
    ylabel={Time (s)},
    xtick={1, 2, 3, 4},
    enlarge y limits={value=0.1, upper},
    enlarge x limits=0.2,
    legend style={
        at={(0.02,0.98)},
        anchor=north west,
        font=\small,
        draw=none,
        fill=none,
    },
    label style={font=\small},
    tick label style={font=\small}
]

\addplot[only marks, color=blue, mark=*, mark options={fill=blue, draw=none}, mark size=2pt, opacity=0.7, xshift=-9.0pt, forget plot] coordinates {(1, 1.51) (1, 1.56) (1, 1.79) (1, 1.97) (1, 2.54) (2, 2.34) (2, 2.40) (2, 4.32) (2, 12.11) (2, 14.51) (3, 8.39) (3, 9.09) (3, 9.65) (3, 10.78) (3, 27.00) (4, 6.98) (4, 7.51) (4, 9.30) (4, 15.58) (4, 36.76)};
\addlegendimage{only marks, mark=*, mark options={fill=blue, draw=none}, color=blue}
\addplot[only marks, color=red, mark=square*, mark options={fill=red, draw=none}, mark size=2pt, opacity=0.7, xshift=-3.0pt, forget plot] coordinates {(1, 1.81) (1, 1.95) (1, 2.27) (1, 3.10) (1, 9.09) (2, 16.01) (2, 17.15) (2, 22.71) (2, 37.71) (2, 40.71) (3, 33.78) (3, 35.94) (3, 38.05) (3, 75.21) (3, 86.89) (4, 40.03) (4, 42.29) (4, 90.60) (4, 96.64) (4, 202.63)};
\addlegendimage{only marks, mark=square*, mark options={fill=red, draw=none}, color=red}
\addplot[only marks, color=green!60!black, mark=triangle*, mark options={fill=green!60!black, draw=none}, mark size=3pt, opacity=0.7, xshift=3.0pt, forget plot] coordinates {(1, 2.45) (1, 2.56) (1, 3.17) (1, 3.48) (1, 4.45) (2, 83.91) (2, 86.80) (2, 195.31) (2, 375.43) (2, 503.07) (3, 84.86) (3, 87.18) (3, 156.93) (3, 385.65) (3, 893.65) (4, 230.38) (4, 249.51) (4, 374.26) (4, 1259.98) (4, 2831.30)};
\addlegendimage{only marks, mark=triangle*, mark size=3pt, mark options={fill=green!60!black, draw=none}, color=green!60!black}
\addplot[only marks, color=orange, mark=diamond*, mark options={fill=orange, draw=none}, mark size=3pt, opacity=0.7, xshift=9.0pt, forget plot] coordinates {(1, 3.32) (1, 3.42) (1, 6.84) (1, 13.23) (1, 107.88) (2, 65.82) (2, 67.39) (2, 110.28) (2, 305.18) (2, 504.18) (3, 20.82) (3, 22.47) (3, 142.48) (3, 142.48) (3, 142.48) (4, 16.00) (4, 17.14) (4, 973.35) (4, 973.35) (4, 973.35)};
\addlegendimage{only marks, mark=diamond*, mark size=3pt, mark options={fill=orange, draw=none}, color=orange}

\legend{$2^6$, $2^7$, $2^8$, $2^9$}

\end{axis}
\begin{axis}[
    at={(timeplot.south east)},
    anchor=south west,
    xshift=1.55cm,
    ymode=log,
    width=0.47\linewidth,
    height=7cm,
    title={Memory vs TT Rank},
    xlabel={TT Rank},
    ylabel={Memory (MB)},
    xtick={1, 2, 3, 4},
    ytick={1e1, 1e2, 1e3, 1e4},
    ymin=1e1,
    enlarge y limits={value=0.1, upper},
    enlarge x limits=0.2,
    label style={font=\small},
    tick label style={font=\small},
]

\addplot[only marks, color=blue, mark=*, mark options={fill=blue, draw=none}, mark size=2pt, opacity=0.7, xshift=-9.0pt, forget plot] coordinates {(1, 36.32) (1, 38.37) (1, 51.25) (1, 52.86) (1, 125.63) (2, 47.09) (2, 49.40) (2, 87.13) (2, 97.91) (2, 161.83) (3, 87.69) (3, 91.34) (3, 102.55) (3, 108.15) (3, 200.81) (4, 83.73) (4, 86.62) (4, 90.81) (4, 94.16) (4, 207.96)};
\addplot[only marks, color=red, mark=square*, mark options={fill=red, draw=none}, mark size=2pt, opacity=0.7, xshift=-3.0pt, forget plot] coordinates {(1, 43.59) (1, 43.65) (1, 43.69) (1, 43.72) (1, 270.63) (2, 27.38) (2, 28.23) (2, 56.45) (2, 139.36) (2, 194.71) (3, 121.88) (3, 124.79) (3, 131.08) (3, 140.96) (3, 256.56) (4, 36.91) (4, 39.83) (4, 47.82) (4, 68.58) (4, 85.51)};
\addplot[only marks, color=green!60!black, mark=triangle*, mark options={fill=green!60!black, draw=none}, mark size=3pt, opacity=0.7, xshift=3.0pt, forget plot] coordinates {(1, 43.44) (1, 43.52) (1, 43.52) (1, 43.60) (1, 139.40) (2, 33.66) (2, 35.56) (2, 154.81) (2, 166.67) (2, 261.52) (3, 20.44) (3, 21.44) (3, 96.50) (3, 356.87) (3, 576.61) (4, 102.31) (4, 106.58) (4, 319.73) (4, 512.60) (4, 1460.36)};
\addplot[only marks, color=orange, mark=diamond*, mark options={fill=orange, draw=none}, mark size=3pt, opacity=0.7, xshift=9.0pt, forget plot] coordinates {(1, 9.01) (1, 9.48) (1, 56.88) (1, 298.44) (1, 356.95) (2, 214.45) (2, 224.16) (2, 330.70) (2, 371.67) (2, 610.21) (3, 41.91) (3, 43.50) (3, 518.28) (3, 518.28) (3, 518.28) (4, 691.75) (4, 713.14) (4, 1515.18) (4, 1515.18) (4, 1515.18)};

\end{axis}
\end{tikzpicture}
\caption{{
Runtime (left) and memory usage (right) for the Maximum Stable Set problem as a function of TT rank across problem sizes \(2^6\)–\(2^9\) using five random seeds; colours denote the problem size
}}

\label{fig:tt-rank-runtime-memory-maxstableset}
\end{figure}

{
Figure~\ref{fig:tt-rank-runtime-memory-maxstableset} shows the runtime (left) and memory usage (right) of TT IPM on Maximum Stable Set problems across varying input TT ranks and problem sizes. Similar to MaxCut, both runtime and memory increase with TT rank and dimension. However, several key differences emerge. Most notably, the runtime and memory variability across seeds is less regular than in the MaxCut experiments, especially for TT ranks 3 and 4. For instance, the plotted summaries for size $2^9$ show large spreads for ranks 1 and 4, and the rank-4 memory reaches about $1.5$ GB. In contrast to MaxCut, where performance changes more gradually, Max Stable Set exhibits more heterogeneous but still informative scaling.
}

\begin{figure}[!htbp]
\centering
\begin{tikzpicture}
\begin{axis}[
    name=timeplot,
    ymode=log,
    width=0.47\linewidth,
    height=7cm,
    title={Runtime vs TT Rank},
    xlabel={TT Rank},
    ylabel={Time (s)},
    xtick={1, 2, 3, 4},
    enlarge y limits={value=0.1, upper},
    enlarge x limits=0.2,
    legend style={
        at={(0.02,0.98)},
        anchor=north west,
        font=\small,
        draw=none,
        fill=none,
    },
    label style={font=\small},
    tick label style={font=\small}
]

\addplot[only marks, color=blue, mark=*, mark options={fill=blue, draw=none}, mark size=2pt, opacity=0.7, xshift=-9.0pt, forget plot] coordinates {(1, 1.01) (1, 1.05) (1, 1.32) (1, 1.54) (1, 3.13) (2, 1.49) (2, 1.53) (2, 2.18) (2, 2.99) (2, 6.37) (3, 4.12) (3, 4.46) (3, 6.69) (3, 15.36) (3, 22.64) (4, 7.57) (4, 8.14) (4, 10.46) (4, 42.42) (4, 53.09)};
\addlegendimage{only marks, mark=*, mark options={fill=blue, draw=none}, color=blue}
\addplot[only marks, color=red, mark=square*, mark options={fill=red, draw=none}, mark size=2pt, opacity=0.7, xshift=-3.0pt, forget plot] coordinates {(1, 3.32) (1, 3.58) (1, 4.46) (1, 6.75) (1, 9.91) (2, 4.61) (2, 4.94) (2, 7.15) (2, 10.70) (2, 11.41) (3, 157.17) (3, 167.20) (3, 278.88) (3, 400.42) (3, 407.57) (4, 87.52) (4, 92.45) (4, 102.65) (4, 103.93) (4, 241.40)};
\addlegendimage{only marks, mark=square*, mark options={fill=red, draw=none}, color=red}
\addplot[only marks, color=green!60!black, mark=triangle*, mark options={fill=green!60!black, draw=none}, mark size=3pt, opacity=0.7, xshift=3.0pt, forget plot] coordinates {(1, 1.44) (1, 1.50) (1, 4.50) (1, 16.63) (1, 53.88) (2, 9.87) (2, 10.21) (2, 22.98) (2, 66.54) (2, 286.72) (3, 68.04) (3, 69.91) (3, 125.83) (3, 604.51) (3, 6395.96) (4, 342.46) (4, 370.89) (4, 556.34) (4, 1645.19) (4, 3235.11)};
\addlegendimage{only marks, mark=triangle*, mark size=3pt, mark options={fill=green!60!black, draw=none}, color=green!60!black}
\addplot[only marks, color=orange, mark=diamond*, mark options={fill=orange, draw=none}, mark size=3pt, opacity=0.7, xshift=9.0pt, forget plot] coordinates {(1, 21.59) (1, 22.26) (1, 33.92) (1, 45.59) (1, 57.99) (2, 29.82) (2, 30.53) (2, 37.35) (2, 48.07) (2, 82.25) (3, 435.09) (3, 469.52) (3, 650.11) (3, 2278.56) (3, 6404.36) (4, 8579.54) (4, 9192.37) (4, 16472.36) (4, 17631.58) (4, 17863.42)};
\addlegendimage{only marks, mark=diamond*, mark size=3pt, mark options={fill=orange, draw=none}, color=orange}

\legend{$2^6$, $2^7$, $2^8$, $2^9$}

\end{axis}
\begin{axis}[
    at={(timeplot.south east)},
    anchor=south west,
    xshift=1.55cm,
    ymode=log,
    width=0.47\linewidth,
    height=7cm,
    title={Memory vs TT Rank},
    xlabel={TT Rank},
    ylabel={Memory (MB)},
    xtick={1, 2, 3, 4},
    ytick={1e1, 1e2, 1e3, 1e4},
    ymin=1e1,
    enlarge y limits={value=0.1, upper},
    enlarge x limits=0.2,
    label style={font=\small},
    tick label style={font=\small},
]

\addplot[only marks, color=blue, mark=*, mark options={fill=blue, draw=none}, mark size=2pt, opacity=0.7, xshift=-9.0pt, forget plot] coordinates {(1, 41.14) (1, 42.80) (1, 43.84) (1, 43.97) (1, 147.24) (2, 43.67) (2, 45.81) (2, 70.93) (2, 78.11) (2, 127.27) (3, 62.15) (3, 64.74) (3, 93.01) (3, 107.15) (3, 213.02) (4, 57.21) (4, 59.19) (4, 86.82) (4, 108.93) (4, 164.70)};
\addplot[only marks, color=red, mark=square*, mark options={fill=red, draw=none}, mark size=2pt, opacity=0.7, xshift=-3.0pt, forget plot] coordinates {(1, 39.36) (1, 40.85) (1, 74.19) (1, 94.45) (1, 126.30) (2, 63.99) (2, 65.97) (2, 129.63) (2, 168.66) (2, 207.69) (3, 88.79) (3, 90.91) (3, 150.08) (3, 253.70) (3, 311.55) (4, 158.42) (4, 170.96) (4, 198.21) (4, 269.07) (4, 309.55)};
\addplot[only marks, color=green!60!black, mark=triangle*, mark options={fill=green!60!black, draw=none}, mark size=3pt, opacity=0.7, xshift=3.0pt, forget plot] coordinates {(1, 42.20) (1, 44.89) (1, 47.53) (1, 195.14) (1, 625.89) (2, 5.21) (2, 5.50) (2, 49.53) (2, 286.62) (2, 359.91) (3, 9.26) (3, 9.71) (3, 43.70) (3, 283.17) (3, 3734.82) (4, 140.59) (4, 146.45) (4, 541.90) (4, 2742.66) (4, 3035.55)};
\addplot[only marks, color=orange, mark=diamond*, mark options={fill=orange, draw=none}, mark size=3pt, opacity=0.7, xshift=9.0pt, forget plot] coordinates {(1, 7.09) (1, 7.46) (1, 161.01) (1, 354.13) (1, 391.45) (2, 33.48) (2, 35.00) (2, 125.99) (2, 278.41) (2, 532.27) (3, 271.68) (3, 282.02) (3, 725.20) (3, 2482.62) (3, 6901.30) (4, 4805.29) (4, 4953.91) (4, 12151.98) (4, 12151.98) (4, 12151.98)};

\end{axis}
\end{tikzpicture}
\caption{{
Runtime (left) and memory usage (right) for the Correlation Clustering problem as a function of the TT rank, evaluated over multiple problem sizes (\(2^6\)–\(2^9\) vertices) using five random seeds; colours denote the problem size
}}

\label{fig:tt-rank-runtime-memory-corr_clust}
\end{figure}

{
Figure~\ref{fig:tt-rank-runtime-memory-corr_clust} reports the runtime (left) and memory usage (right) of TT IPM on Correlation Clustering problems as a function of the TT rank of the input graph. As in MaxCut and Maximum Stable Set, both runtime and memory generally increase with TT rank. However, the scaling pattern here is more instance-dependent. At small sizes ($2^6$ and $2^7$), runtimes and memory footprints remain comparatively modest even at higher TT ranks, whereas for $2^8$ and $2^9$, the spread becomes more pronounced. For example, the rank-4 summaries range from seconds on small instances to roughly $1.8\cdot 10^4$ seconds and over 12 GB on the largest plotted instances. This variability suggests that TT rank alone is not the sole driver of complexity in Correlation Clustering, as particular graph structures within the same rank class can produce substantially different workloads. Towards convergence, this is most likely due to the inequality slack tensor $\mathcal{T}$, which can become close to singular along the path. The resulting increase in TT rank naturally raises memory usage while still preserving a workable TT representation in many completed runs. Interestingly, even when the TT rank of $\mathcal{T}$ stays moderate, $\Delta \mathcal{T}$ can lead to an explosion in memory complexity within the projected KKT system.
}

\section{Discussion and Conclusions}\label{sec:discussion}

The proposed TT based interior-point method demonstrates a strong capability to solve large-scale SDP instances of combinatorial optimisation problems, which remain challenging or intractable for traditional sparse or direct methods. By leveraging the tensor-train format, the method efficiently handles high-dimensional variables and operators without relying on traditional sparsity assumptions, thus opening new possibilities for large-scale relaxations.

Despite lacking the extensive preprocessing routines, heuristic warm-starts, and advanced numerical stabilisation techniques employed by established solvers such as SDPA and SCS, the TT IPM demonstrates impressive performance across a broad range of problem instances. SDPA benefits from sophisticated presolving phases, including problem scaling techniques to improve numerical stability. Similarly, SCS applies equilibration, variable scaling, and adaptive step-size selection within its operator splitting framework to enhance convergence. Both solvers leverage mature, highly optimised C++ implementations with efficient memory management and multithreading capabilities. In contrast, TT IPM operates without such preprocessing enhancements or extensive tuning, yet it {remains effective on larger primal-dual SDP instances and on several higher TT rank cases}. Notably, the iteration count of the TT IPM is often significantly lower than that of SDPA, which is in line with the assumption that the TT format is not prohibitive to the IPM's superlinear convergence and can, in fact, accelerate convergence. This underscores the fundamental advantages of the TT format in handling high-dimensional structured semidefinite programs that are otherwise intractable for classical solvers.

From the computational viewpoint, the solver spends the majority of its runtime in the final iterations, striving to meet the specified absolute tolerance in the duality gap. However, variance across problem instances persists, and some instances are solved significantly faster, pointing towards instance-specific conditioning issues. The observations from the three problem classes indicate that increases in the objective TT rank may disproportionately inflate the TT rank of the primal variables, while increases in the TT rank of the constraint operators cause overall TT rank growth across variables. As the barrier parameter $\mu \to 0$, the method exhibits pronounced increases in variance and computational complexity, primarily due to the progressive ill-conditioning of the operators $\mathbfcal{L}_\mathcal{Z}$ and $\mathbfcal{L}_\mathcal{X}$, which approach rank-deficiency in the terminal phases of the algorithm. This challenge is further intensified by the intrinsic asymmetry of the KKT system within our formulation. Classical interior-point frameworks commonly alleviate these difficulties through symmetrisation techniques, e.g., by choosing the Nesterov-Todd direction \cite{Nesterov1998} or sparse adaptations \cite{Zhang2017}. These remedies are precluded in the tensor-train framework owing to the potential lack of sparsity and the inherent restriction within the TT framework that renders computations beyond extreme eigenvalue problems infeasible. The specific formulation of each problem significantly influences the final solution rank and convergence behaviour, as an incorrect placement of the idle multipliers may lead to higher constraint operator ranks.

At dimensions $2^{10}$ and above, an additional limiting factor arises: the DMRG-step size subroutine becomes a dominant cost, while it is negligible at smaller sizes. In particular, it is sufficient for a single eigenvector to exhibit near full-rank structure for the runtime and memory requirements to deteriorate substantially. This phenomenon highlights a fundamental vulnerability of the TT format, since the efficiency of low-rank methods hinges on the compressibility of all iterates, including eigenvectors needed for the step size calculation.

Despite these challenges, the ability to tackle large-scale problems beyond the reach of existing methods highlights the promise of the TT based approach, motivating further development in robust preconditioning strategies, reformulation strategies, pre-processing techniques, and the exploration of alternative search directions tailored to the TT based KKT system. We have demonstrated that a broad class of structured semidefinite programming problems, including Correlation Clustering, Maximum Cut, and Maximum Stable Set, admit low tensor-train rank solutions. This structural property was exploited to develop and analyse a TT based interior-point method capable of solving instances that are otherwise intractable for conventional sparse or dense approaches. {The analysis gives an inexact-Newton interpretation, local convergence under TT forcing conditions, and a safeguarded global convergence estimate.} Our empirical results establish three key findings. First, the central path for these problems can be traversed without incurring prohibitive memory costs as the TT ranks of the iterates remain moderate throughout the iterations, leading to low TT rank solutions. Second, the TT format does not appear to hinder the superlinear convergence behaviour of the interior-point method on the reported instances. Third, the method often outperforms established baselines in runtime for low TT rank instances.

There are several promising directions for further investigation. Most important among these is the design of preconditioning strategies specifically adapted to the inherently non-symmetric KKT systems arising in the TT format, which could markedly enhance numerical stability, particularly in the final stages of the algorithm where ill-conditioning becomes most pronounced. Equally important is the investigation of alternative problem formulations for constraint-rich settings, including those employing Schur complement-based representations, such as in Graph Matching. {Another natural direction is a TT analogue of SCS, based on first-order operator splitting applied to the homogeneous self-dual conic embedding, where the dominant linear-system operations would be carried out directly in tensor-train form.} Finally, adapting pre-processing techniques such as symmetry \cite{Gatermann2004, Vallentin2009} or facial reductions \cite{Drusvyatskiy2017, Kungurtsev2020} to the TT format could yield substantial improvements.

\section*{Acknowledgements}
This work was supported by the UK Engineering and Physical Sciences Research Council (EPSRC), grant number 2757464.

\section*{Declarations}

\subsection*{Conflict of Interest}
The authors have no relevant financial or non-financial interests to disclose.

\subsection*{Data Availability}
All data analysed during this study are publicly available.

\subsection*{Code Availability}
The full code was made available for review. We remark that a set of packages were used in this study, which were either open source or available for academic use. Specific references are included in this published article.

\bibliography{refs}

\global\let\myoriginalsection\section
\begin{appendices}
\pagebreak

\section{Test Instance Generation}\label{app:test_instance_generation}

To evaluate solver performance, we require graph instances with prescribed TT rank. However, generating such graphs through naive rejection sampling becomes computationally expensive, particularly for large sizes and low TT ranks. Generating a graph of size $2^d \times 2^d$ with TT rank $2$ via rejection sampling is highly inefficient, as the maximum possible rank is $2^d$, making successful samples very rare. Instead, we adopt a method that allows more precise control over rank and offers better scalability. We aim to sample from the set
\begin{align*}
    \mathfrak{S}_r := \{ \mathcal{S} \in \{0, 1\}^{2^d \times 2^d} : \mathcal{S}_{ii} = 0, S = S^\top,\ \textnormal{rank}_\textnormal{TT}(\mathcal{S}) \leq r \}.
\end{align*}

\paragraph{Sampling Binary TT Matrices}
We first generate a random orthonormal basis $Q \in \mathbb{R}^{r \times r}$ via QR-decomposition, i.e., $[Q, \_] = \text{QR}(\Omega)$ with $\Omega_{ij} \sim \mathcal{N}(0, 1)$. This basis is used to construct TT cores. For the first and last cores, we sample row vectors $q_i$ from $Q$ (with possible duplicates). For $d = 2$, an entry $(\alpha; \beta)$ in the matrix results in $1$ if the corresponding row vectors are parallel and $0$ if orthogonal.

For $d > 2$, intermediate cores are generated via random projections that permute the basis vectors. Specifically, for TT matrix $\mathcal{S}$ with cores $P_k$, we generate projection matrices of the form
\[
    P_k(\alpha_k, \beta_k) = I_r + \sum_{(i, j) \in \mathfrak{I}} q_i (q_j - q_i)^\top,
\]
where $\mathfrak{I}$ is a randomly chosen set of index pairs, effectively swapping basis vectors. The sparsity of the resulting matrix can be controlled by skewing the distribution from which the pairs in $\mathfrak{I}$ are sampled: a biased distribution yields fewer orthogonal combinations, resulting in denser matrices.

\paragraph{Sampling Adjacency TT Matrices without Self-Connections}
This follows a similar construction, with an additional constraint: projection matrices associated with diagonal indices must cancel the contribution of previous diagonal components and both off-diagonal projectors must be the same, i.e., the projectors must be $P_k(\alpha_k, \beta_k) = P_k(\beta_k, \alpha_k)$. Let $\mathfrak{D}$ be the set of diagonal index paths; then a diagonal projector is of the form
\[
P_k(\alpha_k, \alpha_k) = I_r + \sum_{(i, j) \in \mathfrak{I},\, i \notin \mathfrak{D}} q_i (q_j - q_i)^\top - \sum_{i \in \mathfrak{D}} q_i q_i^\top,
\]
ensuring that entries on the diagonal are zero. We then proceed with rejection sampling until $\rankTT(\mathcal{S}) = r$ with $\mathcal{S} \sim \mathfrak{S}_r$.

\section{Proofs for the TT IPM Analysis}\label{app:ttipm_proofs}

\subsection{Proof of \Cref{prop:projected-kkt-nonsingular}}\label{app:proof_projected_kkt}
{
\begin{proof}
Since $\mathcal X\succ0$ and $\mathcal Z\succ0$, the projected AHO blocks
$\widehat{\mathbfcal L}_{\mathcal X}$ and
$\widehat{\mathbfcal L}_{\mathcal Z}$ are positive definite. Hence the
elimination of $\Delta\mathcal Z$ is reversible, and
$\widehat{\mathbfcal M}\succ0$. Also
$\widehat{\mathbfcal D}\succeq0$ because
$\mathbfcal K_{\mathcal Y}\succeq0$, $\mathbfcal K_{\mathcal T}\succeq0$, and
$\mathcal D_{\mathrm{ineq}}(\mathcal X)\succeq0$. Moreover,
$\mathcal D_{\mathrm{ineq}}(\mathcal X)+\mathbfcal K_{\mathcal T}\succ0$, so
its projection is positive definite. After eliminating $\Delta\mathcal Z$, multiplying the second block row by
$\widehat{\mathbfcal L}_{\mathcal X}^{-1/2}$, applying the change of variables
induced by $\widehat{\mathbfcal L}_{\mathcal X}^{-1/2}$ to the remaining primal
block, and permuting the column blocks, we obtain the equivalent reduced matrix
\begin{align*}
    \widehat{\mathbfcal K}_{\mathrm{reduced}}
    =
    \begin{bmatrix}
        \widehat{\mathbfcal D}
        &
        -\widehat{\mathbfcal B}\widehat{\mathbfcal L}_{\mathcal X}^{-1/2}\\
        \widehat{\mathbfcal L}_{\mathcal X}^{1/2}\widehat{\mathbfcal A}^{\top}
        &
        \widehat{\mathbfcal M}
    \end{bmatrix}.
\end{align*}
Let $u=[y,t]^\top\in\mathbb R^{2p}$ with $y,t\in\mathbb R^p$, and let
$x\in\mathbb R^p$ satisfy
$\widehat{\mathbfcal K}_{\mathrm{reduced}}[u,x]^{\top}=0$. Taking the inner
product with $[u,x]^{\top}$ gives
\begin{align*}
    0
    =
    u^{\top}\widehat{\mathbfcal D}u
    +
    x^{\top}\widehat{\mathbfcal M}x
    +
    u^{\top}
    \left(
        \widehat{\mathbfcal A}\widehat{\mathbfcal L}_{\mathcal X}^{1/2}
        -
        \widehat{\mathbfcal B}\widehat{\mathbfcal L}_{\mathcal X}^{-1/2}
    \right)x .
\end{align*}
The preceding identity gives
\begin{align*}
    u^{\top}
    \left(
        \widehat{\mathbfcal A}\widehat{\mathbfcal L}_{\mathcal X}^{1/2}
        -
        \widehat{\mathbfcal B}\widehat{\mathbfcal L}_{\mathcal X}^{-1/2}
    \right)x
    =
    -u^{\top}\widehat{\mathbfcal D}u
    -x^{\top}\widehat{\mathbfcal M}x.
\end{align*}
Since both quadratic terms are nonnegative, the absolute value of the mixed
term equals their sum. Assumption (A2) therefore gives
\begin{align*}
    u^{\top}\widehat{\mathbfcal D}u
    +x^{\top}\widehat{\mathbfcal M}x
    \leq
    \gamma
    \bigl(u^{\top}\widehat{\mathbfcal D}u\bigr)^{1/2}
    \bigl(x^{\top}\widehat{\mathbfcal M}x\bigr)^{1/2}.
\end{align*}
By the arithmetic--geometric mean inequality,
\begin{align*}
    \bigl(u^{\top}\widehat{\mathbfcal D}u\bigr)^{1/2}
    \bigl(x^{\top}\widehat{\mathbfcal M}x\bigr)^{1/2}
    \leq
    \frac{1}{2}
    \left(
        u^{\top}\widehat{\mathbfcal D}u
        +x^{\top}\widehat{\mathbfcal M}x
    \right).
\end{align*}
Consequently,
\begin{align*}
    u^{\top}\widehat{\mathbfcal D}u
    +x^{\top}\widehat{\mathbfcal M}x
    \leq
    \frac{\gamma}{2}
    \left(
        u^{\top}\widehat{\mathbfcal D}u
        +x^{\top}\widehat{\mathbfcal M}x
    \right).
\end{align*}
Since $\gamma<2$, both quadratic terms must be zero. Since
$\widehat{\mathbfcal M}\succ0$, we have
$x=0$. Since $\widehat{\mathbfcal D}\succeq0$, we also have
$\widehat{\mathbfcal D}u=0$. The positive definite inequality
block of $\widehat{\mathbfcal D}$ gives $t=0$, while the equality
block gives $\mathcal{W}_{\neq k}^{\top} \mathbfcal K_{\mathcal Y} \mathcal{W}_{\neq k} y = 0$. The second block row gives
$\widehat{\mathbfcal L}_{\mathcal X}^{1/2}\widehat{\mathbfcal A}^{\top}u=0$,
hence $\widehat{\mathbfcal A}^{\top}u=0$. Since $t=0$,
$\left(\mathcal{W}_{\neq k}^{\top}\mathbfcal A_{\mathrm{eq}}\mathcal{W}_{\neq k}\right)^{\top}y=0$. By (A1), $y=0$, and hence $u=0$. Therefore
$\widehat{\mathbfcal K}_{\mathrm{reduced}}$ is nonsingular. Reversibility of the
preceding transformations proves nonsingularity of the projected AMEn KKT system.
\end{proof}
}

\subsection{Proof of \Cref{lem:tt-inexact-newton-direction}}\label{app:proof_tt_inexact_newton}
{
\begin{proof}
For any $\mathcal V$, the TT-rounding bound for $\widetilde{\mathcal K}=\mathrm{round}_\mathrm{TT}(\mathcal K,\frac{\delta_{\mathcal K}}{(1+\delta_{\mathcal K})\kappa_F(\mathcal K)})$ gives
\begin{align*}
    \|(\mathcal K-\widetilde{\mathcal K})\mathcal V\|
    \leq
    \|\mathcal K-\widetilde{\mathcal K} \|_F\|\mathcal V\|
    \leq
    \frac{\delta_{\mathcal K}}{(1+\delta_{\mathcal K})\kappa_F(\mathcal K)}
    \|\mathcal K\|_F\|\mathcal V\|.
\end{align*}
Since $\mathcal K$ is nonsingular, $\|\mathcal V\| \leq \|\mathcal K^{-1}\| \|\mathcal K \mathcal V\|$. Hence
\begin{align*}
    \|(\mathcal K-\widetilde{\mathcal K}) \mathcal V\|
    \leq
    \frac{\delta_{\mathcal K}}{1+\delta_{\mathcal K}}
    \|\mathcal K \mathcal V\|.
\end{align*}
Thus,
\begin{align*}
    \|\widetilde{\mathcal K} \mathcal V\|
    \geq
    \|\mathcal K \mathcal V\|-\|(\mathcal K-\widetilde{\mathcal K})\mathcal V\|
    \geq
    \left(1-\frac{\delta_{\mathcal K}}{1+\delta_{\mathcal K}}\right)\|\mathcal K \mathcal V\|.
\end{align*}
It follows that $\|(\mathcal K-\widetilde{\mathcal K})\mathcal V\| \leq \delta_{\mathcal K}\|\widetilde{\mathcal K} \mathcal V\|$. Choose $\mathcal V=\Delta \mathcal U^{\mathrm{tt}}$; then the rounded-system forcing
condition gives
\begin{align*}
    \|(\mathcal K-\widetilde{\mathcal K}) \Delta \mathcal U^{\mathrm{tt}}\| &\leq \delta_{\mathcal K}\|\widetilde{\mathcal K}\Delta \mathcal U^{\mathrm{tt}}\|  \\
    &\leq \delta_{\mathcal K} \left( \|\widetilde{\mathcal K}\Delta \mathcal U^{\mathrm{tt}}+\widetilde{\mathcal R}\| + \|\widetilde{\mathcal R}\| \right) \\
    &\leq \delta_{\mathcal K}(1+\widetilde{\eta})\|\widetilde{\mathcal R}\|.
\end{align*}
Note that $\mathcal K\Delta \mathcal U^{\mathrm{tt}}+\mathcal R = (\widetilde{\mathcal K}\Delta \mathcal U^{\mathrm{tt}}+\widetilde{\mathcal R}) +(\mathcal R-\widetilde{\mathcal R}) +(\mathcal K-\widetilde{\mathcal K})\Delta \mathcal U^{\mathrm{tt}}$. Thus,
\begin{align*}
    \|\mathcal K\Delta \mathcal U^{\mathrm{tt}}+\mathcal R\|
    &\leq
    \|\widetilde{\mathcal K}\Delta \mathcal U^{\mathrm{tt}}+\widetilde{\mathcal R}\|
    +\|\mathcal R-\widetilde{\mathcal R}\|
    +\|(\mathcal K-\widetilde{\mathcal K})\Delta \mathcal U^{\mathrm{tt}}\|  \\
    &\leq
    \widetilde{\eta}\|\widetilde{\mathcal R}\|
    +\delta_R\|\mathcal R\|
    +\delta_{\mathcal K}(1+\widetilde{\eta})\|\widetilde{\mathcal R}\|.
\end{align*}
Since $\|\widetilde{\mathcal R}\| \leq \|\mathcal R\|+\|\mathcal R-\widetilde{\mathcal R}\| \leq (1+\delta_R)\|\mathcal R\|$, we have
\begin{align*}
    \|\mathcal K\Delta \mathcal U^{\mathrm{tt}}+\mathcal R\| \leq \left[ \delta_R +(1+\delta_R)\widetilde{\eta} +(1+\delta_R)\delta_{\mathcal K}(1+\widetilde{\eta}) \right]\|\mathcal R\|.
\end{align*}
Choosing
\begin{align*}
    \delta_{\mathcal K}
    \leq
    \frac{\eta-\delta_R-(1+\delta_R)\widetilde{\eta}}
    {(1+\delta_R)(1+\widetilde{\eta})}
\end{align*}
therefore gives
$\|\mathcal K\Delta \mathcal U^{\mathrm{tt}}+\mathcal R\|
\leq \eta\|\mathcal R\|$.
\end{proof}
}

\subsection{Proof of \Cref{th:tt-local-convergence}}\label{app:proof_tt_local_convergence}
{
\begin{proof}
By \Cref{lem:tt-inexact-newton-direction},
\begin{align*}
    \mathcal K_k\Delta\mathcal U_k^{\mathrm{tt}}
    =
    -\mathcal R_k+\mathcal E_k,
    \qquad
    \|\mathcal E_k\|
    \leq
    \eta_k\|\mathcal R_k\|,
\end{align*}
with
\begin{align*}
    \eta_k
    &:=
    \delta_{R,k}
    +(1+\delta_{R,k})\widetilde{\eta}_k
    +(1+\delta_{R,k})\delta_{\mathcal K,k}(1+\widetilde{\eta}_k)
    \leq \bar\eta<1 .
\end{align*}
By the local AHO regularity theory \cite[Theorem~3.1]{Alizadeth1998} and
(L1), $\mathcal K_k$ is uniformly nonsingular and $\mathcal R_{\mu_k}$ has
locally Lipschitz continuous derivative. Hence
$\|\Delta\mathcal U_k^{\mathrm{tt}}\|\leq c\|\mathcal R_k\|$ locally. Taylor's theorem at
$\mathcal U_k$ gives, for
$\widehat{\mathcal U}_{k+1}:=\mathcal U_k+\alpha_k\Delta\mathcal U_k^{\mathrm{tt}}$,
\begin{align*}
    \|\mathcal R_{\mu_k}(\widehat{\mathcal U}_{k+1})\|
    &\leq
    \left(1-\alpha_k(1-\eta_k)\right)\|\mathcal R_k\|
    +O(\|\Delta\mathcal U_k^{\mathrm{tt}}\|^2)\\
    &\leq
    \left(1-\alpha_k(1-\bar\eta)\right)\|\mathcal R_k\|
    +O(\|\mathcal R_k\|^2).
\end{align*}
After shrinking the neighbourhood if necessary, the quadratic term is absorbed into
the linear decrease. Since $\alpha_k\geq\underline\alpha>0$, there exists
$q_0\in(0,1)$ such that
\begin{align*}
    \|\mathcal R_{\mu_k}(\widehat{\mathcal U}_{k+1})\|
    \leq
    q_0\|\mathcal R_{\mu_k}(\mathcal U_k)\|.
\end{align*}
The actual iterate is obtained by post-step TT rounding and spectral correction.
By (L3), this contributes only the residual perturbation $\rho_k$:
\begin{align*}
    \|\mathcal R_{\mu_k}(\mathcal U_{k+1})\|
    \leq
    q_0\|\mathcal R_{\mu_k}(\mathcal U_k)\|+\rho_k.
\end{align*}
Finally, locally
$\|\mathcal R_\mu(\mathcal U)-\mathcal R_0(\mathcal U)\|\leq c_\mu\mu$, since the
residuals differ only in the barrier complementarity terms. Therefore
\begin{align*}
    \|\mathcal R_0(\mathcal U_{k+1})\|
    &\leq
    \|\mathcal R_{\mu_k}(\mathcal U_{k+1})\|+c_\mu\mu_k\\
    &\leq
    q_0\|\mathcal R_{\mu_k}(\mathcal U_k)\|+\rho_k+c_\mu\mu_k\\
    &\leq
    q_0\|\mathcal R_0(\mathcal U_k)\|
    +\rho_k+(1+q_0)c_\mu\mu_k .
\end{align*}
This proves the claimed estimate with a sufficiently large constant
$C>0$. Taking limsups in the resulting scalar recursion gives the stated bound.
\end{proof}
}

\subsection{Proof of \Cref{th:tt-global-convergence}}\label{app:proof_tt_global_convergence}
{
\begin{proof}
By (G1), all iterates remain in the strict interior, so the barrier residuals and
Newton systems are well defined. By (G2) and
\Cref{lem:tt-inexact-newton-direction}, each computed TT direction satisfies the
inexact Newton forcing condition with $\eta_k\leq\bar\eta<1$.
Assumptions (G3)--(G4) collect the remaining TT perturbations into the windowed scalar recursion for $\Theta_k$.

Let $\bar\rho:=\limsup_{k\to\infty}\rho_k$. Since $\Omega_k$ is the maximum of finitely many recent perturbation bounds, $\limsup_{k\to\infty}\Omega_k=\bar\rho$. If $\bar\rho=\infty$, the result is immediate. Otherwise, apply (G4) along each residue class modulo $M+1$ and take the limsup to obtain, with $\bar\Theta:=\limsup_{k\to\infty}\Theta_k$,
\begin{align*}
    \bar\Theta
    \leq
    (1-\tau)\bar\Theta+c_\Omega\bar\rho .
\end{align*}
Hence
\begin{align*}
    \limsup_{k\to\infty}\Theta_k
    \leq
    \frac{c_\Omega}{\tau}\bar\rho .
\end{align*}
Since $\Psi_k\leq\Theta_k$, the same bound holds for $\Psi_k$.
Since $\mathcal U_k$ is bounded, the barrier terms give, for some $C_0>0$,
\begin{align*}
    \|\mathcal R_0(\mathcal U_k)\|+\mu_k
    &\le
    C_0\bigl(\|\mathcal R_{\mu_k}(\mathcal U_k)\|+\mu_k\bigr)\\
    &=
    C_0\Psi_k .
\end{align*}
Combining the estimates gives the stated bound, after increasing the constant $C$ if necessary. If $\rho_k\to0$, then $\Psi_k\to0$, so $\mu_k\to0$ and $\|\mathcal R_0(\mathcal U_k)\|\to0$; continuity of the KKT residual then makes every accumulation point a KKT point. The final assertion is exactly \Cref{th:tt-local-convergence} applied to the eventual local tail.
\end{proof}
}

\subsection{Proof of \Cref{cor:psdness_rounding}}\label{app:psdness_rounding}
\begin{proof}
By Weyl's inequality,
\begin{align*}
    \lambda_{\min}(\mathcal M)
    &\leq
    \lambda_{\min}(\widetilde{\mathcal M})
    +
    \|\mathcal M-\widetilde{\mathcal M}\|_2
    \leq
    \lambda_{\min}(\widetilde{\mathcal M})
    +
    \zeta,\\
    \lambda_{\min}(\widetilde{\mathcal M})
    &\leq
    \lambda_{\min}(\mathcal M)
    +
    \|\mathcal M-\widetilde{\mathcal M}\|_2
    \leq
    \lambda_{\min}(\mathcal M)
    +
    \zeta.
\end{align*}
Combining these two inequalities gives
$|\lambda_{\min}(\mathcal M)-\lambda_{\min}(\widetilde{\mathcal M})|
\leq\zeta$ as desired.
\end{proof}

\section{Detailed Numerical Results}\label{app:numerical_tables}

{
The main trends in runtime, memory, and TT-rank behaviour are visualised in Figures~\ref{fig:rank_evolution_maxcut}--\ref{fig:rank_evolution_corr_clust} and Figures~\ref{fig:tt-rank-runtime-memory-maxcut}--\ref{fig:tt-rank-runtime-memory-corr_clust}. The following tables provide the complete numerical values underlying the rank-stratified comparison, including accuracy and missing-run information that cannot be represented compactly in the plots.
}

\begingroup
\setlength{\LTleft}{0pt}
\setlength{\LTright}{0pt}
\fontsize{6}{7}\selectfont
\LTcapwidth=\textwidth
\setlength{\tabcolsep}{0pt}
\renewcommand{\arraystretch}{1.08}
\begin{longtable}{@{\extracolsep{\fill}}cc l c c c c c c c c@{}}
\caption{{Performance and accuracy comparison for MaxCut instances with varying input TT rank $r_{\mathrm{in}}$. Runtime is reported in seconds and memory in MB. Runtime and memory entries show median [IQR] over seed runs; the max columns report the worst observed value. Iteration counts and accuracy entries report arithmetic means over completed seed runs. Accuracy columns report final primal feasibility error, dual feasibility error, and duality gap. For each size, TT IPM is shown only up to input rank \(\min(r_{\mathrm{stop}},5)\), where \(r_{\mathrm{stop}}\) is the first rank at which it becomes OOT. If TT IPM is reported for an instance but SDPA or SCS is absent, the missing solver was OOT/OOM. MANOPT is included as a separate MaxCut low-rank baseline; it was run on the listed instances, but only runtime, memory, iterations, and primal feasibility are reported because it does not return comparable primal-dual certificates. Errors which did not achieve the required target tolerance are printed in bold. }}
\label{tab:comparison_maxcut_ranks}\\
\toprule
\multicolumn{2}{c}{Inst.}
& Method
& \multicolumn{4}{c}{Perf.}
& \multicolumn{4}{c}{Acc.} \\
\cmidrule(lr){1-2}
\cmidrule(lr){4-7}
\cmidrule(lr){8-11}
Size & $r_{\mathrm{in}}$
&
& \shortstack{Time\\med.}
& \shortstack{Time\\max}
& \shortstack{Mem.\\med.}
& \shortstack{Mem.\\max}
& Iter.
& \shortstack{Primal\\err.}
& \shortstack{Dual\\err.}
& Gap \\
\midrule
\endfirsthead
\toprule
\multicolumn{11}{c}{\tablename\ \thetable{} -- cont.}\\
\midrule
Size & $r_{\mathrm{in}}$ & Method & \shortstack{Time\\med.} & \shortstack{Time\\max} & \shortstack{Mem.\\med.} & \shortstack{Mem.\\max} & Iter. & \shortstack{Primal\\err.} & \shortstack{Dual\\err.} & Gap \\
\midrule
\endhead
\midrule
\multicolumn{11}{r}{Continued}\\
\endfoot
\bottomrule
\endlastfoot
${2^{6}}$ & 1 & TT IPM & $1.355\,[0.934]$ & $2.411$ & $55.6\,[15.35]$ & $122.34$ & $8$ & $2.11\mathrm{e}{-9}$ & $1.81\mathrm{e}{-9}$ & $3.64\mathrm{e}{-4}$ \\
 &  & SCS & $0.065\,[0.009]$ & $0.231$ & $71.40\,[0.656]$ & $203.13$ & $155$ & $1.04\mathrm{e}{-15}$ & $3.70\mathrm{e}{-7}$ & $3.17\mathrm{e}{-4}$ \\
 &  & SDPA & $3.14\,[12.48]$ & $16.93$ & \cellcolor[gray]{0.9}$49.68\,[104.47]$ & $155.14$ & $13.4$ & $1.02\mathrm{e}{-9}$ & $6.73\mathrm{e}{-7}$ & $9.48\mathrm{e}{-7}$ \\
 &  & MANOPT & \cellcolor[gray]{0.9}$0.006\,[0.007]$ & $0.011$ & $71.38\,[0.098]$ & $200.20$ & $6.6$ & $3.06\mathrm{e}{-30}$ & $--$ & $--$ \\
${2^{6}}$ & 2 & TT IPM & $1.92\,[0.398]$ & $3.5$ & \cellcolor[gray]{0.9}$48.82\,[11.26]$ & $168.99$ & $7.6$ & $2.59\mathrm{e}{-9}$ & $1.73\mathrm{e}{-9}$ & $3.81\mathrm{e}{-4}$ \\
 &  & SCS & $0.218\,[0.052]$ & $0.259$ & $71.42\,[0.43]$ & $202.93$ & $315$ & $1.84\mathrm{e}{-16}$ & $4.67\mathrm{e}{-9}$ & $1.56\mathrm{e}{-4}$ \\
 &  & SDPA & $2.89\,[2.87]$ & $17.21$ & $49.51\,[36.01]$ & $154.09$ & $12.6$ & $3.39\mathrm{e}{-11}$ & $1.11\mathrm{e}{-8}$ & $3.31\mathrm{e}{-6}$ \\
 &  & MANOPT & \cellcolor[gray]{0.9}$0.008\,[0.002]$ & $0.015$ & $71.10\,[0.137]$ & $199.88$ & $7.4$ & $2.42\mathrm{e}{-30}$ & $--$ & $--$ \\
${2^{6}}$ & 3 & TT IPM & $4.23\,[0.971]$ & $5.61$ & $61.52\,[40.54]$ & $166.12$ & $9.2$ & $6.67\mathrm{e}{-9}$ & $1.07\mathrm{e}{-9}$ & $2.78\mathrm{e}{-4}$ \\
 &  & SCS & $0.24\,[0.076]$ & $0.294$ & \cellcolor[gray]{0.9}$44.51\,[0.406]$ & $125.40$ & $285$ & $2.61\mathrm{e}{-16}$ & $9.68\mathrm{e}{-11}$ & $3.58\mathrm{e}{-4}$ \\
 &  & SDPA & $3.06\,[0.329]$ & $19.58$ & $49.82\,[3.77]$ & $154.37$ & $13$ & $4.60\mathrm{e}{-11}$ & $1.45\mathrm{e}{-12}$ & $3.97\mathrm{e}{-6}$ \\
 &  & MANOPT & \cellcolor[gray]{0.9}$0.015\,[0.001]$ & $0.062$ & $71.34\,[0.203]$ & $199.50$ & $11.4$ & $2.35\mathrm{e}{-30}$ & $--$ & $--$ \\
${2^{6}}$ & 4 & TT IPM & $9.19\,[3.59]$ & $21.98$ & $68.16\,[15.56]$ & $214.4$ & $10.6$ & $9.43\mathrm{e}{-9}$ & $1.88\mathrm{e}{-11}$ & $3.39\mathrm{e}{-4}$ \\
 &  & SCS & $0.293\,[0.135]$ & $0.416$ & \cellcolor[gray]{0.9}$44.77\,[0.645]$ & $125.47$ & $335$ & $2.24\mathrm{e}{-17}$ & $1.57\mathrm{e}{-8}$ & $6.08\mathrm{e}{-5}$ \\
 &  & SDPA & $3.04\,[0.261]$ & $3.67$ & $49.63\,[0.348]$ & $50.52$ & $12.8$ & $1.96\mathrm{e}{-10}$ & $3.42\mathrm{e}{-8}$ & $3.55\mathrm{e}{-6}$ \\
 &  & MANOPT & \cellcolor[gray]{0.9}$0.02\,[0.005]$ & $0.056$ & $71.12\,[0.082]$ & $199.12$ & $9.6$ & $2.71\mathrm{e}{-30}$ & $--$ & $--$ \\
${2^{6}}$ & 5 & TT IPM & $10.96\,[1.05]$ & $12.37$ & $93.29\,[14.03]$ & $199.74$ & $8.4$ & $1.65\mathrm{e}{-8}$ & $3.25\mathrm{e}{-11}$ & $4.84\mathrm{e}{-4}$ \\
 &  & SCS & $0.267\,[0.072]$ & $0.578$ & $71.22\,[0.676]$ & $203.23$ & $425$ & $5.03\mathrm{e}{-16}$ & $1.11\mathrm{e}{-7}$ & $9.01\mathrm{e}{-5}$ \\
 &  & SDPA & $3.26\,[0.306]$ & $16$ & \cellcolor[gray]{0.9}$49.69\,[0.922]$ & $154.51$ & $13.2$ & $6.42\mathrm{e}{-12}$ & $5.06\mathrm{e}{-8}$ & $5.24\mathrm{e}{-6}$ \\
 &  & MANOPT & \cellcolor[gray]{0.9}$0.018\,[0.009]$ & $0.021$ & $71.14\,[0.086]$ & $199.26$ & $8.8$ & $2.68\mathrm{e}{-30}$ & $--$ & $--$ \\

\midrule
${2^{7}}$ & 1 & TT IPM & $1.359\,[0.448]$ & $3.97$ & \cellcolor[gray]{0.9}$46.13\,[1.11]$ & $150.17$ & $8.6$ & $8.54\mathrm{e}{-9}$ & $1.09\mathrm{e}{-9}$ & $3.37\mathrm{e}{-4}$ \\
 &  & SCS & $0.785\,[0.221]$ & $1.35$ & $76.24\,[1.22]$ & $212.04$ & $230$ & $4.24\mathrm{e}{-15}$ & $1.04\mathrm{e}{-8}$ & $1.45\mathrm{e}{-3}$ \\
 &  & SDPA & $143.81\,[466.75]$ & $602.26$ & $551.47\,[3107.66]$ & $3395.64$ & $13.6$ & $1.21\mathrm{e}{-9}$ & $4.64\mathrm{e}{-8}$ & $2.59\mathrm{e}{-6}$ \\
 &  & MANOPT & \cellcolor[gray]{0.9}$0.006\,[0.007]$ & $0.015$ & $71.64\,[0.957]$ & $200.97$ & $6.2$ & $5.66\mathrm{e}{-30}$ & $--$ & $--$ \\
${2^{7}}$ & 2 & TT IPM & $2.56\,[1.23]$ & $4.78$ & \cellcolor[gray]{0.9}$58.31\,[87.63]$ & $168.41$ & $8.2$ & $8.40\mathrm{e}{-9}$ & $8.40\mathrm{e}{-10}$ & $4.14\mathrm{e}{-4}$ \\
 &  & SCS & $0.598\,[0.638]$ & $1.14$ & $76.19\,[0.059]$ & $210.75$ & $290$ & $9.71\mathrm{e}{-16}$ & $7.97\mathrm{e}{-9}$ & $4.13\mathrm{e}{-4}$ \\
 &  & SDPA & $20.29\,[156.8]$ & $584$ & $272.58\,[396.66]$ & $3395.3$ & $13.8$ & $9.32\mathrm{e}{-10}$ & $3.24\mathrm{e}{-8}$ & $3.53\mathrm{e}{-6}$ \\
 &  & MANOPT & \cellcolor[gray]{0.9}$0.014\,[0.004]$ & $0.032$ & $71.78\,[0.539]$ & $200.93$ & $8.2$ & $4.89\mathrm{e}{-30}$ & $--$ & $--$ \\
${2^{7}}$ & 3 & TT IPM & $13.82\,[28.73]$ & $37.85$ & $156.94\,[151.81]$ & $258.54$ & $10$ & $1.42\mathrm{e}{-8}$ & $1.09\mathrm{e}{-10}$ & $3.82\mathrm{e}{-4}$ \\
 &  & SCS & $1.15\,[0.551]$ & $1.65$ & \cellcolor[gray]{0.9}$50.20\,[0.148]$ & $133.75$ & $435$ & $2.83\mathrm{e}{-16}$ & $1.05\mathrm{e}{-7}$ & $1.80\mathrm{e}{-4}$ \\
 &  & SDPA & $170.05\,[165.56]$ & $583.26$ & $551.56\,[279.25]$ & $3395.82$ & $14$ & $1.74\mathrm{e}{-9}$ & $5.18\mathrm{e}{-8}$ & $3.24\mathrm{e}{-6}$ \\
 &  & MANOPT & \cellcolor[gray]{0.9}$0.02\,[0.008]$ & $0.04$ & $71.53\,[0.234]$ & $201.51$ & $9.2$ & $6.10\mathrm{e}{-30}$ & $--$ & $--$ \\
${2^{7}}$ & 4 & TT IPM & $7.75\,[10.08]$ & $37.93$ & $120.36\,[49.3]$ & $258.72$ & $9.6$ & $2.37\mathrm{e}{-8}$ & $1.09\mathrm{e}{-9}$ & $4.03\mathrm{e}{-4}$ \\
 &  & SCS & $1.08\,[0.032]$ & $1.45$ & \cellcolor[gray]{0.9}$49.97\,[3.66]$ & $133.46$ & $330$ & $7.42\mathrm{e}{-16}$ & $2.50\mathrm{e}{-9}$ & $9.43\mathrm{e}{-5}$ \\
 &  & SDPA & $170.44\,[14.93]$ & $183.57$ & $551.45\,[0.086]$ & $551.53$ & $13$ & $6.38\mathrm{e}{-11}$ & $2.29\mathrm{e}{-7}$ & $3.50\mathrm{e}{-6}$ \\
 &  & MANOPT & \cellcolor[gray]{0.9}$0.024\,[0.004]$ & $0.034$ & $71.81\,[0.301]$ & $201.21$ & $9.6$ & $5.60\mathrm{e}{-30}$ & $--$ & $--$ \\
${2^{7}}$ & 5 & TT IPM & $41.95\,[65.06]$ & $184.5$ & $149.91\,[76.22]$ & $242.57$ & $10.2$ & $1.59\mathrm{e}{-8}$ & $8.74\mathrm{e}{-10}$ & $5.05\mathrm{e}{-4}$ \\
 &  & SCS & $1.51\,[0.213]$ & $1.65$ & $76.40\,[0.328]$ & $210.62$ & $400$ & $3.21\mathrm{e}{-17}$ & $6.90\mathrm{e}{-8}$ & $1.93\mathrm{e}{-5}$ \\
 &  & SDPA & $169.38\,[0.601]$ & $195.2$ & $551.58\,[0.148]$ & $551.77$ & $13.4$ & $1.53\mathrm{e}{-11}$ & $1.43\mathrm{e}{-7}$ & $5.18\mathrm{e}{-6}$ \\
 &  & MANOPT & \cellcolor[gray]{0.9}$0.04\,[0.017]$ & $0.044$ & \cellcolor[gray]{0.9}$71.47\,[0.133]$ & $201.25$ & $11.8$ & $5.42\mathrm{e}{-30}$ & $--$ & $--$ \\

\midrule
${2^{8}}$ & 1 & TT IPM & $1.856\,[0.324]$ & $7.53$ & \cellcolor[gray]{0.9}$50.82\,[78.79]$ & $131.96$ & $9.2$ & $3.49\mathrm{e}{-9}$ & $2.84\mathrm{e}{-11}$ & $2.45\mathrm{e}{-4}$ \\
 &  & SCS & $2.21\,[1.69]$ & $5.47$ & $92.25\,[1.87]$ & $241.06$ & $210$ & $1.11\mathrm{e}{-14}$ & $6.40\mathrm{e}{-9}$ & $\mathbf{1.36\mathrm{e}{-2}}$ \\
 &  & SDPA & $1327.3\,[1965.11]$ & $8280.88$ & $6879.77\,[1078.77]$ & $8284.94$ & $14$ & $1.12\mathrm{e}{-9}$ & $5.57\mathrm{e}{-8}$ & $5.32\mathrm{e}{-7}$ \\
 &  & MANOPT & \cellcolor[gray]{0.9}$0.019\,[0.005]$ & $0.045$ & $73.09\,[0.492]$ & $206.09$ & $5.8$ & $1.57\mathrm{e}{-29}$ & $--$ & $--$ \\
${2^{8}}$ & 2 & TT IPM & $2.01\,[0.504]$ & $10.37$ & \cellcolor[gray]{0.9}$45.6\,[86.63]$ & $230.24$ & $8.6$ & $4.24\mathrm{e}{-9}$ & $9.94\mathrm{e}{-11}$ & $1.98\mathrm{e}{-4}$ \\
 &  & SCS & $5.35\,[1.93]$ & $7.75$ & $91.99\,[0.594]$ & $241.12$ & $310$ & $3.25\mathrm{e}{-16}$ & $6.84\mathrm{e}{-7}$ & $9.18\mathrm{e}{-4}$ \\
 &  & SDPA & $1607.15\,[2707.08]$ & $8032.45$ & $3398.61\,[532.91]$ & $8284.11$ & $14.2$ & $8.62\mathrm{e}{-10}$ & $3.89\mathrm{e}{-8}$ & $7.25\mathrm{e}{-7}$ \\
 &  & MANOPT & \cellcolor[gray]{0.9}$0.027\,[0.007]$ & $0.044$ & $72.26\,[0.66]$ & $206.54$ & $8.8$ & $1.41\mathrm{e}{-29}$ & $--$ & $--$ \\
${2^{8}}$ & 3 & TT IPM & $36.71\,[31.99]$ & $208.1$ & $196.73\,[145.33]$ & $446.27$ & $10$ & $2.53\mathrm{e}{-8}$ & $3.12\mathrm{e}{-9}$ & $5.46\mathrm{e}{-4}$ \\
 &  & SCS & $8.79\,[3.61]$ & $9.93$ & \cellcolor[gray]{0.9}$66.26\,[0.812]$ & $165.40$ & $460$ & $7.24\mathrm{e}{-17}$ & $2.49\mathrm{e}{-8}$ & $5.67\mathrm{e}{-4}$ \\
 &  & SDPA & $1568.87\,[2322.76]$ & $8015.89$ & $6880.89\,[1078.95]$ & $8285.38$ & $14.4$ & $1.61\mathrm{e}{-9}$ & $6.22\mathrm{e}{-8}$ & $6.65\mathrm{e}{-7}$ \\
 &  & MANOPT & \cellcolor[gray]{0.9}$0.027\,[0.008]$ & $0.044$ & $72.52\,[0.395]$ & $206.16$ & $8.6$ & $1.48\mathrm{e}{-29}$ & $--$ & $--$ \\
${2^{8}}$ & 4 & TT IPM & $17.34\,[19.36]$ & $662.74$ & $79.59\,[345.38]$ & $703.06$ & $10.6$ & $2.60\mathrm{e}{-8}$ & $6.58\mathrm{e}{-9}$ & $5.25\mathrm{e}{-4}$ \\
 &  & SCS & $7.28\,[2.87]$ & $10.07$ & \cellcolor[gray]{0.9}$65.92\,[0.504]$ & $165.73$ & $435$ & $8.19\mathrm{e}{-17}$ & $4.41\mathrm{e}{-8}$ & $1.36\mathrm{e}{-4}$ \\
 &  & SDPA & $1472.85\,[2328.66]$ & $2525.67$ & $6879.52\,[1073.53]$ & $13421.16$ & $13.4$ & $5.90\mathrm{e}{-11}$ & $2.75\mathrm{e}{-7}$ & $7.19\mathrm{e}{-7}$ \\
 &  & MANOPT & \cellcolor[gray]{0.9}$0.027\,[0.002]$ & $0.032$ & $72.18\,[0.102]$ & $206.11$ & $9$ & $1.25\mathrm{e}{-29}$ & $--$ & $--$ \\
${2^{8}}$ & 5 & TT IPM & $303.82\,[46.39]$ & $435.1$ & $261.63\,[159.74]$ & $618.64$ & $11.4$ & $2.68\mathrm{e}{-8}$ & $3.80\mathrm{e}{-9}$ & $3.85\mathrm{e}{-4}$ \\
 &  & SCS & $33.78\,[43.71]$ & $144.20$ & $92.23\,[0.758]$ & $243.41$ & $2565$ & $1.19\mathrm{e}{-16}$ & $1.73\mathrm{e}{-6}$ & $9.92\mathrm{e}{-5}$ \\
 &  & SDPA & $1569.56\,[2314.90]$ & $2683.00$ & $6881.12\,[1081.11]$ & $13422.74$ & $13.8$ & $1.42\mathrm{e}{-11}$ & $1.72\mathrm{e}{-7}$ & $1.06\mathrm{e}{-6}$ \\
 &  & MANOPT & \cellcolor[gray]{0.9}$0.04\,[0.005]$ & $0.068$ & \cellcolor[gray]{0.9}$72.64\,[0.637]$ & $206.89$ & $11$ & $1.30\mathrm{e}{-29}$ & $--$ & $--$ \\
\midrule
${2^{9}}$ & 1 & TT IPM & $3.67\,[0.45]$ & $4.51$ & \cellcolor[gray]{0.9}$46.25\,[1.36]$ & $158.13$ & $9.8$ & $1.45\mathrm{e}{-8}$ & $2.69\mathrm{e}{-9}$ & $6.21\mathrm{e}{-4}$ \\
 &  & SCS & $36.63\,[44.41]$ & $67.21$ & $161.81\,[26.29]$ & $359.48$ & $300$ & $2.04\mathrm{e}{-13}$ & $3.24\mathrm{e}{-5}$ & $\mathbf{1.73\mathrm{e}{-2}}$ \\
 &  & MANOPT & \cellcolor[gray]{0.9}$0.041\,[0.019]$ & $0.081$ & $81.01\,[3.69]$ & $217.98$ & $8.6$ & $3.05\mathrm{e}{-29}$ & $--$ & $--$ \\
${2^{9}}$ & 2 & TT IPM & $10.63\,[34.2]$ & $51.53$ & $155.11\,[332.34]$ & $407.91$ & $10.6$ & $5.86\mathrm{e}{-8}$ & $1.01\mathrm{e}{-10}$ & $5.49\mathrm{e}{-4}$ \\
 &  & SCS & $43.96\,[44.27]$ & $518.82$ & $162.91\,[34.16]$ & $366.53$ & $1180$ & $3.56\mathrm{e}{-16}$ & $1.57\mathrm{e}{-6}$ & $1.39\mathrm{e}{-3}$ \\
 &  & MANOPT & \cellcolor[gray]{0.9}$0.468\,[0.233]$ & $2.60$ & \cellcolor[gray]{0.9}$78.61\,[8.27]$ & $213.98$ & $9.2$ & $3.06\mathrm{e}{-29}$ & $--$ & $--$ \\
${2^{9}}$ & 3 & TT IPM & $12.97\,[38.15]$ & $88.83$ & $261.3\,[226.5]$ & $275.55$ & $10.4$ & $5.44\mathrm{e}{-8}$ & $7.75\mathrm{e}{-9}$ & $7.14\mathrm{e}{-4}$ \\
 &  & SCS & $79.29\,[8.88]$ & $103.02$ & $136.94\,[31.21]$ & $288.81$ & $570$ & $8.11\mathrm{e}{-17}$ & $3.78\mathrm{e}{-6}$ & $1.48\mathrm{e}{-3}$ \\
 &  & MANOPT & \cellcolor[gray]{0.9}$0.055\,[0.015]$ & $0.141$ & \cellcolor[gray]{0.9}$76.91\,[2.46]$ & $216.80$ & $9.8$ & $3.29\mathrm{e}{-29}$ & $--$ & $--$ \\
${2^{9}}$ & 4 & TT IPM & $47.91\,[65.87]$ & $3750.51$ & $310.04\,[178.38]$ & $4609.53$ & $12$ & $5.74\mathrm{e}{-8}$ & $8.75\mathrm{e}{-10}$ & $6.70\mathrm{e}{-4}$ \\
 &  & SCS & $572.30\,[983.48]$ & $1076.18$ & $136.27\,[34.78]$ & $287.87$ & $3720$ & $1.45\mathrm{e}{-17}$ & $1.88\mathrm{e}{-5}$ & $5.66\mathrm{e}{-4}$ \\
 &  & MANOPT & \cellcolor[gray]{0.9}$0.063\,[0.031]$ & $0.087$ & \cellcolor[gray]{0.9}$79.34\,[2.34]$ & $221.32$ & $9.6$ & $3.53\mathrm{e}{-29}$ & $--$ & $--$ \\
${2^{9}}$ & 5 & TT IPM & $376.01\,[1819.11]$ & $9921.33$ & $423.74\,[979.37]$ & $8171.28$ & $13$ & $7.79\mathrm{e}{-8}$ & $2.43\mathrm{e}{-10}$ & $6.62\mathrm{e}{-4}$ \\
 &  & SCS & $1133.35\,[1377.12]$ & $4239.73$ & $159.84\,[67.13]$ & $366.72$ & $7575$ & $6.99\mathrm{e}{-17}$ & $8.14\mathrm{e}{-6}$ & $6.23\mathrm{e}{-4}$ \\
 &  & MANOPT & \cellcolor[gray]{0.9}$0.127\,[0.103]$ & $0.228$ & \cellcolor[gray]{0.9}$77.10\,[6.49]$ & $221.59$ & $14.6$ & $3.38\mathrm{e}{-29}$ & $--$ & $--$ \\

\midrule
${2^{10}}$ & 1 & TT IPM & $28.18\,[16.69]$ & $123.78$ & $200.7\,[189.14]$ & $594.2$ & $11.8$ & $6.17\mathrm{e}{-8}$ & $3.26\mathrm{e}{-8}$ & $6.45\mathrm{e}{-4}$ \\
 &  & SCS & $208.65\,[184.45]$ & $289.11$ & $582.09\,[4.96]$ & $847.04$ & $315$ & $3.10\mathrm{e}{-16}$ & $8.60\mathrm{e}{-7}$ & $9.21\mathrm{e}{-3}$ \\
 &  & MANOPT & \cellcolor[gray]{0.9}$0.143\,[0.025]$ & $0.184$ & \cellcolor[gray]{0.9}$102.91\,[12.63]$ & $256.39$ & $6$ & $7.21\mathrm{e}{-29}$ & $--$ & $--$ \\
${2^{10}}$ & 2 & TT IPM & $29.14\,[37.58]$ & $180.99$ & \cellcolor[gray]{0.9}$45.61\,[458.38]$ & $607.42$ & $10.6$ & $6.63\mathrm{e}{-8}$ & $1.03\mathrm{e}{-8}$ & $7.97\mathrm{e}{-4}$ \\
 &  & SCS & $296.27\,[171.77]$ & $602.72$ & $633.85\,[24.39]$ & $847.96$ & $570$ & $4.82\mathrm{e}{-15}$ & $1.04\mathrm{e}{-4}$ & $4.86\mathrm{e}{-3}$ \\
 &  & MANOPT & \cellcolor[gray]{0.9}$0.314\,[1.59]$ & $1.98$ & $112.56\,[10.32]$ & $256.50$ & $9.4$ & $7.64\mathrm{e}{-29}$ & $--$ & $--$ \\
${2^{10}}$ & 3 & TT IPM & $390.9\,[278.34]$ & $424.87$ & $236.61\,[47.12]$ & $898.12$ & $12.8$ & $1.49\mathrm{e}{-7}$ & $4.83\mathrm{e}{-10}$ & $7.64\mathrm{e}{-4}$ \\
 &  & SCS & $383.20\,[3626.39]$ & $17464.23$ & $605.85\,[125.70]$ & $780.32$ & $7525$ & $1.31\mathrm{e}{-16}$ & $1.10\mathrm{e}{-3}$ & $1.87\mathrm{e}{-3}$ \\
 &  & MANOPT & \cellcolor[gray]{0.9}$0.203\,[0.014]$ & $0.215$ & \cellcolor[gray]{0.9}$102.94\,[6.59]$ & $260.21$ & $9.8$ & $7.08\mathrm{e}{-29}$ & $--$ & $--$ \\
${2^{10}}$ & 4 & TT IPM & $279.56\,[1335.33]$ & $2809.52$ & $763.38\,[2763.02]$ & $6669.44$ & $13.4$ & $3.17\mathrm{e}{-7}$ & $2.91\mathrm{e}{-8}$ & $1.66\mathrm{e}{-3}$ \\
 &  & MANOPT & \cellcolor[gray]{0.9}$0.178\,[0.015]$ & $0.346$ & \cellcolor[gray]{0.9}$111.20\,[1.73]$ & $231.19$ & $10.2$ & $7.26\mathrm{e}{-29}$ & $--$ & $--$ \\

\midrule
${2^{11}}$ & 1 & TT IPM & $72.12\,[53.56]$ & $90.26$ & $284.76\,[281.16]$ & $647.97$ & $10.6$ & $4.20\mathrm{e}{-8}$ & $4.39\mathrm{e}{-9}$ & $7.02\mathrm{e}{-4}$ \\
 &  & SCS & $1969.90\,[2284.00]$ & $3493.32$ & $2336.45\,[277.94]$ & $2621.32$ & $405$ & $2.77\mathrm{e}{-16}$ & $\mathbf{1.87\mathrm{e}{-2}}$ & $\mathbf{9.32\mathrm{e}{-2}}$ \\
 &  & MANOPT & \cellcolor[gray]{0.9}$0.773\,[0.231]$ & $0.901$ & \cellcolor[gray]{0.9}$231.42\,[7.20]$ & $409.02$ & $6$ & $1.96\mathrm{e}{-28}$ & $--$ & $--$ \\
${2^{11}}$ & 2 & TT IPM & $104.07\,[91.78]$ & $124.25$ & $355.27\,[452.91]$ & $658.66$ & $11.6$ & $3.37\mathrm{e}{-8}$ & $9.40\mathrm{e}{-11}$ & $4.50\mathrm{e}{-4}$ \\
 &  & MANOPT & \cellcolor[gray]{0.9}$2.35\,[2.20]$ & $4.64$ & \cellcolor[gray]{0.9}$231.39\,[10.17]$ & $427.75$ & $10$ & $1.87\mathrm{e}{-28}$ & $--$ & $--$ \\
${2^{11}}$ & 3 & TT IPM & $343.74\,[1424.27]$ & $6149.85$ & $1140.94\,[5924.15]$ & $12019.87$ & $18.2$ & $4.11\mathrm{e}{-7}$ & $7.01\mathrm{e}{-8}$ & $2.47\mathrm{e}{-3}$ \\
 &  & MANOPT & \cellcolor[gray]{0.9}$1.06\,[0.429]$ & $2.32$ & \cellcolor[gray]{0.9}$237.99\,[27.57]$ & $409.21$ & $11$ & $2.18\mathrm{e}{-28}$ & $--$ & $--$ \\

\midrule
${2^{12}}$ & 1 & TT IPM & $142.55\,[4405.55]$ & $8821.16$ & \cellcolor[gray]{0.9}$316.27\,[1308.11]$ & $2790.73$ & $12.6$ & $9.18\mathrm{e}{-8}$ & $7.52\mathrm{e}{-10}$ & $7.68\mathrm{e}{-4}$ \\
 &  & MANOPT & \cellcolor[gray]{0.9}$6.39\,[0.767]$ & $6.50$ & $751.80\,[47.44]$ & $917.99$ & $7.4$ & $3.89\mathrm{e}{-28}$ & $--$ & $--$ \\

\end{longtable}
\endgroup

\begingroup
\setlength{\LTleft}{0pt}
\setlength{\LTright}{0pt}
\fontsize{6}{7}\selectfont
\setlength{\tabcolsep}{0pt}
\renewcommand{\arraystretch}{1.08}
\begin{longtable}{@{\extracolsep{\fill}}cc l c c c c c c c c@{}}
\caption{{Performance and accuracy comparison for Maximum Stable Set instances with varying input TT rank $r_{\mathrm{in}}$. Runtime is reported in seconds and memory in MB. Runtime and memory entries show median [IQR] over seed runs; the max columns report the worst observed value. Iteration counts and accuracy entries report arithmetic means over completed seed runs. Accuracy columns report final primal feasibility error, dual feasibility error, and duality gap. For each size, TT IPM is shown only up to input rank \(\min(r_{\mathrm{stop}},5)\), where \(r_{\mathrm{stop}}\) is the first rank at which it becomes OOT. If TT IPM is reported for an instance but SDPA or SCS is absent, the missing solver was OOT/OOM. Errors which did not achieve the required target tolerance are printed in bold. }}
\label{tab:comparison_stable_set_ranks}\\
\toprule
\multicolumn{2}{c}{Inst.}
& Method
& \multicolumn{4}{c}{Perf.}
& \multicolumn{4}{c}{Acc.} \\
\cmidrule(lr){1-2}
\cmidrule(lr){4-7}
\cmidrule(lr){8-11}
Size & $r_{\mathrm{in}}$
&
& \shortstack{Time\\med.}
& \shortstack{Time\\max}
& \shortstack{Mem.\\med.}
& \shortstack{Mem.\\max}
& Iter.
& \shortstack{Primal\\err.}
& \shortstack{Dual\\err.}
& Gap \\
\midrule
\endfirsthead
\toprule
\multicolumn{11}{c}{\tablename\ \thetable{} -- cont.}\\
\midrule
Size & $r_{\mathrm{in}}$ & Method & \shortstack{Time\\med.} & \shortstack{Time\\max} & \shortstack{Mem.\\med.} & \shortstack{Mem.\\max} & Iter. & \shortstack{Primal\\err.} & \shortstack{Dual\\err.} & Gap \\
\midrule
\endhead
\midrule
\multicolumn{11}{r}{Continued}\\
\endfoot
\bottomrule
\endlastfoot
${2^{6}}$ & 1 & TT IPM & $1.79\,[0.409]$ & $2.54$ & $51.25\,[14.49]$ & $125.63$ & $7.8$ & $3.16\mathrm{e}{-8}$ & $1.25\mathrm{e}{-9}$ & $3.54\mathrm{e}{-4}$ \\
 &  & SCS & \cellcolor[gray]{0.9}$0.134\,[0.074]$ & $48.44$ & \cellcolor[gray]{0.9}$44.85\,[0.164]$ & $125.82$ & $20160$ & $1.37\mathrm{e}{-15}$ & $3.12\mathrm{e}{-6}$ & $2.75\mathrm{e}{-5}$ \\
 &  & SDPA & $1.68\,[0.103]$ & $1.76$ & $51.34\,[0.246]$ & $52.5$ & $15.2$ & $3.00\mathrm{e}{-10}$ & $1.13\mathrm{e}{-6}$ & $3.01\mathrm{e}{-14}$ \\
${2^{6}}$ & 2 & TT IPM & $4.32\,[9.71]$ & $14.51$ & $87.13\,[48.51]$ & $161.83$ & $8.6$ & $1.58\mathrm{e}{-7}$ & $3.50\mathrm{e}{-9}$ & $5.94\mathrm{e}{-4}$ \\
 &  & SCS & \cellcolor[gray]{0.9}$0.266\,[0.051]$ & $1.03$ & \cellcolor[gray]{0.9}$44.85\,[0.422]$ & $126.21$ & $595$ & $8.14\mathrm{e}{-18}$ & $9.25\mathrm{e}{-7}$ & $7.61\mathrm{e}{-6}$ \\
 &  & SDPA & $4.2\,[0.071]$ & $4.23$ & $49.44\,[0.035]$ & $49.46$ & $17$ & $1.18\mathrm{e}{-11}$ & $1.06\mathrm{e}{-6}$ & $1.71\mathrm{e}{-15}$ \\
${2^{6}}$ & 3 & TT IPM & $9.65\,[1.69]$ & $27.00$ & $102.55\,[16.81]$ & $200.81$ & $8.4$ & $1.22\mathrm{e}{-7}$ & $2.98\mathrm{e}{-9}$ & $6.27\mathrm{e}{-4}$ \\
 &  & SCS & \cellcolor[gray]{0.9}$0.273\,[0.664]$ & $2.93$ & \cellcolor[gray]{0.9}$44.59\,[0.651]$ & $126.18$ & $1155$ & $5.71\mathrm{e}{-17}$ & $8.14\mathrm{e}{-8}$ & $7.40\mathrm{e}{-7}$ \\
 &  & SDPA & $4.24\,[0.324]$ & $4.42$ & $49.6\,[0.094]$ & $49.65$ & $16.8$ & $4.10\mathrm{e}{-11}$ & $1.34\mathrm{e}{-6}$ & $1.84\mathrm{e}{-15}$ \\
${2^{6}}$ & 4 & TT IPM & $9.30\,[8.07]$ & $36.76$ & $90.81\,[7.54]$ & $207.96$ & $9.4$ & $1.41\mathrm{e}{-7}$ & $1.43\mathrm{e}{-9}$ & $4.98\mathrm{e}{-4}$ \\
 &  & SCS & \cellcolor[gray]{0.9}$0.342\,[0.092]$ & $0.533$ & \cellcolor[gray]{0.9}$44.50\,[0.828]$ & $125.26$ & $545$ & $1.42\mathrm{e}{-18}$ & $2.31\mathrm{e}{-7}$ & $1.70\mathrm{e}{-6}$ \\
 &  & SDPA & $4.13\,[0.225]$ & $4.34$ & $49.41\,[0.156]$ & $50.34$ & $16.8$ & $3.89\mathrm{e}{-11}$ & $2.63\mathrm{e}{-6}$ & $3.16\mathrm{e}{-15}$ \\
${2^{6}}$ & 5 & TT IPM & $25.48\,[21.13]$ & $59.05$ & $97.26\,[10.43]$ & $209.02$ & $10.2$ & $2.88\mathrm{e}{-7}$ & $2.79\mathrm{e}{-9}$ & $6.42\mathrm{e}{-4}$ \\
 &  & SCS & \cellcolor[gray]{0.9}$0.288\,[0.079]$ & $1.09$ & \cellcolor[gray]{0.9}$44.61\,[0.691]$ & $125.84$ & $635$ & $3.92\mathrm{e}{-18}$ & $7.15\mathrm{e}{-7}$ & $6.09\mathrm{e}{-6}$ \\
 &  & SDPA & $4.14\,[0.11]$ & $4.44$ & $49.65\,[0.117]$ & $50.52$ & $17.2$ & $1.93\mathrm{e}{-11}$ & $2.99\mathrm{e}{-6}$ & $1.91\mathrm{e}{-15}$ \\

\midrule
${2^{7}}$ & 1 & TT IPM & $2.27\,[1.15]$ & $9.09$ & \cellcolor[gray]{0.9}$43.69\,[0.066]$ & $270.63$ & $7.8$ & $5.04\mathrm{e}{-8}$ & $1.04\mathrm{e}{-9}$ & $6.43\mathrm{e}{-4}$ \\
 &  & SCS & \cellcolor[gray]{0.9}$0.619\,[0.217]$ & $0.887$ & $50.57\,[4.04]$ & $132.88$ & $230$ & $1.26\mathrm{e}{-20}$ & $4.86\mathrm{e}{-6}$ & $4.14\mathrm{e}{-5}$ \\
 &  & SDPA & $58.37\,[3.78]$ & $62.12$ & $556.31\,[0.07]$ & $557.44$ & $15.8$ & $8.74\mathrm{e}{-12}$ & $\mathbf{1.12\mathrm{e}{-2}}$ & $1.32\mathrm{e}{-13}$ \\
${2^{7}}$ & 2 & TT IPM & $22.71\,[20.56]$ & $40.71$ & $56.45\,[111.13]$ & $194.71$ & $11$ & $2.50\mathrm{e}{-7}$ & $7.70\mathrm{e}{-9}$ & $6.52\mathrm{e}{-4}$ \\
 &  & SCS & \cellcolor[gray]{0.9}$0.92\,[0.215]$ & $1.01$ & \cellcolor[gray]{0.9}$49.58\,[1.66]$ & $133.66$ & $295$ & $1.35\mathrm{e}{-18}$ & $3.57\mathrm{e}{-6}$ & $2.96\mathrm{e}{-5}$ \\
 &  & SDPA & $235.78\,[13.65]$ & $248.76$ & $551.71\,[0.684]$ & $553.64$ & $17.6$ & $7.78\mathrm{e}{-12}$ & $2.15\mathrm{e}{-6}$ & $2.82\mathrm{e}{-14}$ \\
${2^{7}}$ & 3 & TT IPM & $38.05\,[39.27]$ & $86.89$ & $131.08\,[16.17]$ & $256.56$ & $12.2$ & $3.97\mathrm{e}{-7}$ & $6.04\mathrm{e}{-9}$ & $7.39\mathrm{e}{-4}$ \\
 &  & SCS & \cellcolor[gray]{0.9}$1.10\,[0.687]$ & $3.13$ & \cellcolor[gray]{0.9}$49.78\,[0.645]$ & $133.09$ & $445$ & $2.85\mathrm{e}{-18}$ & $1.84\mathrm{e}{-6}$ & $1.49\mathrm{e}{-5}$ \\
 &  & SDPA & $235.43\,[12.35]$ & $248.17$ & $551.74\,[0.266]$ & $551.99$ & $17.8$ & $1.94\mathrm{e}{-11}$ & $2.94\mathrm{e}{-6}$ & $2.46\mathrm{e}{-14}$ \\
${2^{7}}$ & 4 & TT IPM & $90.60\,[54.35]$ & $202.63$ & \cellcolor[gray]{0.9}$47.82\,[28.75]$ & $85.51$ & $13.4$ & $4.81\mathrm{e}{-7}$ & $2.62\mathrm{e}{-9}$ & $8.23\mathrm{e}{-4}$ \\
 &  & SCS & \cellcolor[gray]{0.9}$2.12\,[19.50]$ & $415.44$ & $50.10\,[0.34]$ & $132.93$ & $21335$ & $4.88\mathrm{e}{-16}$ & $5.20\mathrm{e}{-6}$ & $4.13\mathrm{e}{-5}$ \\
 &  & SDPA & $247.22\,[10.08]$ & $261.45$ & $551.62\,[0.461]$ & $552.02$ & $18.8$ & $1.63\mathrm{e}{-10}$ & $1.57\mathrm{e}{-6}$ & $5.58\mathrm{e}{-15}$ \\
${2^{7}}$ & 5 & TT IPM & $108.35\,[37.85]$ & $1324.92$ & $98.75\,[58.80]$ & $172.07$ & $13.4$ & $6.42\mathrm{e}{-7}$ & $4.57\mathrm{e}{-9}$ & $8.26\mathrm{e}{-4}$ \\
 &  & SCS & \cellcolor[gray]{0.9}$2.91\,[1.79]$ & $483.04$ & \cellcolor[gray]{0.9}$49.84\,[0.453]$ & $134.34$ & $20460$ & $8.70\mathrm{e}{-17}$ & $4.95\mathrm{e}{-6}$ & $3.86\mathrm{e}{-5}$ \\
 &  & SDPA & $235.4\,[23.44]$ & $247.73$ & $551.73\,[0.559]$ & $552.38$ & $18$ & $9.54\mathrm{e}{-11}$ & $3.19\mathrm{e}{-6}$ & $8.54\mathrm{e}{-15}$ \\

\midrule
${2^{8}}$ & 1 & TT IPM & \cellcolor[gray]{0.9}$3.17\,[0.92]$ & $4.45$ & \cellcolor[gray]{0.9}$43.52\,[0.082]$ & $139.40$ & $9.8$ & $8.65\mathrm{e}{-8}$ & $6.09\mathrm{e}{-10}$ & $7.93\mathrm{e}{-4}$ \\
 &  & SCS & $3.28\,[2.69]$ & $2487.72$ & $66.78\,[7.19]$ & $176.93$ & $20160$ & $4.96\mathrm{e}{-16}$ & $6.09\mathrm{e}{-5}$ & $4.97\mathrm{e}{-4}$ \\
${2^{8}}$ & 2 & TT IPM & $195.31\,[288.63]$ & $503.07$ & $154.81\,[131.11]$ & $261.52$ & $13.8$ & $3.57\mathrm{e}{-7}$ & $1.09\mathrm{e}{-8}$ & $7.40\mathrm{e}{-4}$ \\
 &  & SCS & \cellcolor[gray]{0.9}$20.42\,[8.89]$ & $110.07$ & \cellcolor[gray]{0.9}$72.68\,[6.32]$ & $166.96$ & $915$ & $1.95\mathrm{e}{-17}$ & $1.15\mathrm{e}{-7}$ & $9.74\mathrm{e}{-7}$ \\
${2^{8}}$ & 3 & TT IPM & $156.93\,[298.47]$ & $893.65$ & $96.50\,[335.43]$ & $576.61$ & $15.2$ & $6.15\mathrm{e}{-7}$ & $1.40\mathrm{e}{-8}$ & $8.38\mathrm{e}{-4}$ \\
 &  & SCS & \cellcolor[gray]{0.9}$6.55\,[3.50]$ & $8.08$ & \cellcolor[gray]{0.9}$72.78\,[7.07]$ & $166.18$ & $360$ & $8.70\mathrm{e}{-18}$ & $6.30\mathrm{e}{-7}$ & $5.09\mathrm{e}{-6}$ \\
${2^{8}}$ & 4 & TT IPM & $374.26\,[1010.47]$ & $2831.30$ & $319.73\,[406.02]$ & $1460.36$ & $14.4$ & $5.20\mathrm{e}{-7}$ & $2.23\mathrm{e}{-8}$ & $8.60\mathrm{e}{-4}$ \\
 &  & SCS & \cellcolor[gray]{0.9}$10.28\,[91.72]$ & $1884.95$ & \cellcolor[gray]{0.9}$73.94\,[2.21]$ & $165.16$ & $21635$ & $5.95\mathrm{e}{-16}$ & $2.14\mathrm{e}{-5}$ & $1.99\mathrm{e}{-4}$ \\
${2^{8}}$ & 5 & TT IPM & $421.38\,[545.22]$ & $1485.40$ & $227.03\,[35.40]$ & $585.50$ & $11$ & $7.54\mathrm{e}{-7}$ & $3.57\mathrm{e}{-8}$ & $8.97\mathrm{e}{-4}$ \\
 &  & SCS & \cellcolor[gray]{0.9}$13.38\,[104.14]$ & $2468.28$ & \cellcolor[gray]{0.9}$69.08\,[3.72]$ & $166.25$ & $21380$ & $2.02\mathrm{e}{-15}$ & $3.12\mathrm{e}{-6}$ & $2.75\mathrm{e}{-5}$ \\

\midrule
${2^{9}}$ & 1 & TT IPM & \cellcolor[gray]{0.9}$6.84\,[9.81]$ & $107.88$ & \cellcolor[gray]{0.9}$56.88\,[288.96]$ & $356.95$ & $10.4$ & $7.62\mathrm{e}{-8}$ & $4.25\mathrm{e}{-8}$ & $6.67\mathrm{e}{-4}$ \\
 &  & SCS & $34.05\,[7.24]$ & $42.43$ & $178.65\,[51.76]$ & $331.03$ & $310$ & $1.12\mathrm{e}{-17}$ & $2.62\mathrm{e}{-6}$ & $2.24\mathrm{e}{-5}$ \\
${2^{9}}$ & 2 & TT IPM & $110.28\,[237.79]$ & $504.18$ & $330.70\,[147.51]$ & $610.21$ & $15$ & $1.70\mathrm{e}{-5}$ & $3.86\mathrm{e}{-8}$ & $9.28\mathrm{e}{-3}$ \\
 &  & SCS & \cellcolor[gray]{0.9}$52.77\,[242.18]$ & $316.55$ & \cellcolor[gray]{0.9}$175.36\,[21.38]$ & $295.39$ & $1330$ & $1.05\mathrm{e}{-16}$ & $4.15\mathrm{e}{-6}$ & $3.66\mathrm{e}{-5}$ \\
${2^{9}}$ & 3 & TT IPM & \cellcolor[gray]{0.9}$142.48\,[120.01]$ & $142.48$ & $518.28\,[474.78]$ & $518.28$ & $15$ & $8.13\mathrm{e}{-7}$ & $1.01\mathrm{e}{-7}$ & $2.50\mathrm{e}{-3}$ \\
 &  & SCS & $177.44\,[62.16]$ & $560.43$ & \cellcolor[gray]{0.9}$179.52\,[21.89]$ & $306.25$ & $2880$ & $1.82\mathrm{e}{-16}$ & $7.82\mathrm{e}{-7}$ & $6.12\mathrm{e}{-6}$ \\
${2^{9}}$ & 4 & TT IPM & $973.35\,[956.21]$ & $973.35$ & $1515.18\,[802.04]$ & $1515.18$ & $16$ & $1.41\mathrm{e}{-5}$ & $1.92\mathrm{e}{-7}$ & $6.08\mathrm{e}{-3}$ \\
 &  & SCS & \cellcolor[gray]{0.9}$45.41\,[4.80]$ & $404.32$ & \cellcolor[gray]{0.9}$179.50\,[4.20]$ & $304.49$ & $1665$ & $1.21\mathrm{e}{-16}$ & $2.74\mathrm{e}{-6}$ & $2.32\mathrm{e}{-5}$ \\

\midrule
${2^{10}}$ & 1 & TT IPM & \cellcolor[gray]{0.9}$23.25\,[18.88]$ & $63.37$ & \cellcolor[gray]{0.9}$49.98\,[99.93]$ & $178.59$ & $10.6$ & $8.91\mathrm{e}{-8}$ & $3.61\mathrm{e}{-8}$ & $9.05\mathrm{e}{-4}$ \\
 &  & SCS & $196.22\,[88.25]$ & $385.51$ & $694.55\,[170.08]$ & $827.41$ & $480$ & $2.96\mathrm{e}{-14}$ & $3.11\mathrm{e}{-6}$ & $2.44\mathrm{e}{-5}$ \\

\midrule
${2^{11}}$ & 1 & TT IPM & \cellcolor[gray]{0.9}$37.50\,[60.02]$ & $119.75$ & \cellcolor[gray]{0.9}$193.10\,[246.89]$ & $548.33$ & $9.4$ & $1.60\mathrm{e}{-8}$ & $6.56\mathrm{e}{-8}$ & $8.50\mathrm{e}{-3}$ \\
 &  & SCS & $1743.34\,[1925.44]$ & $3908.13$ & $2702.14\,[198.94]$ & $2761.59$ & $790$ & $3.98\mathrm{e}{-17}$ & $1.90\mathrm{e}{-5}$ & $1.52\mathrm{e}{-4}$ \\

\midrule
${2^{12}}$ & 1 & TT IPM & \cellcolor[gray]{0.9}$1070.71\,[150.25]$ & $1250.32$ & \cellcolor[gray]{0.9}$1800.55\,[80.12]$ & $1920.45$ & $13$ & $8.37\mathrm{e}{-5}$ & $7.05\mathrm{e}{-7}$ & $1.44\mathrm{e}{-3}$ \\

\end{longtable}
\endgroup

\begingroup
\setlength{\LTleft}{0pt}
\setlength{\LTright}{0pt}
\fontsize{6}{7}\selectfont
\setlength{\tabcolsep}{0pt}
\renewcommand{\arraystretch}{1.08}
\begin{longtable}{@{\extracolsep{\fill}}cc l c c c c c c c c@{}}
\caption{{Performance and accuracy comparison for Correlation Clustering instances with varying input TT rank $r_{\mathrm{in}}$. Runtime is reported in seconds and memory in MB. Runtime and memory entries show median [IQR] over seed runs; the max columns report the worst observed value. Iteration counts and accuracy entries report arithmetic means over completed seed runs. Accuracy columns report final primal feasibility error, dual feasibility error, and duality gap. For each size, TT IPM is shown only up to input rank \(\min(r_{\mathrm{stop}},5)\), where \(r_{\mathrm{stop}}\) is the first rank at which it becomes OOT. If TT IPM is reported for an instance but SDPA or SCS is absent, the missing solver was OOT/OOM. Errors which did not achieve the required target tolerance are printed in bold. }}
\label{tab:comparison_corr_clust_ranks}\\
\toprule
\multicolumn{2}{c}{Inst.}
& Method
& \multicolumn{4}{c}{Perf.}
& \multicolumn{4}{c}{Acc.} \\
\cmidrule(lr){1-2}
\cmidrule(lr){4-7}
\cmidrule(lr){8-11}
Size & $r_{\mathrm{in}}$
&
& \shortstack{Time\\med.}
& \shortstack{Time\\max}
& \shortstack{Mem.\\med.}
& \shortstack{Mem.\\max}
& Iter.
& \shortstack{Primal\\err.}
& \shortstack{Dual\\err.}
& Gap \\
\midrule
\endfirsthead
\toprule
\multicolumn{11}{c}{\tablename\ \thetable{} -- cont.}\\
\midrule
Size & $r_{\mathrm{in}}$ & Method & \shortstack{Time\\med.} & \shortstack{Time\\max} & \shortstack{Mem.\\med.} & \shortstack{Mem.\\max} & Iter. & \shortstack{Primal\\err.} & \shortstack{Dual\\err.} & Gap \\
\midrule
\endhead
\midrule
\multicolumn{11}{r}{Continued}\\
\endfoot
\bottomrule
\endlastfoot
${2^{6}}$ & 1 & TT IPM & $1.32\,[0.49]$ & $3.13$ & \cellcolor[gray]{0.9}$43.84\,[1.17]$ & $147.24$ & $9.4$ & $2.23\mathrm{e}{-9}$ & $2.46\mathrm{e}{-12}$ & $2.69\mathrm{e}{-4}$ \\
 &  & SCS & \cellcolor[gray]{0.9}$0.109\,[0.09]$ & $0.24$ & $45.14\,[0.695]$ & $125.93$ & $180$ & $6.21\mathrm{e}{-12}$ & $4.10\mathrm{e}{-12}$ & $4.68\mathrm{e}{-5}$ \\
 &  & SDPA & $4.86\,[3.49]$ & $5.28$ & $49.61\,[17.58]$ & $50.06$ & $20.8$ & $5.62\mathrm{e}{-10}$ & $3.69\mathrm{e}{-12}$ & $4.01\mathrm{e}{-6}$ \\
${2^{6}}$ & 2 & TT IPM & $2.18\,[1.46]$ & $6.37$ & $70.93\,[32.30]$ & $127.27$ & $10.4$ & $4.01\mathrm{e}{-9}$ & $1.18\mathrm{e}{-8}$ & $2.74\mathrm{e}{-4}$ \\
 &  & SCS & \cellcolor[gray]{0.9}$0.118\,[0.157]$ & $0.238$ & \cellcolor[gray]{0.9}$44.96\,[0.293]$ & $125.83$ & $220$ & $1.06\mathrm{e}{-10}$ & $2.98\mathrm{e}{-11}$ & $2.74\mathrm{e}{-5}$ \\
 &  & SDPA & $5.23\,[0.147]$ & $6.69$ & $49.77\,[0.449]$ & $50.42$ & $22.4$ & $1.52\mathrm{e}{-9}$ & $\mathbf{5.16\mathrm{e}{-2}}$ & $2.40\mathrm{e}{-6}$ \\
${2^{6}}$ & 3 & TT IPM & $6.69\,[10.90]$ & $22.64$ & $93.01\,[42.41]$ & $213.02$ & $9.4$ & $2.91\mathrm{e}{-9}$ & $8.95\mathrm{e}{-9}$ & $1.76\mathrm{e}{-4}$ \\
 &  & SCS & \cellcolor[gray]{0.9}$0.226\,[0.038]$ & $0.269$ & \cellcolor[gray]{0.9}$44.87\,[0.434]$ & $126.32$ & $265$ & $8.27\mathrm{e}{-14}$ & $8.87\mathrm{e}{-11}$ & $3.04\mathrm{e}{-5}$ \\
 &  & SDPA & $5.12\,[0.436]$ & $5.54$ & $49.65\,[0.242]$ & $53.66$ & $20.6$ & $1.07\mathrm{e}{-9}$ & $4.43\mathrm{e}{-11}$ & $1.87\mathrm{e}{-6}$ \\
${2^{6}}$ & 4 & TT IPM & $10.46\,[34.28]$ & $53.09$ & $86.82\,[49.74]$ & $164.70$ & $10.8$ & $1.86\mathrm{e}{-9}$ & $1.63\mathrm{e}{-8}$ & $1.67\mathrm{e}{-4}$ \\
 &  & SCS & \cellcolor[gray]{0.9}$0.329\,[0.157]$ & $0.535$ & \cellcolor[gray]{0.9}$44.94\,[0.676]$ & $126.03$ & $265$ & $2.02\mathrm{e}{-11}$ & $7.69\mathrm{e}{-12}$ & $2.86\mathrm{e}{-5}$ \\
 &  & SDPA & $5.19\,[0.539]$ & $5.75$ & $49.7\,[0.102]$ & $53.88$ & $21.2$ & $4.52\mathrm{e}{-10}$ & $2.68\mathrm{e}{-3}$ & $3.02\mathrm{e}{-6}$ \\
${2^{6}}$ & 5 & TT IPM & $22.90\,[24.67]$ & $50.13$ & $117.65\,[28.04]$ & $225.23$ & $12$ & $1.38\mathrm{e}{-9}$ & $2.52\mathrm{e}{-8}$ & $1.79\mathrm{e}{-4}$ \\
 &  & SCS & \cellcolor[gray]{0.9}$0.163\,[0.128]$ & $1.36$ & \cellcolor[gray]{0.9}$44.82\,[0.668]$ & $126.54$ & $565$ & $8.81\mathrm{e}{-12}$ & $6.81\mathrm{e}{-11}$ & $9.99\mathrm{e}{-6}$ \\
 &  & SDPA & $5.16\,[0.365]$ & $6.28$ & $50.02\,[0.348]$ & $50.15$ & $21.2$ & $8.14\mathrm{e}{-10}$ & $2.86\mathrm{e}{-11}$ & $2.38\mathrm{e}{-6}$ \\

\midrule
${2^{7}}$ & 1 & TT IPM & $4.46\,[3.17]$ & $9.91$ & $74.19\,[53.60]$ & $126.30$ & $12.8$ & $1.93\mathrm{e}{-9}$ & $1.59\mathrm{e}{-9}$ & $2.65\mathrm{e}{-4}$ \\
 &  & SCS & \cellcolor[gray]{0.9}$0.987\,[0.347]$ & $1.13$ & \cellcolor[gray]{0.9}$50.45\,[1.73]$ & $133.77$ & $265$ & $1.02\mathrm{e}{-10}$ & $4.63\mathrm{e}{-11}$ & $5.92\mathrm{e}{-5}$ \\
 &  & SDPA & $275.06\,[268.56]$ & $327.81$ & $552.71\,[500.28]$ & $553.83$ & $22.8$ & $3.55\mathrm{e}{-9}$ & $\mathbf{1.80\mathrm{e}{-2}}$ & $1.79\mathrm{e}{-6}$ \\
${2^{7}}$ & 2 & TT IPM & $7.15\,[5.76]$ & $11.41$ & $129.63\,[102.69]$ & $207.69$ & $12$ & $1.39\mathrm{e}{-8}$ & $1.40\mathrm{e}{-10}$ & $3.19\mathrm{e}{-4}$ \\
 &  & SCS & \cellcolor[gray]{0.9}$0.435\,[0.208]$ & $0.892$ & \cellcolor[gray]{0.9}$49.93\,[0.461]$ & $133.44$ & $160$ & $3.02\mathrm{e}{-11}$ & $1.43\mathrm{e}{-10}$ & $1.45\mathrm{e}{-5}$ \\
 &  & SDPA & $262.26\,[246.32]$ & $335.2$ & $551.8\,[279.07]$ & $552.01$ & $20.6$ & $2.32\mathrm{e}{-9}$ & $4.83\mathrm{e}{-11}$ & $1.79\mathrm{e}{-6}$ \\
${2^{7}}$ & 3 & TT IPM & $278.88\,[233.22]$ & $407.57$ & $150.08\,[162.79]$ & $311.55$ & $15.8$ & $2.23\mathrm{e}{-8}$ & $3.28\mathrm{e}{-7}$ & $4.19\mathrm{e}{-4}$ \\
 &  & SCS & \cellcolor[gray]{0.9}$0.745\,[0.861]$ & $29.60$ & \cellcolor[gray]{0.9}$50.29\,[1.06]$ & $133.72$ & $2050$ & $1.09\mathrm{e}{-11}$ & $3.27\mathrm{e}{-10}$ & $1.96\mathrm{e}{-4}$ \\
 &  & SDPA & $274.72\,[89.12]$ & $365.51$ & $552.25\,[0.512]$ & $552.72$ & $23.8$ & $2.83\mathrm{e}{-9}$ & $\mathbf{1.30\mathrm{e}{-1}}$ & $4.22\mathrm{e}{-6}$ \\
${2^{7}}$ & 4 & TT IPM & $102.65\,[11.48]$ & $241.40$ & $198.21\,[98.11]$ & $309.55$ & $11.6$ & $7.45\mathrm{e}{-9}$ & $1.17\mathrm{e}{-7}$ & $4.60\mathrm{e}{-4}$ \\
 &  & SCS & \cellcolor[gray]{0.9}$1.03\,[2.25]$ & $3.61$ & \cellcolor[gray]{0.9}$49.99\,[1.19]$ & $132.72$ & $185$ & $3.04\mathrm{e}{-12}$ & $4.19\mathrm{e}{-11}$ & $3.10\mathrm{e}{-5}$ \\
 &  & SDPA & $275.12\,[1.17]$ & $336.84$ & $551.72\,[0.32]$ & $551.97$ & $21$ & $1.20\mathrm{e}{-9}$ & $2.66\mathrm{e}{-11}$ & $4.57\mathrm{e}{-6}$ \\
${2^{7}}$ & 5 & TT IPM & $241.81\,[162.28]$ & $646.28$ & $452.39\,[490.53]$ & $2034.34$ & $17.4$ & $6.10\mathrm{e}{-9}$ & $1.37\mathrm{e}{-7}$ & $2.91\mathrm{e}{-4}$ \\
 &  & SCS & \cellcolor[gray]{0.9}$1.36\,[10.67]$ & $61.14$ & \cellcolor[gray]{0.9}$50.35\,[0.363]$ & $133.91$ & $3595$ & $2.61\mathrm{e}{-11}$ & $5.77\mathrm{e}{-10}$ & $7.69\mathrm{e}{-5}$ \\
 &  & SDPA & $327.42\,[65.07]$ & $364.48$ & $552.31\,[0.141]$ & $552.76$ & $25$ & $2.28\mathrm{e}{-9}$ & $2.27\mathrm{e}{-3}$ & $5.07\mathrm{e}{-6}$ \\

\midrule
${2^{8}}$ & 1 & TT IPM & $4.50\,[15.13]$ & $53.88$ & \cellcolor[gray]{0.9}$47.53\,[150.25]$ & $625.89$ & $14$ & $1.65\mathrm{e}{-8}$ & $1.58\mathrm{e}{-7}$ & $3.65\mathrm{e}{-4}$ \\
 &  & SCS & \cellcolor[gray]{0.9}$1.76\,[3.55]$ & $5.81$ & $65.89\,[1.45]$ & $161.79$ & $185$ & $1.37\mathrm{e}{-9}$ & $2.45\mathrm{e}{-10}$ & $6.89\mathrm{e}{-4}$ \\
${2^{8}}$ & 2 & TT IPM & $22.98\,[56.33]$ & $286.72$ & \cellcolor[gray]{0.9}$49.53\,[281.12]$ & $359.91$ & $12.4$ & $1.64\mathrm{e}{-8}$ & $8.91\mathrm{e}{-8}$ & $4.49\mathrm{e}{-4}$ \\
 &  & SCS & \cellcolor[gray]{0.9}$4.26\,[3.56]$ & $6.30$ & $75.48\,[11.34]$ & $159.98$ & $230$ & $4.28\mathrm{e}{-11}$ & $8.89\mathrm{e}{-11}$ & $9.20\mathrm{e}{-4}$ \\
${2^{8}}$ & 3 & TT IPM & $125.83\,[534.60]$ & $6395.96$ & \cellcolor[gray]{0.9}$43.70\,[273.46]$ & $3734.82$ & $15.6$ & $2.86\mathrm{e}{-8}$ & $2.23\mathrm{e}{-7}$ & $3.76\mathrm{e}{-4}$ \\
 &  & SCS & \cellcolor[gray]{0.9}$5.08\,[0.364]$ & $6.95$ & $69.12\,[8.78]$ & $165.44$ & $260$ & $3.70\mathrm{e}{-11}$ & $4.09\mathrm{e}{-10}$ & $1.19\mathrm{e}{-4}$ \\
${2^{8}}$ & 4 & TT IPM & $556.34\,[1274.30]$ & $3235.11$ & $541.90\,[2596.21]$ & $3035.55$ & $16.2$ & $1.60\mathrm{e}{-7}$ & $5.31\mathrm{e}{-7}$ & $4.94\mathrm{e}{-4}$ \\
 &  & SCS & \cellcolor[gray]{0.9}$4.92\,[7.45]$ & $28.68$ & \cellcolor[gray]{0.9}$67.43\,[3.85]$ & $168.02$ & $355$ & $1.06\mathrm{e}{-10}$ & $4.61\mathrm{e}{-10}$ & $9.22\mathrm{e}{-4}$ \\
${2^{8}}$ & 5 & TT IPM & $1168.48\,[1385.60]$ & $2754.63$ & $391.55\,[432.08]$ & $918.65$ & $17.2$ & $3.70\mathrm{e}{-7}$ & $9.58\mathrm{e}{-8}$ & $3.90\mathrm{e}{-4}$ \\
 &  & SCS & \cellcolor[gray]{0.9}$5.63\,[2.34]$ & $17.47$ & \cellcolor[gray]{0.9}$67.43\,[0.93]$ & $169.61$ & $485$ & $4.63\mathrm{e}{-11}$ & $5.50\mathrm{e}{-10}$ & $1.74\mathrm{e}{-4}$ \\

\midrule
${2^{9}}$ & 1 & TT IPM & $33.92\,[23.33]$ & $57.99$ & $161.01\,[346.67]$ & $391.45$ & $13.6$ & $7.31\mathrm{e}{-8}$ & $1.70\mathrm{e}{-8}$ & $4.92\mathrm{e}{-4}$ \\
 &  & SCS & \cellcolor[gray]{0.9}$7.62\,[2.16]$ & $9.69$ & \cellcolor[gray]{0.9}$126.45\,[41.90]$ & $290.62$ & $90$ & $5.01\mathrm{e}{-10}$ & $2.92\mathrm{e}{-12}$ & $4.25\mathrm{e}{-4}$ \\
${2^{9}}$ & 2 & TT IPM & $37.35\,[17.54]$ & $82.25$ & \cellcolor[gray]{0.9}$125.99\,[243.41]$ & $532.27$ & $16.2$ & $4.57\mathrm{e}{-8}$ & $3.80\mathrm{e}{-8}$ & $5.73\mathrm{e}{-4}$ \\
 &  & SCS & \cellcolor[gray]{0.9}$30.24\,[19.51]$ & $600.09$ & $194.03\,[58.98]$ & $331.31$ & $915$ & $2.06\mathrm{e}{-11}$ & $2.06\mathrm{e}{-8}$ & $1.07\mathrm{e}{-4}$ \\
${2^{9}}$ & 3 & TT IPM & $650.11\,[1809.04]$ & $6404.36$ & $725.20\,[2200.60]$ & $6901.30$ & $16.8$ & $5.90\mathrm{e}{-7}$ & $3.40\mathrm{e}{-7}$ & $1.42\mathrm{e}{-3}$ \\
 &  & SCS & \cellcolor[gray]{0.9}$44.10\,[35.05]$ & $165.08$ & \cellcolor[gray]{0.9}$183.60\,[13.26]$ & $297.26$ & $495$ & $1.25\mathrm{e}{-10}$ & $2.08\mathrm{e}{-9}$ & $6.30\mathrm{e}{-4}$ \\
${2^{9}}$ & 4 & TT IPM & $16472.36\,[8439.21]$ & $17863.42$ & $12151.98\,[7198.07]$ & $12151.98$ & $22$ & $7.48\mathrm{e}{-6}$ & $3.03\mathrm{e}{-6}$ & $1.75\mathrm{e}{-3}$ \\
 &  & SCS & \cellcolor[gray]{0.9}$25.33\,[157.90]$ & $620.86$ & \cellcolor[gray]{0.9}$180.71\,[31.16]$ & $306.14$ & $1120$ & $1.27\mathrm{e}{-11}$ & $4.21\mathrm{e}{-10}$ & $6.86\mathrm{e}{-5}$ \\

\midrule
${2^{10}}$ & 1 & TT IPM & \cellcolor[gray]{0.9} $32.20\,[28.28]$ & $79.95$ & \cellcolor[gray]{0.9}$128.09\,[175.55]$ & $295.54$ & $11$ & $3.99\mathrm{e}{-8}$ & $3.38\mathrm{e}{-10}$ & $8.97\mathrm{e}{-4}$ \\
 &  & SCS & $68.50\,[12.37]$ & $225.03$ & $636.55\,[38.60]$ & $728.28$ & $140$ & $6.36\mathrm{e}{-9}$ & $2.30\mathrm{e}{-10}$ & $2.10\mathrm{e}{-5}$ \\

 \midrule
${2^{11}}$ & 1 & TT IPM & \cellcolor[gray]{0.9} $5493.76\,[4715.61]$ & $11209.38$ & \cellcolor[gray]{0.9} $3699.79\,[3655.91]$ & $7355.71$ & $14.6$ & $1.57\mathrm{e}{-6}$ & $3.06\mathrm{e}{-8}$ & $4.56\mathrm{e}{-4}$ \\

\midrule
${2^{12}}$ & 1 & TT IPM & \cellcolor[gray]{0.9} $7173.93\,[6312.16]$ & $12173.93$ & \cellcolor[gray]{0.9} $3842.36\,[3144.27]$ & $7388.62$ & $14.6$ & $9.62\mathrm{e}{-7}$ & $5.09\mathrm{e}{-8}$ & $8.89\mathrm{e}{-3}$ \\

\end{longtable}
\endgroup

\end{appendices}

\end{document}